\documentclass[moor,nonblindrev]{informs1} % for review, not blinded
\usepackage{natbib}
 \NatBibNumeric
 \def\bibfont{\small}%
 \def\BIBand{and}%
 \bibpunct[, ]{[}{]}{,}{n}{}{,}%

\usepackage[colorlinks=true,breaklinks=true,bookmarks=true,urlcolor=blue,
     citecolor=blue,linkcolor=blue,bookmarksopen=false,draft=false]{hyperref}

\def\EMAIL#1{\href{mailto:#1}{#1}}% When hyperref is used, otherwise outcomment 
\def\URL#1{\href{#1}{#1}}         % When hyperref is used, otherwise outcomment 

\TheoremsNumberedThrough     % Preferred (Theorem 1, Lemma 1, Theorem 2)
\EquationsNumberedBySection % (1.1), (1.2), ...

\usepackage{setspace}
\usepackage{amsmath}
\usepackage{cleveref}
\usepackage{braket}
\usepackage{amsopn}
\usepackage{todonotes}
\usepackage{algorithm}
\usepackage{algpseudocode}
\usepackage{dcolumn,multirow,ccaption}
\usepackage{here}
\usepackage{comment}
\usepackage{cite}
\usepackage{mathtools}
\usepackage{amssymb}
\usepackage{multirow}
\usepackage{txfonts}
\usepackage{color}
\usepackage{graphicx}
\newcommand{\wtGast}{\widetilde{G}(\xast)}
\newcommand{\rast}{r_{\ast}}
\renewcommand{\O}{{\rm O}}
\renewcommand{\o}{{\rm o}}

\newcommand{\bmuone}{\bmu_1}
\newcommand{\bmutwo}{\bmu_2}
\newcommand{\bmuthree}{\bmu_3}
\newcommand{\bmufour}{\bmu_4}

\newcommand{\bmurho}{\bmu_{\rho}}
\newcommand{\brhoone}{\brho_1}
\newcommand{\brhotwo}{\brho_2}
\newcommand{\brhothree}{\brho_3}
\newcommand{\brhofour}{\brho_4}

\newcommand{\dimh}{s}
\newcommand{\vast}{v^{\ast}}

\newcommand{\hball}{\mathcal{B}}

\newcommand{\UV}{V}
\newcommand{\uv}{v}
\newcommand{\uvmu}{\uv_{\mu}}

\newcommand{\wtGell}{\widetilde{G}_{\ell}}

\newcommand{\dimeq}{s}
\newcommand{\calC}{\mathcal{C}}

\newcommand{\impfun}{\Phi}

\newcommand{\wt}{\widetilde}

\newcommand{\fr}{\frac{1}{2}}

\newcommand{\Gi}{\mathcal{G}_i}

\newcommand{\K}{\mathcal{K}}            
\newcommand{\R}{\mathbb{R}}

\newcommand{\wl}{w^{\ell}}
\newcommand{\xl}{x^{\ell}}

\newcommand{\W}{\mathcal{W}}

\newcommand{\Pl}{P_{\ell}}

\renewcommand{\epsilon}{\varepsilon}

\newcommand{\wk}{w^k}

\renewcommand{\phi}{\varphi}

\newcommand{\G}{\mathcal{G}}

\newcommand{\bd}{{\mathop{\mathrm{bd}\,}}}
\newcommand{\cl}{{\mathop{\mathrm{cl}\,}}}

\newcommand{\dist}{{\mathop{\mathrm{dist}\,}}}

\def\ru[#1][#2]{{#1}^{#2}}
\def\rl[#1][#2]{{#1}_{#2}}

\newcommand{\tred}{}

\newcommand{\cred}{}

\newcommand{\xiast}{\xi^{\ast}}

\newcommand{\calL}{\mathcal{L}}
\newcommand{\xmu}{x(\mu)}
\newcommand{\Uell}{\UV_{\ell}}
\newcommand{\bK}{\overline{K}}
\newcommand{\tuell}{\theta^{\ell}}
\newcommand{\txell}{y^{\ell}}

\newcommand{\muk}{\mu_k}

\newcommand{\calG}{\mathcal{G}}

\newcommand{\bmu}{\bar{\mu}}

\newcommand{\wmu}{w(\mu)}

\newcommand{\umuell}{\uv_{\mul}}
\newcommand{\xmudell}{\chx(\mul)}
\newcommand{\bxmudell}{\chx_{\calM}(\mul)}
\newcommand{\chx}{\check{x}}
\newcommand{\xmudel}{\chx(\mu)}

\newcommand{\ZAA}{Z_{11}}
\newcommand{\ZAB}{Z_{12}}

\newcommand{\ZBB}{Z_{22}}
\newcommand{\SSOSC}{ESOSC}
\newcommand{\DeltawtGell}{\wt{\Delta G}_{\ell}}
\newcommand{\DeltaGell}{\Delta\Gell}

\newcommand{\DGxzero}{D_{\Gell}^0}
\newcommand{\DGxpp}{D_{\Gell}^{++}}

\newcommand{\Yrho}{\Yast}
\newcommand{\zrho}{\zast}
\newcommand{\wrho}{\wast}
\newcommand{\brho}{\bar{\rho}}

\newcommand{\calM}{\mathcal{M}}
\newcommand{\Sc}{S_{\rm c}}
\newcommand{\wtYell}{\wt{Y}_{\ell}}
\newcommand{\uast}{\uv_{\ast}}
\newcommand{\Scxu}{\Sc(x,\mu)}

\newcommand{\VAA}{Z_{11}}
\newcommand{\WAA}{W_{11}}
\newcommand{\mul}{\mu_{\ell}}

\newcommand{\VAB}{Z_{12}}
\newcommand{\WAB}{W_{12}}

\newcommand{\VBB}{Z_{22}}

\newcommand{\WBB}{W_{22}}

\newcommand{\calHl}{\mathcal{H}_{\ell}}

\newcommand{\calA}{\mathcal{A}}
\newcommand{\xell}{x^{\ell}}
\newcommand{\NN}{\rm FF}
\newcommand{\Un}{\rm EF}
\newcommand{\NU}{\rm FE}
\newcommand{\UU}{\rm EE}
\newcommand{\GNN}{G^{\NN}}
\newcommand{\GUU}{G^{\UU}}
\newcommand{\GUN}{G^{\Un}}
\newcommand{\GNU}{G^{\NU}}
\newcommand{\YNN}{Y^{\NN}}
\newcommand{\YUU}{Y^{\UU}}
\newcommand{\YUN}{Y^{\Un}}
\newcommand{\YNU}{Y^{\NU}}

\newcommand{\etal}{\eta^{\ell}}
\newcommand{\dk}{d^k}
\newcommand{\UR}{U_{\xast}}

\newcommand{\UN}{V_{\ast}}
\newcommand{\xk}{x^k}
\newcommand{\Pell}{P_{\ell}}
\newcommand{\calPc}{\mathcal{P}_{\rho}}

\newcommand{\calPcbone}{\mathcal{P}_{\brhoone}}

\newcommand{\Gk}{G_k}
\newcommand{\Yk}{Y_k}
\newcommand{\zk}{z^k}
\newcommand{\aell}{a_{\ell}}
\newcommand{\bell}{b_{\ell}}
\newcommand{\cell}{c_{\ell}}

\newcommand{\Gell}{G_{\ell}}
\newcommand{\well}{w^{\ell}}
\newcommand{\yell}{y^{\ell}}
\newcommand{\Yell}{Y_{\ell}}
\newcommand{\zell}{z^{\ell}}

\newcommand{\uell}{\uv^{\ell}}
\newcommand{\calJ}{\mathcal{J}}
\newcommand{\vell}{v^{\ell}}
\newcommand{\JPHI}{\mathcal{A}}

\newcommand{\YkUN}{\YUN_k}
\newcommand{\YkNU}{\YNU_k}
\newcommand{\YkNN}{\YNN_k}
\newcommand{\YaUU}{\YUU_{\rm a}}

\newcommand{\JGast}{\mathcal{J}G(\xast)}
\newcommand{\Gast}{G_{\ast}}
\newcommand{\Past}{P_{\ast}}
\newcommand{\past}{p_{\ast}}
\newcommand{\mathUast}{\mathcal{U}_{\ast}}
\newcommand{\hYa}{\widetilde{\Ya}}
\newcommand{\hza}{\widetilde{\za}}
\newcommand{\balpha}{\bar{\alpha}}

\newcommand{\Tr}{{\rm Tr}}

\newcommand{\olY}{\overline{Y}}
\newcommand{\DeltaG}{\Delta G}
\newcommand{\PUN}{P_{\ast}}
\newcommand{\PU}{E_{\ast}}
\newcommand{\PN}{F_{\ast}}
\newcommand{\PUell}{E_{\ell}}
\newcommand{\PNell}{F_{\ell}}
\newcommand{\DeltahGUU}{\Delta {G}^{\UU}}
\newcommand{\DeltahGUN}{\Delta {G}^{\Un}}
\newcommand{\DeltahGNU}{\Delta {G}^{\NU}}
\newcommand{\DeltahGNN}{\Delta {G}^{\NN}}

\newcommand{\DeltaGUU}{\Delta G^{\UU}}
\newcommand{\DeltaGUN}{\Delta G^{\Un}}
\newcommand{\DeltaGNU}{\Delta G^{\NU}}
\newcommand{\DeltaGNN}{\Delta G^{\NN}}

\newcommand{\DeltaYUU}{Y^{\UU}}
\newcommand{\DeltaYUN}{Y^{\Un}}
\newcommand{\DeltaYNU}{Y^{\NU}}
\newcommand{\DeltaYNN}{Y^{\NN}}

\newcommand{\Ya}{Y_{\rm a}}
\newcommand{\za}{z^{\rm a}}
\newcommand{\wa}{w^{\rm a}}

\newcommand{\hd}{\widetilde{d}}
\newcommand{\Yast}{Y_{\ast}}
\newcommand{\zast}{z^{\ast}}
\newcommand{\wast}{w^\ast}
\newcommand{\xast}{x^\ast}

\newcommand{\Ker}{\mathop{\rm Ker}}

\makeatletter
\crefname{theorem}{Theorem}{Theorems}
\crefname{proposition}{Proposition}{Propositions}
\crefname{definition}{Definition}{Definitions}
\crefname{lemma}{Lemma}{Lemmas}
\crefname{corollary}{Corollary}{Corollaries}
\crefname{assumption}{Assumption}{Assumptions}
\begin{document}
%%%%%%%%%%%%%%%%

%\setcounter{page}{1}
%\DOI{10.1287/xxxx.0000.0000}%


%\RUNTITLE{}

% \TITLE{Bundling Information Goods of Decreasing Value}
\TITLE{Analysis of the primal-dual central path for nonlinear semidefinite optimization without the nondegeneracy condition}

\ARTICLEAUTHORS{%
\AUTHOR{Takayuki Okuno}
\AFF{
Faculty of Science and Technology Department of Science and Technology, Seikei University, 3-3-1, Kita, Kichijoji, Musashino, Tokyo, 180-8633, Japan
\\
Center for Advanced Intelligence Project, RIKEN, Nihonbashi 1-chome Mitsui Building, 15th floor, 1-4-1 Nihonbashi, Chuo-ku, Tokyo 103-0027, Japan, \EMAIL{takayuki-okuno@st.seikei.ac.jp}, \URL{}}
} % end of the block

\ABSTRACT{%
We study properties of the central path underlying a nonlinear semidefinite optimization problem, called NSDP for short.
The latest radical work on this topic was contributed by Yamashita and Yabe (2012):
they proved that the Jacobian of a certain equation-system derived from the Karush-Kuhn-Tucker (KKT) conditions of the NSDP is nonsingular at a KKT point under the second-order sufficient condition (SOSC), the strict complementarity condition (SC), and the nondegeneracy condition (NC). This yields uniqueness and existence of the central path 
through the implicit function theorem. 
In this paper, we consider the following three assumptions on a KKT point:      
the \tred{enhanced} SOSC, the SC, and the Mangasarian-Fromovitz constraint qualification. 
Under the absence of the NC, the Lagrange multiplier set is not necessarily a singleton and 
the nonsingularity of the above-mentioned Jacobian is no longer valid.  Nonetheless, we establish that the central path exists uniquely, and moreover prove that the dual component of the path converges to the so-called analytic center of the Lagrange multiplier set. 
As another notable result, we clarify a region around the central path where Newton's equations relevant to primal-dual interior point methods are uniquely solvable. }%

%\KEYWORDS{deterministic inventory theory; infinite linear programming duality; 
%  existence of optimal policies; semi-Markov decision process; cyclic schedule}
%\MSCCLASS{Primary: 90B05; secondary: 90C40, 90C90}
%\ORMSCLASS{Primary: Inventory/production: deterministic multi-item;
%  secondary: dynamic programming/optimal control: deterministic 
%  semi-Markov; programming: infinite dimensional}
%\HISTORY{Received November 20, 2003; revised March 8, 2004, and March 26, 2004.}

\KEYWORDS{
nonlinear semidefinite optimization, primal-dual interior-point method, central path, nondegeneracy condition
}
%\MSCCLASS{}
%\ORMSCLASS{Primary: ; secondary: }

\maketitle
%%%%%%%%%%%%%%%%%%%%%%%%%%%%%%%%%%%%%%%%%%%%%%%%%%%%%%%%%%%%%%%%%%%%%%

% Samples of sectioning (and labeling) in MOOR.
% NOTE: (1) all section levels end with a period,
%       (2) capitalization is as shonw (sentence style, not title style).
%
%\section{Introduction.}\label{intro} %%1.
%\subsection{Duality and the classical EOQ problem.}\label{class-EOQ} %% 1.1.
%\subsection{Outline.}\label{outline1} %% 1.2.
%\subsubsection{Cyclic schedules for the general deterministic SMDP.}
%  \label{cyclic-schedules} %% 1.2.1
%\section{Problem description.}\label{problemdescription} %% 2.

% Text of your paper here

\section{Introduction}\label{sec:1}
We consider the following nonlinear semidefinite optimization problem:
\begin{align}\label{al:nsdp}
\begin{array}{rcl}
\displaystyle\mathop{\rm Minimize}& &f(x)\\
                              \mbox{subject~to}& &G(x)\in \mathbb{S}^m_+,\\
                                               & &h(x) = 0,%  \\
                                                %& &x\in X:=[a,b],
\end{array}
\end{align}
where
$f:\R^n\to \R$, $G:\R^n\to \mathbb{S}^m$, and $h:\R^n\to \R^{\dimh}$
are twice continuously differentiable functions. Moreover, $\mathbb{S}^m$ denotes the set of real 
$m\times m$ symmetric matrices
and $\mathbb{S}^m_{++}$ (resp. $\mathbb{S}^m_+$) stands for the set of 
$m\times m$ real symmetric positive definite (resp. semidefinite) matrices.
Throughout the paper, we
often refer to problem\,\eqref{al:nsdp} as NSDP.
%and suppose that $f, G,$ and $h$ are three-times continuously differentiable. 
%$f:\R^n\to \R$, $G:\R^n\to \mathbb{S}^m$, and $h:\R^n\to \R^{\ell}$ are three-times continuously differentiable functions.
NSDP\,\eqref{al:nsdp} contains a wide class of optimization problems.
Indeed, when all the functions are affine with respect to $x$, it reduces to a linear semidefinite optimization problem\,(\citet{vandenberghe1996semidefinite,wolkowicz2012handbook}). 
When the function $G$ is of the diagonal matrix form, it is regarded as a conventional nonlinear optimization problem\,(\citet{
mangasarian1994nonlinear,luen2008}).
%When $G$ takes arrow-matrix
%By letting $G$ 
Moreover, it contains nonlinear second-order cone optimization problems\,(\citet{kato2007sqp,bonnans2005perturbation})
by restricting the form of $G$ appropriately.  
%LSDPs have been extensively
%The LSDP has attracted many researchers, as it can be very powerful tools in various fields including combinatorial optimization. 
%See a comprehensive survey on LSDPs, for example, by \citet{vandenberghe1996semidefinite} and also refer to the handbook by \citet{wolkowicz2012handbook}. 
%In comparison with LSDPs, studies of NSDPs are fers 
%Compared with the LSDP, the studies of the NSDP are not sufficient 
%in the last decade,
%On the other hand, studies on the NSDP itself have also advanced significantly in the 2000s.

The recent advance of researches on the NSDP is remarkable. 
%from aspects of 
%applications in the real world and theories. 
Abundant practical applications of the NSDP can be found in a wide variety of fields, for example,
structural optimization\,(\citet{kovcvara2004solving,thore2017general,takezawa2011topology,thore2022worst}), 
control\,(\citet{scherer1995multiobjective,kovcvara2005nonlinear,hoi2003nonlinear,leibfritz2006reduced}),
\tred{statistics}\,(\citet{qi2006quadratically}), finance\,(\citet{konno2003cutting,leibfritz2009successive}),
positive semidefinite factorization\,(\citet{vandaele2018algorithms}), and so on.
%  Logit model   konno2003cutting 
% Finance leibfritz2009successive
% Controle {H}2/{H}$\infty$control scherer1995multiobjective 
% static output feedback problems  kovcvara2005nonlinear hoi2003nonlinear leibfritz2006reduced
% nearest correlation matrix qi2006quadratically  
%although they are still fewer than LSDPs. 
Elegant theoretical results on optimality conditions for the NSDP have been also developed. For example, the Karush-Kuhn-Tucker (KKT) conditions and the second-order conditions for the NSDP were studied in detail by \citet{shapiro1997first} and \citet{forsgren2000optimality}. 
Further examples are: the strong second-order conditions by \citet{sun2006strong}, sequential optimality conditions by \citet{andreani2018optimality}, the local duality by \citet{qi2009local}, and the optimality conditions via squared slack variables by \citet{lourencco2018optimality}.
Along with such theoretical results, various algorithms have been proposed for solving the NSDP, for example, augmented Lagrangian methods\,(\citet{kovcvara2004solving,sun2006properties,sun2008rate,andreani2018optimality,andreani2021use,fukuda2018exact,huang2006approximate,wu2013global}),
sequential linear semidefinite optimization methods\,(\citet{kanzow2005successive}), sequential quadratic semidefinite optimization methods\,(\citet{correa2004global,freund2007nonlinear,zhao2016superlinear,zhao2018sqp,yamakawa2019stabilized}),
sequential quadratically constrained quadratic semidefinite optimization methods\,(\citet{auslender2013extended}),
exact penalty methods\,(\citet{auslender2015exact}), interior point-type methods\,(\citet{arahata2021interior, jarre2000interior,kato2015interior,leibfritz2002interior,okuno2018primal,okuno2020interior,okuno2020local,yamashita2012local,yamashita2012primal,yamashita2021primal, yamakawa1,yamakawa2}), homotopy methods\,(\citet{yang2013homotopy}), and so forth. 

In this paper, we study properties of the central path for the NSDP. 
The central path is a path formed by stationary points of the log-barrier penalized problem, and is a key concept of interior-point methods, abbreviated as IPMs, in solving a wide class of optimization problems including the NSDP.
Many IPMs share the strategy of approaching a KKT point by following the central path approximately. 
Since the geometry of the central path is related to the performance of IPMs,  
%Behavior and curvature of the central path has been studied extensively under various problem setting.  
%Since performance of IPMs is closely related to properties of the central path, 
it has been well studied under various settings. 
For example, \citet{megiddo1989pathways} presented an early work in this line for linear optimization or linear programming.
\citet{kojima1990limiting} and \citet{monteiro1996limiting} studied the central path for monotone complementarity problems under the absence of strict complementarity condition. 
\citet{monteiro1998existence} worked with the existence of the central path for convex optimization problems.  
\citet{wright2002properties} considered nonlinear optimization problems and analyzed the properties of the central path under the absence of linear independence constraint qualification.  
%In fact, the theory in this paper is inspired by Wright and Orban's work.
%These studies are motivated from primal-dual IPMs (PDIPMs) for SDPs as limiting behavior of the central path is closely related to 
%the performance of the PDIPM.  
%Many researches of the central path of LSDPs are found. 
%As regards a class of semidefinite optimization problems (SDPs), 

We briefly review the history of the central path of semidefinite optimization problems (SDPs). 
Concerning linear SDPs, \citet{luo1998superlinear} showed that the (primal-dual) central path converges to the analytic center under the presence of the strict complementarity condition.  \citet{sturm2001sensitivity} further proved that the derivative of the central path is convergent.
%The same result was obtained by \citet{goldfarb1998interior} without assuming this strict complementarity condition.
\citet{halicka2002convergence} proved that the central path is convergent 
regardless of the strict complementarity, by means of the curve selection lemma from algebraic geometry, although it can fail to converge to the analytic center in the absence of the strict complementarity. 
%In \citep{halicka2005limiting}, the same authors further clarified that 
%the limit point is the analytic center of a certain subset of the optimal set.
\citet{halicka2002analyticity} established that the central path is analytic including the boundary.
%\citet{kakihara2014curvature} related the curvature integral along the central path to the iteration complexity of the Mizuno-Todd-Ye predictor-corrector IPM for the LSDP.  
See also other works by \citet{goldfarb1998interior,halicka2005limiting,kakihara2013information,kakihara2014curvature,da2005asymptotic}, and so forth.
More generally,  
\citet{grana2002central} worked with the existence and convergence of the central path
of convex smooth SDP by assuming that the functions organizing the problem are analytic. While there are many such studies concerning linear and convex SDPs, those for the general NSDP\,\eqref{al:nsdp} are very scarce. 

The latest radical work for NSDP\,\eqref{al:nsdp} along this research-topic was presented by \citet{yamashita2012local}. 
%In \citep{yamashita2012primal}, 
The authors analyzed the local convergence property of the primal-dual IPM, called PDIPM for short, 
that was proposed in another article of theirs\,(\citet{yamashita2012primal}). This PDIPM is explained briefly as follows: in the algorithm, the barrier KKT (BKKT) conditions are derived by perturbing the KKT conditions, and 
the degree of perturbation is controlled by the so-called barrier parameter. See Section\,\ref{subsec:bkkt} for the precise definition of the BKKT conditions.  The PDIPM approaches a KKT point by generating a sequence of approximate BKKT points while driving the barrier parameter to zero.  To compute a BKKT point, the Newton method combined with scaling techniques is applied to an equation-system equivalent to the BKKT conditions. 
In \citep{yamashita2012local}, Yamashita and Yabe proved that the Jacobian of this equation-system is nonsingular at a KKT point 
under the following three conditions: the strict complementarity condition (SC), 
the second-order sufficient condition (SOSC), and the nondegeneracy condition (NC).
Along with the classical implicit function theorem, this fact yields that 
there exists a unique smooth path, i.e., a central path, passing through the focused KKT point, and this path is formed by BKKT points.
%Given those three conditions at the KKT point, the set of corresponding Lagrange multiplier matrix and vectors is a singleton.  
\subsection*{Contribution}
The main contribution of this paper is summarized as follows:\vspace{0.5em}
\begin{enumerate}
\item We prove that there exists a smooth central path under the SC, the \tred{enhanced} SOSC, and the Mangasarian-Fromovitz constraint qualification (MFCQ) at a KKT point of the NSDP.
  We also prove that the central path converges to the KKT point and the analytic center of the corresponding Lagrange multiplier set.  Since the NC is not assumed therein, the Lagrange multiplier set is compact and convex, but not necessarily a singleton, although
the KKT point is a strict local optimum due to the \tred{enhanced} SOSC. In such a situation, 
it is difficult (or impossible) to prove existence of the central path straightforwardly by means of the implicit function theorem.
\item Under the same conditions as above, we give a region around the central path where the Newton equation is solvable uniquely when applying the PDIPM.     
\end{enumerate}\vspace{0.5em}
%One may ask about novelty of the analysis of the present paper, compared to the existing ones.
Many of the analyses in literature on SDPs exploit the fact that 
the functions are analytic and thus so is the underlying central path. 
However, this methodology is no longer available in our setting since the functions of the NSDP are not assumed to be analytic.
The manner of our analysis conducted in this paper is motivated from \citet{wright2002properties} for nonlinear optimization, but ours is more complicated because the SOSC of the NSDP involves difficulty arising from the so-called sigma term. Furthermore, we deal with the nonlinear equality constraints together, whereas \citep{wright2002properties} does not.

\subsection*{Notations and terminologies}
%\noindent{\bf Notations}\\
Throughout the paper, we use the following notations as necessary:
\tred{for a set $S$, we denote by ${\rm int}\,S$, $\cl S$, and $\bd S$ the topological interior, closure, and boundary of $S$, respectively.}
We denote the identity matrix in $\R^{m\times m}$ by $I$, and for $A\in \R^{m\times m}$, we define ${\rm Sym}(A):=(A+A^{\top})/2$ and $\|A\|_{\rm F}:=\sqrt{{\rm trace}(A^{\top}A)}$.
\tred{
For $B\in \R^{m\times n}$, we denote the kernel and image spaces of $B$ by ${\rm Ker}\,B$ and ${\rm Im}\,B$, respectively, that is, ${\rm Ker}\,B:=\{x\in \R^n\mid Bx=0\}$ and ${\rm Im}\,B:=\{By\mid y\in \R^n\}$.}
For $X,Y\in \mathbb{S}^m$, we define the inner product $X\bullet Y$ by $X\bullet Y:= {\rm trace}(XY)$. We also define the linear operator $\mathcal{L}_X:\mathbb{S}^m\to \mathbb{S}^m$ by 
$$\mathcal{L}_X(Y):=XY+YX.$$ 
Denote the smallest eigenvalue of $X\in \mathbb{S}^m$ by $\lambda_{\rm min}(X)$.
For $X\in \mathbb{S}^m_{+}$ and $r>0$, we denote by $X^{\frac{1}{r}}$ the unique solution $U\in \mathbb{S}^m_+$ of $U^r=X$. 
For a function $g:\R^n\to \R$, we denote by $\nabla g(x)$ or $\nabla_x g(x)$ the gradient of $g$, namely, $\nabla g(x):=(\frac{\partial g(x)}{\partial x_1},\ldots,
\frac{\partial g(x)}{\partial x_n}
)^{\top}\in \R^n$ and, also denote by $\nabla^2_{xx}g(x)$ the hessian of $g$, namely, $\nabla^2_{xx}g(x)=(\frac{\partial^2 g(x)}{\partial x_i\partial x_j})_{1\le i,j\le n}\in \R^{n\times n}$.
For $\{A_{k}\}$ in a \tred{normed} vector space with norm $\|\cdot\|$ and $\{b_k\}\subseteq \R$, 
we write $A_k=\O(b_k)$ if there exists some $M>0$ such that 
$\|A_k\|\le M|b_k|$ for all $k$ sufficiently large, and write $A_k=\o(b_k)$ if there exists some \tred{nonnegative} sequence $\{\alpha_k\}\subseteq \R$ such that 
$\lim_{k\to \infty}\alpha_k=0$ and 
$\|A_k\|\le \alpha_k|b_k|$ for all $k$ sufficiently large. We also say $A_k=\Theta(b_k)$ if there exist $M_1,M_2>0$ such that
$M_1|b_k|\le \|A_k\|\le M_2|b_k|$ for all $k$ sufficiently large.  

%which is often abbreviated as $\Gi$ for simplicity.
We also denote $\R_{++}:=\{a\in \R\mid a>0\}$, $\mathcal{W}:=\R^n\times \mathbb{S}^m\times \R^{\dimh}$,
%\begin{equation*}
\begin{align*}
  %\mathcal{W}_{+({\rm resp.,}\ ++)}:=\R^n\times \mathbb{S}^m_{+({\rm resp.,}\ ++)}\times \R^{\ell}.
  \mathcal{W}_{++}:=\{
(x,Y,z)\in \mathcal{W}\mid G(x)\in \mathbb{S}^m_{++}, Y\in \mathbb{S}^m_{++}
  \},\ \W_+:=\{(x,Y,z)\in \mathcal{W}\mid G(x)\in \mathbb{S}^m_{+}, Y\in \mathbb{S}^m_{+}\}.
\end{align*}
%Additionally, let $\W_+$ be the set obtained by replacing $\mathbb{S}^m_{++}$ with $\mathbb{S}^m_+$ in $\W_{++}$.
For $w:=(x,Y,z)\in \W$, we define $\|w\|:=\sqrt{\|x\|_2^2+\|Y\|_{\rm F}^2+\|z\|_2^2}$, where $\|\cdot\|_2$ denotes the Euclidean norm. 

Lastly, relevant to the function $G$ in NSDP\,\eqref{al:nsdp}, we define the following notations.
For $i=1,2,\ldots,n$, we write
$$\Gi(x):=\frac{\partial G(x)}{\partial x_i}.$$
For any $x,d\in \R^n$ and $Y\in \mathbb{S}^m$, we write 
$$
\DeltaG(x;d) :=\sum_{i=1}^nd_i\Gi(x)\in \mathbb{S}^m,\ \tred{\mathcal{J}G(x)^{\ast}Y:=\left[\mathcal{G}_1(x)\bullet Y,\mathcal{G}_2(x)\bullet Y,\ldots,
\mathcal{G}_n(x)\bullet Y
\right]^{\top}\in \R^n.}
$$
Some more notations and symbols 
will be introduced for the main analysis. 
See the paragraph {\it Additional notations and symbols used hereafter} at the end of subsection\,\ref{sec:3-1}.

\subsection*{Organization of the paper}
The rest of the paper is organized as follows.
In \cref{sec:2}, we review some important concepts related to the NSDP such as the KKT conditions.
In \cref{sec:mainresults}, the main analysis is presented.
In \cref{sec:con}, we conclude this paper with some remarks. 

\section{Preliminaries}\label{sec:2}
\subsection{KKT conditions for NSDP}
We introduce the KKT conditions for NSDP\,\eqref{al:nsdp}.
\begin{definition}
We say that the the Karush-Kuhn-Tucker (KKT) conditions for NSDP\,\eqref{al:nsdp}
hold at $x\in \R^n$ if there exist a Lagrange multiplier matrix $Y\in \mathbb{S}^m$ and vector $z\in \R^{\dimh}$ such that 
\begin{eqnarray}
&\nabla_xL(w)=\nabla f(x) - \mathcal{J}G(x)^{\ast}Y+\nabla h(x)z=0,\label{al:kkt1}\\
&G(x)\bullet Y = 0,\ G(x)\in \mathbb{S}^m_+,\ Y\in \mathbb{S}^m_+,\label{al:comp}\\
&h(x) = 0,\label{al:kkt3}
\end{eqnarray}
where
$w:=(x,Y,z)\in \mathcal{W}$ and $L:\mathcal{W}\to \R$ denotes the Lagrange function for the NSDP, that is, 
\begin{equation*}
L(w):=f(x) - G(x)\bullet Y + h(x)^{\top}z \label{eq:lo}
\end{equation*}
for any $w\in\mathcal{W}$.
Particularly, we call a triplet $w=(x,Y,z)$ satisfying the KKT conditions a KKT triplet of NSDP\,\eqref{al:nsdp}, and also call $x$ a KKT point of the NSDP.
Moreover, given a KKT point $x$, we denote by $\Lambda(x)$ 
the set of Lagrange multiplier pairs $(Y,z)$ satisfying the KKT conditions at $x$, namely,
$$
\Lambda(x):=\{(Y,z)\mbox{ satisfying \eqref{al:kkt1}-\eqref{al:kkt3}}\}.
$$
\end{definition}
\tred{As can be checked easily,
$\Lambda(x)$ is convex.}
Below, we define the Mangasarian-Fromovitz constraint qualification (MFCQ), under which 
the KKT conditions are ensured to be necessary optimality conditions for the NSDP. 
\begin{definition}(\citep[Definition~2.8.6]{bonnans2013perturbation})
\label{def:mfcq}
Let ${x}\in \R^n$ be a feasible point of NSDP\,\eqref{al:nsdp}.
We say that
the Mangasarian-Fromovitz constraint qualification (MFCQ) holds at ${x}$
if 
$\nabla h({x})$ is of full column rank and 
there exists a vector $d\in \R^n$ such that
$G({x})+\DeltaG({x};d)\in \mathbb{S}^m_{++}$ and $\nabla h(x)^{\top}d=0$.
\end{definition}
\begin{remark}
      {The MFCQ is equivalent to the following Robinson's constraint qualification at a feasible point $x\in \R^n$ \citep[Corollary~2.101]{bonnans2013perturbation}:
  $$
  \begin{bmatrix}O\\
    0
  \end{bmatrix}\in {\rm int}\left(\Set{\begin{bmatrix}
  G(x)+\DeltaG(x;d)\\
     \nabla h(x)^{\top}d
\end{bmatrix}
    | d\in \R^n}-\begin{bmatrix}\mathbb{S}^m_{++}\\ \{0\}\end{bmatrix}\right).$$}
\end{remark}
\begin{remark}
Let ${x}\in \R^n$ be a local optimum of NSDP\,\eqref{al:nsdp}.
Under the MFCQ, the KKT conditions hold at ${x}$, thus $\Lambda(x)\neq \emptyset$. 
In particular, the MFCQ implies that $\Lambda(x)$ is \tred{compact}.
Conversely, when $f$ is convex, $h$ is affine, and $G$ is matrix-convex in the sense of \citet[Section~5.3.2]{bonnans2013perturbation}, a KKT point is a global optimum of \eqref{al:nsdp}.
\end{remark}
%The relation \eqref{al:comp} represents the complementary condition associated to the semi-definite constraint $G(x)\in \mathbb{S}^m_+$. As known well,  it yields the crucial property that $G(x)$ and $Y$ commute, namely, $G(x)Y=YG(x)$.
%According to this fact, we can obtain several conditions equivalent to \eqref{al:comp}.
%Note that  if and only if $G(x)\bullet Y=0$. 
%Since 
%$G(x)Y=O$ and $G(x)$ and $Y$ commute when $G(x)\bullet Y=0$,

There are several equivalent reformulations for the semidefinite complementarity condition\,\eqref{al:comp}, among which the simplest one is 
\begin{equation}
  G(x)Y = O,\ G(x)\in \mathbb{S}^m_{+},\ Y\in \mathbb{S}^m_{+},\label{al:xy}
\end{equation}
and two other formulations are
\begin{align}
  %\frac{G(x)Y+YG(x)}{2}=O,\ G(x)\in \mathbb{S}^m_+,\ Y\in \mathbb{S}^m_+.\label{eq:aho}
  {\rm Sym}\left(G(x)Y\right)=O,\ G(x)\in \mathbb{S}^m_+,\ Y\in \mathbb{S}^m_+,\label{eq:aho}\\
  G(x)^{\frac{1}{2}}YG(x)^{\frac{1}{2}}=O,\ G(x)\in \mathbb{S}^m_+,\ Y\in \mathbb{S}^m_+.\label{eq:mt}
\end{align}
Based on the above two formulations, 
primal-dual interior point methods (PDIPMs) have been developed for solving NSDPs so far. 
For example,
see \citet{yamashita2012primal} and \citet{yamashita2012local} for PDIPM with \eqref{eq:aho} and also see \citet{okuno2020local} for that with \eqref{eq:mt}.
%This was originally presented by Alizadeh, Haeberly, and Overton\,\citep{alizadeh1998primal} for LSDPs.
%The Monteiro-Zhang (MZ) family 
%is the set of search-directions obtained by 
%solving certain scaled Newton equations associated with condition\,\eqref{eq:aho}. See \eqref{eq:mz}.
%On the other hand,  the MT family was
%For solving LSDPs, many variants of primal-dual interior point methods (called PDIPMS for short) based on the MZ family have been developed so far. 
%In the context of the NSDP, 
%we believe that Yamashita et al.\,

\paragraph{Other fundamental properties of the complementarity condition}
Let $\xast\in \R^n$ be a KKT point for the NSDP. 
With an appropriate orthogonal matrix $P_{\ast}\in \R^{m\times m}$, 
the matrix $G(x^{\ast})$ and an arbitrary dual matrix
\tred{$\tred{Y}\in \mathbb{S}^m_+$ such that $\Gast\tred{Y}=O$ holds can be factorized as}
\begin{equation}
G(x^{\ast})=P_{\ast}\begin{bmatrix}
O&O\\
O&\Gast^{\NN} 
\end{bmatrix}P_{\ast}^{\top},\ \tred{Y}=P_{\ast}\begin{bmatrix}
\tred{\YUU}&O\\
O&O
\end{bmatrix}P_{\ast}^{\top},
\label{eq:past}
\end{equation}
where $\Gast^{\NN}\in \mathbb{S}^{r_{\ast}}_{++}$ is a {\it diagonal} matrix 
with $\rast:={\rm rank}G(\xast)$, the positive real eigenvalues of $G(x^{\ast})$ are aligned on the diagonal line, and 
$\tred{\YUU}\in \mathbb{S}^{m-r_{\ast}}_+$.
\tred{$\YUU\in \mathbb{S}^{m-\rast}_{++}$ does not necessarily hold}.
Without loss of generality, we may assume that the eigenvalues are placed in the ascending order on the diagonal.  
Needless to say, $P_{\ast}$ is a matrix whose columns are eigenvectors of $G(x^{\ast})$.
Partition the matrix $P_{\ast}$ as
$$P_{\ast }=[\PU,\PN],$$
where $\PU\in \R^{m\times (m-r_{\ast})}$ and $\PN\in \R^{m\times r_{\ast}}$.
Note that each column of $\PU$
represents an eigenvector of $G_{\ast}:=G(x^{\ast})$ which corresponds to the zero-eigenvalue of $G(x^{\ast})$, while that of $\PN$ does to a positive eigenvalue of $G_{\ast}$.
In terms of $\PU$ and $\PN$, the two equations in \eqref{eq:past} are transformed as 
\begin{equation}
\begin{bmatrix}
  \PU^{\top}G_{\ast}\PU&  \PU^{\top}G_{\ast}\PN\\
  \PN^{\top}G_{\ast}\PU&  \PN^{\top}G_{\ast}\PN
\end{bmatrix}=\begin{bmatrix}
O&O\\
O&\Gast^{\NN}
\end{bmatrix},\ \ 
\begin{bmatrix}
  \PU^{\top}\tred{Y}\PU&  \PU^{\top}\tred{Y}\PN\\
  \PN^{\top}\tred{Y}\PU&  \PN^{\top}\tred{Y}\PN
\end{bmatrix}=\begin{bmatrix}
\tred{\YUU}&O\\
O&O
\end{bmatrix}.
\label{eq:anotherform}
\end{equation}
We will often make use of formulation\,\eqref{eq:anotherform}.
For later use, we define the following notations: for the above $\Past = [\PU,\PN]$ and given $x,d\in \R^n$ and $\tred{\Yast}\in \mathbb{S}^m$, we write
\begin{eqnarray}
&\tred{\begin{bmatrix}
   \YUU_{\ast}&\YUN_{\ast}\\
\YNU_{\ast}&\YNN_{\ast}
\end{bmatrix}:=\begin{bmatrix}
\PU^{\top}\Yast\PU&\PU^{\top}\Yast\PN\\
\PN^{\top}\Yast\PU&\PN^{\top}\Yast\PN
\end{bmatrix}},\ 
\begin{bmatrix}
\GUU&\GUN\\
\GNU&\GNN
\end{bmatrix}:=\begin{bmatrix}
\PU^{\top}G(x)\PU&\PU^{\top}G(x)\PN\\
\PN^{\top}G(x)\PU&\PN^{\top}G(x)\PN
\end{bmatrix},&\label{eqn:1212-1}\\
&\begin{bmatrix}
\DeltaGUU(x;d)&\DeltaGUN(x;d)\\
\DeltaGNU(x;d)&\DeltaGNN(x;d)
\end{bmatrix}:=
\begin{bmatrix}
\PU^{\top}\DeltaG(x;d)\PU&
\PU^{\top}\DeltaG(x;d)\PN
\\
\PN^{\top}\DeltaG(x;d)\PU
&
\PN^{\top}\DeltaG(x;d)\PN
\end{bmatrix}.\label{eqn:1212-2}&
\end{eqnarray}

\subsection{Second-order optimality conditions and relevant properties}\label{subsec:sosc}
In this subsection, we review the second-order necessary/sufficient conditions for the NSDP.  Subsequently, we will describe the relevant properties briefly.  For more detailed explanations, we refer readers to, e.g., \citep{yamashita2012local, shapiro1997first} or \citep{bonnans2013perturbation}.

%To this end, we need to define the three concepts: Nondegenerancy condition, second-order conditions, and strict complemenarity condition.
\begin{definition}
  Let $x^{\ast}$ be a KKT point for the NSDP and consider the corresponding Lagrange multiplier set $\Lambda(x^{\ast})$. 
  Then, the nondegeneracy condition, strict complementarity conditions, and second-order condition are defined as follows:
  \begin{description}
  \item[\bf Nondegeneracy condition]
Let 
$r_{\ast}:={\rm rank}\,G(x^{\ast})$ and let $\{e_1,e_2,\ldots,e_{m-r_{\ast}}\}$ be an orthonormal basis of the null space of $G(x^{\ast})$.
Moreover, denote
$$
v_{ij}:=(e_i^{\top}\G_1(x^{\ast})e_j,\cdots, e_i^{\top}\G_n(x^{\ast})e_j)^{\top}\in \R^n\ (1\le i\le j\le m-r_{\ast}).
$$
We say that the nondegeneracy condition holds at $x^{\ast}$ if the vectors 
$v_{ij}\in \R^n\ (1\le i\le j\le m-r_{\ast})$ and $\nabla h_i(x^{\ast})\ (i=1,2,\ldots,\ell)$ are linearly independent.%\vspace{0.5em}
\item[\bf Strict complementarity condition]
  Let $\tred{Y}\in \mathbb{S}^m_{+}$ be a Lagrange multiplier matrix at $x^{\ast}$, which means that $G(x^{\ast})$ and $\tred{Y}$ satisfies the complementarity condition\,\eqref{al:comp}. 
We say that the strict complementarity condition holds at $(x^{\ast}, \tred{Y})$ if 
$G(x^{\ast})+\tred{Y}\in \mathbb{S}^m_{++}$, which is equivalent to ${\rm rank}\,G(x^{\ast})+{\rm rank}\,\tred{Y}=m$
under \eqref{al:comp}.
\item[\bf Second-order conditions]
  We say that the second-order necessary (resp., sufficient) condition holds at $x^{\ast}$ if 
\begin{equation}
  \sup_{(Y,z)\in \Lambda(x^{\ast})}d^{\top}\left(\nabla_{xx}^2L(x^{\ast},Y,z)+\Omega(x^{\ast},Y)\right)d\ge(\mbox{resp., }>)0,\ \ \forall
  \tred{d\in C(x^{\ast})\setminus\{0\}}, \label{eq:soc}
\end{equation}
where $C(x^{\ast})$ is the critical cone at $x^{\ast}$ and specifically represented as 
\begin{equation} 
C(x^{\ast})= \left\{d\in \R^n\mid \nabla f(x^{\ast})^{\top}d=0, \nabla h(x^{\ast})^{\top}d=0,\ \DeltaG(x^{\ast};d)\in T_{\mathbb{S}^m_+}(G(x^{\ast}))
\right\}.\label{eq:critical}
\end{equation}
Here, $T_{\mathbb{S}^m_+}(G(x^{\ast}))$ denotes the tangent cone of $\mathbb{S}^m_+$ at $G(x^{\ast})$ and is represented specifically as
\begin{equation*}
T_{\mathbb{S}^m_+}(G(x^{\ast}))=\left\{X\in \mathbb{S}^m\mid \PU^{\top}X\PU(=X^{\UU})\in S^{r_{\ast}}_{+}\right\}.\label{eq:tangentS}
\end{equation*}
Moreover, for any $x\in \R^n$ and $Y\in \mathbb{S}^m$, 
$\Omega(x,Y)$ denotes the matrix in $S^n$ whose $(i,j)$-th entry is given as 
$$
(\Omega(x,Y))_{i,j}:= 2Y\bullet \Gi(x)G(x)^{\dag}\G_j(x)
$$
for $i,j=1,2,\ldots,n$, where $G(x)^{\dag}$ denotes the Moore-Penrose inverse matrix of $G(x)$.
\end{description}
\vspace{0.5em}
\end{definition}
\begin{remark}\label{rem:0314-1}
The nondegeneracy condition at $x^{\ast}$ is a constraint qualification for the NSDP and yields the MFCQ.
It reduces to the linear independence constraint qualification (LICQ) when nonlinear optimization is considered. As with the LICQ, 
the Lagrange multiplier set $\Lambda(x^{\ast})$ is a singleton under the nondegeneracy condition.
\end{remark}

The term $d^{\top}\Omega(x^{\ast},Y)d$ in \eqref{eq:soc} is called the {\it sigma term} for the semi-definite constraint $G(x)\in \mathbb{S}^m_+$. We refer readers to \citep{bonnans2013perturbation} for a precise description of its background and properties.
In the following lemma, the sigma term is expressed more specifically, thereby being ensured to be nonnegative.
%The following lemma 
\begin{lemma}\label{lem:sigma}
  %For a KKT triplet $w^{\ast}=(\xast,\Yast,\zast)$
\tred{For $\tred{Y}\in \mathbb{S}^m_+$ such that $\Gast \tred{Y}=O$ and a direction $d\in \R^n$, it holds that}
\begin{align*}
d^{\top}\Omega(x^{\ast},Y)d&=2
\Tr\left(
\YUU\DeltahGNU(x^{\ast};d)
(\Gast^{\NN})^{-1}\DeltahGUN(x^{\ast};d)
\right)\\
&=2\left\|
(\YUU)^{\fr}\DeltahGNU(x^{\ast};d)
(\Gast^{\NN})^{-\fr}
\right\|^2_{\rm F},
\end{align*}
where $\YUU$ and $\Gast^{\NN}$ are defined in \eqref{eq:anotherform}, and moreover $\DeltahGNU$ and $\DeltahGUN$ in \eqref{eqn:1212-2}.
\end{lemma}
\proof{{\rm Proof}.}
By straightforward calculation. See Appendix\,\ref{subsec:lemma-1} for details.
$\hfill\Box$
\endproof\vspace{0.5em}

When we consider the standard nonlinear optimization where the nonnegative cone is set in the NSDP in place of the semidefinite cone, the sigma term always vanishes because $\DeltahGNU(x^{\ast};d)=O$ holds for any $d$ in the above lemma, and thus it never appears in the second-order conditions. 
In contrast, in the NSDP, the sigma term reflects curvature of $\mathbb{S}^m_+$ and is nonnegative for any $d\neq 0$ and $Y\in \mathbb{S}^m_{+}$ as shown in \cref{lem:sigma}. 
With the help of this term, the second-order condition is more likely to hold
even when $\nabla_{xx}^2L$ is not positive semidefinite over the critical cone. 
%Meanwhile, its complicated structure often brings about difficulty of analyzing several properties which were already shown only for nonlinear optimization.
\tred{However, this term makes the analysis for the NSDP more complicated than in nonlinear optimization.}
%In order to analyze the limiting behavior of the central path of the NSDP, it is crucial to study a relationship  between the sigma term and the Hessian of the log-barrier function.

Lastly, we mention useful facts associated with the second-order conditions in the following two necessary and sufficient optimality conditions.
\paragraph{Second-order necessary optimality for the NSDP\,\tred{\citep[Theorem~3.45,5.88]{bonnans2013perturbation}}}
Let $x^{\ast}\in \R^n$ be a local optimum of NSDP\,\eqref{al:nsdp} and suppose that the MFCQ holds there.
  Then, the second-order necessary condition holds at $x^{\ast}$.   
  \paragraph{Second-order sufficient optimality for the NSDP\,\tred{\citep[Theorem~5.89]{bonnans2013perturbation}}}
%  Let $x^{\ast}\in \R^n$ be local optimum of NSDP\,\eqref{al:nsdp} and suppose that the MFCQ holds there.
  %Then, the second-order sufficient condition holds at $x^{\ast}$.
  Suppose that $x^{\ast}$ is a KKT point of NSDP\,\eqref{al:nsdp} and, furthermore, the second-order sufficient condition holds.  
Then, $x^{\ast}$ is a strict local optimum of NSDP\,\eqref{al:nsdp}. In particular, the quadratic growth condition holds, that is,
there exists some $q>0$ and vicinity $\mathcal{N}(x^{\ast})$ of $x^{\ast}$ such that 
$f(x)-f(x^{\ast})\ge q\|x-x^{\ast}\|^2$ for all $x\in \mathcal{N}(x^{\ast})\cap \mathcal{F}$.

\subsection{BKKT conditions and central path}\label{subsec:bkkt}
In this section, we introduce the barrier KKT (BKKT) conditions for the NSDP.
The BKKT conditions are composed of \eqref{al:kkt1}, \eqref{al:kkt3}, and the following perturbed conditions for \eqref{al:xy}: for $\mu>0$,
%obtained by perturbing the complementary condition of the form\,\eqref{eq:xy} in terms of $\mu>0$ as
\begin{equation}
G(x)Y=\mu I,\ G(x)\in \mathbb{S}^m_{++},\ Y\in \mathbb{S}^m_{++}.\label{eq:bkkt}
\end{equation}
The parameter $\mu$ is often referred to as barrier parameter, and $x$ and $(x,Y,z)$ satisfying the BKKT conditions are called a BKKT point and BKKT triplet, respectively. 
It is worth mentioning that condition\,\eqref{eq:bkkt} is equivalent to the condition obtained by replacing $O$ with $\mu I$ in \eqref{eq:aho} or \eqref{eq:mt}.
As $\mu$ gets closer to 0, BKKT points are expected to approach the set of KKT points for the NSDP.
%which motivates us to track the central path with driving $\mu\to 0$.
A basic algorithmic policy of primal-dual interior point methods is to track BKKT triplets while driving $\mu$ to 0, so as to reach a KKT triplet.
%this idea as a basic algorithmic policy.
%Now, the central path is defined via the BKKT conditions as follows:
In this paper, we will refer to a path formed by BKKT triplets as a {\it central path}.    
%% \begin{definition}
%%   Let
%%   $
%% \mathcal{C}:=\{w\in\mathcal{W}\mid w\mbox{ satisfies the BKKT conditions with $\mu>0$}\}.
%%   $
%% The set $\mathcal{C}$ is called a central path of NSDP\,\eqref{al:nsdp}.\footnote{Strictly speaking, 
%% $\mathcal{C}$ may not be a path unlike linear SDPs.
%% %In fact, as will be proved, $\mathcal{C}$ forms a path locally.
%% Nevertheless, we call it {\it path} for convention. 
%% }
%% \end{definition}

%As $\mu$ gets closer to 0, BKKT points are expected to approach the set of KKT points for the NSDP, which motivates us to track the central path with driving $\mu\to 0$.
%Many primal-dual interior point methods employs this idea as a basic algorithmic policy.
%\section{Preriliminary results: }

\section{Main analysis}\label{sec:mainresults}
\subsection{Assumptions and outline of analysis}\label{sec:3-1}
Throughout Section\,\ref{sec:mainresults}, $x^{\ast}$ denotes a KKT point of the NSDP, and 
%In this section, We will study the limiting behavior of the central path around $w^{\ast}$ under the following assumptions:
%The first set concerns the assumptions on $w^{\ast}$. 
%Now, let us specify the assumptions on $\xast$ in our setting.
%We impose the following assumptions on $\xast$.
is assumed to satisfy the following:
\begin{assumption}\label{assumA}
The KKT point $x^{\ast}$ satisfies the following three conditions:
\begin{enumerate}
\item 
There exists a Lagrange multiplier matrix $\tred{Y}\in \mathbb{S}^m_{+}$ satisfying the strict complementarity condition.
\item The \tred{enhanced} second-order sufficient condition (\tred{\SSOSC}) holds:
for all $(Y,z)\in \Lambda(\xast)$, it holds that
%\tred{there exists some $\kappa>0$ such that} 
{
%\begin{equation}
  $
  d^{\top}\left(\nabla_{xx}^2L(x^{\ast},Y,z)+\Omega(x^{\ast},Y)\right)d>0,\ \forall d\in C(x^{\ast})\setminus\{0\}
  $.
  }
%\label{eq:ssoc0}
%\end{equation}}
%\item the full column rank of $\nabla h(x^{\ast})$;
\item The MFCQ holds at $x^{\ast}$. 
%\item $\Gast\in \mathbb{S}^m_{+}\setminus \mathbb{S}^m_{++}$, namely, 
%$\Gast$ has at least one zero-eigenvalue. 
\end{enumerate}
\end{assumption}
%Under the above strict complementrity condition, 
%the critical cone $C(x^{\ast})$ defined by \eqref{eq:critical} reduces to a linear space expressed as 
%\begin{equation} 
%C(x^{\ast})= \left\{d\in \R^n\mid \nabla f(x^{\ast})^{\top}d=0, \nabla h(x^{\ast})^{\top}d=0,\ \mathcal{J}G(x^{\ast})d\in T_{\mathbb{S}^m_+}(G(x^{\ast}))
%\right\},\label{eq:criticalsc}
%\end{equation}\tblue{[[Wrong. Later fixed!]]}
%where $T_{\mathbb{S}^m_+}(G(x^{\ast}))$ denotes the tangent cone of $\mathbb{S}^m_+$ at $G(x^{\ast})$.
%Even in the absence of 
%Mind that the nondegeneracy condition is not included in the above assumptions.
\tred{The above {\SSOSC} is indeed stronger than the second-order sufficient condition\,(SOSC) defined in \eqref{eq:soc}, because, with arbitrarily chosen $(\overline{Y},\bar{z})\in \Lambda(\xast)$, we have 
$\sup_{(Y,z)\in \Lambda(x^{\ast})}d^{\top}\left(\nabla_{xx}^2L(x^{\ast},Y,z)+\Omega(x^{\ast},Y)\right)d
\ge d^{\top}\left(\nabla_{xx}^2L(x^{\ast},\overline{Y},\bar{z})+\Omega(x^{\ast},\overline{Y})\right)d>0
$
for any $d\in C(\xast)\setminus\{0\}$, where the last inequality is due to the {\SSOSC}.
}
Under the \tred{\SSOSC}, $x^{\ast}$ is a strict local optimum of the NSDP since the SOSC follows from the {\SSOSC} as shown above.
\tred{See also {\it Second-order sufficient optimality for the NSDP} at the end of subsection\,\ref{subsec:sosc}.}
The \tred{\SSOSC} holds, for example, when $f$ is strongly convex and $G$ and $h$ are affine. 
%It is worth mentioning that this \tred{\SSOSC} differs from the one proposed by \citet{sun2006strong}. 
It \tred{can be seen} as a straightforward generalization of the strong second-order condition (SSOSC) considered by \citet{wright2002properties} for nonlinear optimization.
\tred{Though one may think it natural to refer to the condition as SSOSC,
we call it {\SSOSC} so as to distinguish it from the SSOSC for the NSDP studied by \citet{sun2006strong}.}
%and differs from \tred{the SSOSC} proposed by \citet{sun2006strong} \tred{for the NSDP}. 
Under the presence of the MFCQ, we ensure compactness and convexity of the Lagrange multiplier set $\Lambda(\xast)$, but $\Lambda(\xast)$ is not necessarily a singleton(cf. Remark\,\ref{rem:0314-1}). %This is a quite different situation from the case where the nondegeneracy condition is supposed. (cf. Remark\,\ref{rem:0314-1}) 
{Note that $\nabla_{xx}^2L$ is continuous, and so is $\Omega(\xast,Y)$ with respect to $Y\in \mathbb{S}^m_{+}$ such that $G(\xast)Y=O$ from \cref{lem:sigma}.
  This fact, the compactness of $\Lambda(\xast)$, and the {\SSOSC} guarantee that there exists some $\kappa>0$ such that
\begin{equation}
\inf_{(Y,z)\in \Lambda(\xast)}d^{\top}\left(\nabla_{xx}^2L(x^{\ast},Y,z)+\Omega(x^{\ast},Y)\right)d\tred{\ge \kappa\|d\|^2},\ \ \forall d\in C(x^{\ast})\setminus \{0\}.\label{eq:ssoc}
\end{equation}
%which will be often used.
}
%$x$\{(Y,z)\in \mathbb{S}^m_+\times \R^{p}\mid \mbox{$(Y,z)$ is a pair of Lagrange multiplier matrix and vector at $x^{\ast}$}\}.$$
\paragraph{Goal and outline of the analysis:} The goal of the whole analysis we will conduct is to prove that under the above assumptions, 
there exists a unique and smooth central path converging to the KKT triplet 
\begin{equation}
\wa:=(\xast,\Ya,\za),\label{eq:analyticKKT}
\end{equation}
where $(\Ya,\za)\in\Lambda(\xast)$ is called an analytic center at $\xast$, defined formally in the next subsection. 
In order to achieve this goal, we will prove the following claims in order:\vspace{0.1em}
\begin{description}
\item[\bf Claim\,(i)] There exists a sequence of BKKT triplets
$\left\{w^k=(x^k,\Yk,\zk)\right\}$ converging to the KKT triplet $\wa$ (cf. \cref{thm:analytic} in subsection\,\ref{subsec:claim1}).\vspace{0.1em}
\item[\bf Claim\,(ii)] 
%% A direction of BKKT point $x^k$ approaching the KKT point $\xast$ converges to some direction $\xiast$.
%% This $\xiast$ is a unique $x$-component solution of a certain linear equation related to tangential directions of the central path (cf. \cref{prop:xiast} and \cref{cor1} in subsection\,\ref{subsec:claim2}).
  \tred{Any sequence of BKKT points $\{\xk\}$ approaches the KKT point $\xast$ asymptotically along a certain nonzero direction $\xiast\in \R^n$ in the sense that $\lim_{k\to\infty}\frac{\xast-\xk}{\|\xast-\xk\|}=\frac{\xiast}{\|\xiast\|}$ (cf. \cref{prop:xiast} and \cref{cor1} in subsection\,\ref{subsec:claim2}). This $\xiast$ is a unique $x$-component solution of a certain linear equation system related to the BKKT conditions.
  }
\vspace{0.1em}
\item[\bf Claim~(iii)] For any sufficiently small barrier parameter $\mu$, a corresponding BKKT point $\xmu$ exists uniquely \tred{in the open ball
$\{x\in \R^n\mid \|x-\xast-\mu\xiast\|<\rho\mu\|\xiast\|\}$, where $\rho>0$ is a certain small constant. Moreover, the Hessian of a certain barrier function is nonsingular at $\xmu$.}
  %as a vertex and $\xiast$ as an axis
  (cf. \cref{prop:0603-1} in subsection\,\ref{subsec:claim4}).\vspace{0.1em}
\end{description}

With the help of the above claims and the classical implicit function theorem, 
we will prove our main claim of the goal (cf. \cref{lastthm} in subsection\,\ref{sec:claim5} and \cref{lastthm2} in subsection\,\ref{sec:unique}).
Mind that henceforth, several proofs are deferred to the Appendix for the sake of readability. 

\paragraph{Additional notations and symbols used hereafter} In the remaining of Section\,\ref{sec:mainresults} and the Appendix, we will use the symbols and the notations defined in \eqref{eq:past}-\eqref{eqn:1212-2} in addition to those introduced at the end of 
Section\,\ref{sec:1}. In particular, $\Past=[\PU,\PN]$ is an {\it arbitrarily} chosen orthogonal matrix defined for $G(\xast)$ so that \eqref{eq:past} holds.
For the sake of simplicity, we often write 
$$G_{\ast}:=G(\xast),\ G_k:=G(\xk).$$
Besides, we will make use of  
$G_{\ast}^{ind}$ and $G_k^{ind}$
$(ind\in \{{\UU},{\NU},{\Un},{\NN}\})$ defined  by replacing $G$ and $Y$ in \eqref{eqn:1212-1} with $G_{\ast}$and $G_k$, respectively. 
Furthermore,  $\Yk^{ind}$ and $\Ya^{ind}$ $(ind\in \{{\UU},{\NU},{\Un},{\NN}\})$ are defined in the same way using $\Yk$ and $\Ya$.
%we aim to prove main results about a central path of the NSDP. 
%The organization of the section is three-fold.
%\subsection{Assumptions}\label{subsec:assum}

\subsection{Existence of analytic center for NSDP}\label{subsec:analytic}
The analytic center for the NSDP at $x^{\ast}$ is formally defined as follows:  
\begin{definition}[Analytic center for NSDP\,\eqref{al:nsdp}]
%Let $P_{\ast}=[\PU,\PN]\in \R^{m\times m}$ be the orthogonal matrix that is chosen in \eqref{eq:past}. 
We say that $(\Ya,\za)\in \Lambda(x^{\ast})$ is an analytic center of NSDP\,\eqref{al:nsdp} at $x^{\ast}$ if it is an optimum of 
%it holds that 
%$$
%\Ya=P_{\ast}\begin{bmatrix}
%\YUUa&O\\
%O&O
%\end{bmatrix}P_{\ast}^{\top},
%$$
%where 
%$(\Ya^{\UU},\za)$ solves the following convex problem with linear constraints:
\begin{equation}
%\min_{\YUU,z}
%\min_{\YUU,z}
%\ -\log\det \YUU \mbox{ ${\rm s.t}$. }\nabla f(x^{\ast})-\mathcal{J}\GUU(x^{\ast})^{\ast}\YUU + \nabla h(x^{\ast})z=0,\ \YUU\in \mathbb{S}^{m-\rast}_{++}
\min\ -\log\det \YUU \mbox{ ${\rm s.t}$. }{(Y,z)\in \Lambda(\xast)}.\label{eq:analy}
\end{equation}
%with 
%$\mathcal{J}\GUU(x^{\ast})^{\ast}\YUU:=[\PU^{\top}\Gi(x^{\ast})\PU\bullet \YUU]_{i=1}^{n}\in \R^{n}$.
Here, we define $\log 0 := - \infty$ by convention.
\end{definition}
%In view of relation\,\eqref{eq:past}, problem\,\eqref{eq:analy} is rewritten as 
%\begin{align}
%\begin{array}{cc}
%\displaystyle{\min_{(Y,z)\in \mathbb{S}^m\times \R^s}}&-\log\det \YUU \\
%\mbox{ ${\rm s.t}$. }
%%&\nabla f(x^{\ast})-\mathcal{J}         \GUU(x^{\ast})^{\ast}\YUU + \nabla h(x^{\ast})z=0,\\ 
%&\nabla f(x^{\ast})-\mathcal{J}G(x^{\ast})^{\ast}Y+ \nabla h(x^{\ast})z=0,\\
%      &\YUN=\YNU=O,\ \YNN=O,\\
%      &\YUU\in \mathbb{S}^{m-\rast}_{+}.\\
%\end{array}
%\label{eq:analy2}
%\end{align}
%where $\mathcal{J}\GUU(x^{\ast})^{\ast}\YUU:=\left[\left(\PU^{\top}\Gi(x^{\ast})\PU\right)\bullet \YUU\right]_{i=1}^{n}\in \R^{n}$.

%From these insights, the following proposition is redobtained.
%We summarize these insights as the following proposition. 
%From these insights, we have the following proposition.

In the next proposition, we ensure existence and uniqueness of the analytic center at $\xast$.
In other words, the KKT triplet $\wa$ defined in \eqref{eq:analyticKKT} is well-defined. 
\begin{proposition}\label{prop:analyticunique} 
Suppose that Assumption\,\ref{assumA} holds. 
Then, an analytic center of NSDP\,\eqref{al:nsdp} at $x^{\ast}$ exists uniquely.
In particular, $(\Ya,\za)\in \mathbb{S}^m\times \R^s$ is the analytic center at $\xast$ if and only if $(\Ya,\za)\in \Lambda(\xast)$ and 
there exists some vector $v\in\R^n$ such that
\begin{equation}	
\DeltaGUU(\xast;v)=(\Ya^{\UU})^{-1},\ \nabla h(x^{\ast})^{\top}v=0.\label{eq:proofanalytic}
\end{equation}

%In particular, $(\Ya^{\UU},\za)$ is a unique optimum of 
%the following convex problem with linear constraints:
%%with a feasible region wider than \eqref{eq:analy0}:
%\begin{equation}
%\min_{\YUU,z}
%\ -\log\det \YUU \mbox{ ${\rm s.t}$. }\nabla f(x^{\ast})-\mathcal{J}\GUU(x^{\ast})^{\ast}Z + \nabla h(x^{\ast})z=0,\ Z\in \mathbb{S}^{m-\rast}_{++}
%\label{eq:analy}
%\end{equation}
%with 
%$\mathcal{J}\GUU(x^{\ast})^{\ast}Z:=[\PU^{\top}\Gi(x^{\ast})\PU\bullet Z]_{i=1}^{n}\in \R^{n}$.
\end{proposition}
\proof{\rm Proof.}
See Appendix\,\ref{app:prop:analyticunique}.
\hfill\Halmos
\endproof\vspace{0.5em}

%For later use, we write down the necessary and sufficeint conditions of \eqref{eq:analy} via the KKT conditions of \eqref{eq:analy2}: 
%$(\Ya^{\UU},\za)$ is the global optimum of \eqref{eq:analy} if and only if it is feasible and moreover
%there exists some vector $v\in\R^n$ such that
%%\begin{equation}	
%%\sum_{i=1}^n v_iE_{\ast}^{\top}\Gi(x^{\ast}) E_{\ast}=(\Ya^{\UU})^{-1},\ \nabla h(x^{\ast})^{\top}v=0.\label{eq:proofanalytic}
%%\end{equation}
%\begin{equation}	
%\DeltaGUU(\xast;v)=(\Ya^{\UU})^{-1},\ \nabla h(x^{\ast})^{\top}v=0.\label{eq:proofanalytic}
%\end{equation}

\subsection{{\bf Proof of Claim~(i):} convergence of BKKT triplets to KKT triplet with analytic center}\label{subsec:claim1}
In this subsection, we will prove that there exists a sequence of BKKT points which converges to the KKT point $\xast$. Moreover, we will show that 
the corresponding dual sequence converges to the analytic center $(\Ya,\za)$.  
%First, we consider the relationship between a central path and the following log-barrier function for the NSD

Let us define the following log-barrier function for the NSDP: for each $\mu>0$
\begin{equation}
  \psi_{\mu}(x):=f(x)-\mu\log\det G(x).\label{eq:231213-4}
\end{equation}
%For later use, we give the following lemma regarding the hessian of $\psi_{\mu}$.  The proof is omitted as it is done by straightforward calculation.
%\begin{lemma}\label{lem:positive}
%  For simplicity of the expression, we write $G:=G(x)$,
%%\begin{align*}
%$\G_i:=\G_i(x)$, and $\H_{ij}:=\frac{\partial^2G(x)}{\partial x_i\partial x_j}$
%%\end{align*}
%for each $i,j\in \{1,2,\ldots,n\}$. Then, it holds that 
%  \begin{equation}
%    \nabla^2\psi_{\mu}(x)
%    =\nabla^2f(x)+\mu\left[\left(G^{-1}\G_iG^{-1}\right)\bullet \G_j\right]_{i,j}-\mu\left[G^{-1}\bullet \H_{ij}\right]_{ij}.
%    \label{eq:hessipsimu}
%  \end{equation}
%  In particular, for any $v\in \R^n$
%\begin{align}
%  v^{\top}\nabla^2\psi_{\mu}(x)v=
%  v^{\top}\nabla^2f(x)v
%  -\mu\sum_{i=1}^n\sum_{j=1}^nv_iv_jG^{-1}\bullet \H_{ij}+\mu \left\|G^{-\fr}\DeltaG(x;v) G^{-\fr}\right\|_{\rm F}^2.\label{eq:hessipsimudd}
%\end{align}
%\end{lemma}

%\begin{theorem}
%  There exist some $\ol{\mu}>0$ and a unique path $x(\cdot):[0,\bar{\mu}]\to \R$ such that
%\begin{enumerate}
%  \item $x(\cdot)$ is smooth on $(0,\ol{\mu})$;
%  \item it is continuous at $\mu=0$ and also $x(0)=x^{\ast}$.
%\end{enumerate}
%\end{theorem}
The following proposition states that 
there exists a sequence of local optima of barrier penalized NSDPs converging to $\xast$.
Such local optima are BKKT points of the NSDP locally around $\xast$.
%the relationship between local optima of the barrier penalized NSDP with $\psi_{\mu}$
%and BKKT points. 
\begin{proposition}\label{thm:barrier}
Let Assumption\,\ref{assumA} hold and $\{\mu_k\}\subseteq \R_{++}$ be an arbitrary {\it decreasing} sequence converging to 0.
  Then, there exists a sequence $\{x^k\}\subseteq \R^n$ such that $\lim_{k\to\infty}x^k=x^{\ast}$ and, for any $k\ge \bK$ with
  $\bK$ sufficiently large, $x^k$ is a local optimum of
  \begin{equation}
  \min\ \psi_{\muk}(x)\ \ {\rm s.t.}\ \ h(x)=0,\ G(x)\in \mathbb{S}^m_{++}.\label{eq:log}
  \end{equation}
  %$\nabla h(\xk)$ is of full column rank, and 
%  and furthermore $\xk$ is a BKKT point with $\mu=\muk$ of NSDP\,\eqref{al:nsdp} by re-taking larger $\bK$ if necessary.
  %In particular, 
  %by re-taking $\bK$ if necessary, 
  %for any $k\ge \bK$.     
\end{proposition}
\proof{\rm Proof.}
The proof is analogous to those of classical results as to penalty methods\,\citep{luen2008}, although it is different in dealing with the log determinant function and the semidefinite constraint. 
Nonetheless, the precise proof is given in Appendix\,\ref{sec:appendix1} for completeness.
\hfill\Halmos\endproof\vspace{0.5em}

In Proposition\,\ref{thm:barrier}, as $\nabla h(\xast)$ is of full column rank, so is $\nabla h(\xk)$ for any $k\ge \bK$ with $\bK$ large enough, and thus the KKT conditions for \eqref{eq:log} holds at $\xk$. From the KKT conditions together with  
$\nabla \psi_{\mu}(x)=\nabla f(x)-\mu \mathcal{J}G(x)^{\ast}G(x)^{-1}$, there exists $\zk\in \R^s$ such that  
$$
\nabla f(\xk)-\muk \mathcal{J}G(\xk)^{\ast}G_k^{-1}+ \nabla h(\xk)\zk = 0,\ h(\xk)=0,\ G_k\in\mathbb{S}^m_{++}, 
$$
which together with $\Yk:=\muk \Gk^{-1}\in \mathbb{S}^m_{++}$ implies that 
$\xk$ and 
$\wk:=(\xk,\Yk,\zk)\in \W_{++}$ are BKKT point and BKKT triplet for each $k\ge \bK$, respectively.

In summary, as a consequence of \cref{thm:barrier}, 
given a decreasing sequence $\{\muk\}\subseteq \R_{++}$ converging to 0, there exists an integer $\bK>0$ and 
$\left\{\wk\right\}\subseteq \W_{++}$ with $\wk=(\xk,\Yk,\zk)$
such that 
\begin{equation}
  \lim_{k\to \infty}\xk=\xast,\ \Yk=\muk G_k^{-1}\in \mathbb{S}^m_{++}\label{eq:231230-2}
\end{equation}
and, for each $k\ge \bK$, 
\begin{eqnarray*}
\nabla h(\xk):\mbox{ full column rank},\ \wk:\mbox{ BKKT triplet with barrier parameter $\muk$}.
\end{eqnarray*}
%$\lim_{k\to \infty}(\muk,\xk)=(0,\xast)$ and moreover, for each $k$,  
%$\Yk=\muk \Gk^{-1}\in \mathbb{S}^m_{++}$, $\nabla h(\xk)\in \R^{n\times \dimeq}$ is of full column rank, and 
%$\xk$ and $\wk$ are a BKKT point and a BKKT triplet with $\mu=\muk$, respectively.
In what follows, for the sake of brevity, we assume $\bK = 0$. Moreover, we suppose that $\xk\neq \xast$ for all $k$ without loss of generality and define  
\begin{equation}
d^k:=\xk-\xast.\label{eq:dk}
\end{equation}
%% Notice that
%% \begin{equation}
%%   \tred{\lim_{k\to\infty}d^k=0}\label{eq:dkconv}
%% \end{equation}
%% by construction.
Hereafter, we focus on those sequences $\{w^k\}$ and $\{d^k\}$.
\begin{remark}\label{rem:20211217}
In a quite similar manner to the proof of \citet[Theorem~1]{yamashita2012primal}, we ensure that, under the MFCQ at $x^{\ast}$, the sequence
$\{(Y_k,z^k)\}$ is bounded, and its accumulation point together with $x^{\ast}$ fulfills the KKT conditions of NSDP\,\eqref{al:nsdp}.
\end{remark}\vspace{0em}

In fact, the whole sequence $\{(Y_k,z^k)\}$ converges to the analytic center $(\Ya,\za)$.
%and thus has at least one accumulation point, say
%In the remainder of this section, we devote ourselves into proving this assertion by showing that 
%any accumulation point of $\{(Y_k,z^k)\}$, denoted by $(Y_{\ast},z^{\ast})$, is idential to $(\Ya,\za)$.
%We next compare the convergence speeds of $\{\muk\}$ and $\{\|d^k\|\}$. 
%As verified in later on, we have  
%\begin{equation*}
%\mu_k=\Theta(\|d^k\|)
%\end{equation*}
%under the present setting. 
To prove this claim, we first present the following proposition, which claims that the convergence speeds of $\{\muk\}$ and $\{\|d^k\|\}$ towards zero are equivalent.
\begin{proposition}\label{prop:mukdk}
Suppose that Assumption\,\ref{assumA} holds. Then, we have 
\begin{equation*}
\mu_k=\Theta(\|d^k\|).
\end{equation*}
\end{proposition}
\proof{\rm Proof.}
See Appendix\,\ref{app:prop:mukdk}.
\hfill\Halmos\endproof\vspace{0.5em}

Using this proposition, the convergence to the analytic center can be established.
%\proof{\rm Proof.}
%In \cref{prop:liminfmuk}, we have already proved
%$\liminf_{k\to\infty}\frac{\mu_k}{\|d^k\|}>0$. In order to show the desired claim, there remains to prove 
%$\frac{\mu_k}{\|d^k\|}=\O(1)$. 
%
%As $w^k$ satisfies the BKKT conditions, we obtain, for each $k$,
%\begin{align}
%\frac{\mu_kI_{r_{\ast}}}{\|d^k\|}&=\frac{\PU^{\top}G_kY_k\PU}{\|d^k\|}\notag \\
%                            &=\frac{\PU^{\top}\left(G(x^{\ast})+\mathcal{J}G(x^{\ast})(d^k)+\O(\|d^k\|^2)\right)\begin{bmatrix}\PU\ \PN
%                            \end{bmatrix}
%                            \begin{bmatrix}
%                            \PU^{\top}\\
%                            \PN^{\top}
%                            \end{bmatrix}
%                            Y_k\PU}{\|d^k\|},\notag
%\end{align}
%which together with the fact from \eqref{eq:anotherform} that 
%$$
%\PU^{\top}G(x^{\ast})\PU=O,\ 
%\begin{bmatrix}
%                            \PU^{\top}\Yast\PU\\
%                            \PN^{\top}\Yast\PU
%                            \end{bmatrix}=\begin{bmatrix}
%                            {\YUU_{\ast}}\\
%                            O
%                            \end{bmatrix}$$
% yields 
%\begin{equation*}
%\lim_{k\to\infty}\frac{\mu_kI_{r_{\ast}}}{\|d^k\|}=\begin{bmatrix}\DeltahGUU(x^{\ast};\hd^{\ast})\YUU_{\ast}\\O\end{bmatrix}.
%\end{equation*}
%This means that the sequence $\left\{\frac{\mu_k}{\|d^k\|}\right\}$ is bounded, and we thus obtain the desired consequence. The proof is complete.
%\hfill\Halmos\endproof
%

\begin{theorem}\label{thm:analytic}
Suppose that Assumption\,\ref{assumA} holds. Then, 
the whole sequence $\left\{(Y_k,z^k)\right\}$ converges to the analytic center $(\Ya,\za)$ of NSDP\,\eqref{al:nsdp} at $x^{\ast}$,
that is, $\lim_{k\to\infty}\wk=\wa$.
\end{theorem}
\proof{\rm Proof.}
\tred{Note that $\lim_{k\to \infty}d^k=0$ from \eqref{eq:231230-2} and \eqref{eq:dk}.
Also, note that $\{(\Yk,\zk)\}$ is bounded as was explained in Remark\,\ref{rem:20211217}, and let $(\Yast,\zast)$ be an arbitrary accumulation point of $\{(\Yk,\zk)\}$.}
For each $k\ge 0$, define $\hd^k:=\frac{d^k}{\|d^k\|}.$
Since $\{\hd^k\}$ is bounded, it has at least one accumulation point, say $\hd^{\ast}$. 
Choose an arbitrary subsequence $\{\hd^k\}_{k\in \K}$ which converges to $\hd^{\ast}$.
From \cref{prop:mukdk}, $\left\{\frac{\mu_k}{\|d^k\|}\right\}_{k\in\K}$ is bounded and any accumulation point, say $\balpha\in \R$, is positive. Without loss of generality, we assume that
$\left\{\frac{\mu_k}{\|d^k\|}\right\}_{k\in\K}$ \tred{and $\{(\Yk,\zk)\}_{k\in\K}$
  to $\balpha>0$ and $(\Yast,\zast)$, respectively,}
  by taking a subsequence further if necessary.

%Note that $\|Y_k\|_{\rm F}=\O(1)$ as was explained in Remark\,\ref{rem:20211217}.
Recalling that $\PUN=\begin{bmatrix}\PU,\ \PN\end{bmatrix}$ is orthogonal, we have 
\begin{align}
\frac{\mu_kI_m}{\|d^k\|} &= \frac{\PUN^{\top}G_kY_k\PUN}{\|d^k\|}\notag \\
                              &= \frac{
                              \PUN^{\top}\left(G(x^{\ast}) + \DeltaG(x^{\ast};d^k)+\O(\|d^k\|^2)\right)\PUN\PUN^{\top}Y_k\PUN 
                              }{\|d^k\|}\notag \\
                              &=
                              \begin{bmatrix} O& O\\
\frac{1}{\|d^k\|}\Gast^{\NN}\YUN_k &\frac{1}{\|d^k\|}\Gast^{\NN}\YNN_k 
\end{bmatrix}
+ \PUN^{\top} \DeltaG(x^{\ast};\hd^k)\PUN\left(\PUN^{\top}Y_k\PUN\right)
+\O(\|d^k\|).\label{al:0409}
\end{align}
Taking into account that $\lim_{k\in\K\to\infty}\PUN^{\top}Y_k\PUN =\begin{bmatrix}
\YUU_{\ast}&O\\
O& O
\end{bmatrix}$ 
with $\YUU_{\ast}=\PU^{\top} Y_{\ast}\PU$
and driving $k\in \K\to \infty$ in the $(1,1)$-block component of \eqref{al:0409},
we obtain $\DeltaGUU(x^{\ast};\hd^{\ast})\YUU_{\ast}=\balpha I_{r_{\ast}}$, implying 
\begin{equation}
(\Yast^{\UU})^{-1} = \DeltaGUU(x^{\ast};\balpha^{-1}\hd^{\ast}).
\label{eq:0409-2}
\end{equation}
Moreover, for each $k\in\K$, it holds that 
\begin{align*}
0&=\frac{h(x^k)}{\|d^k\|}=\frac{h(x^{\ast}) + \nabla h(x^{\ast})^{\top}d^k+\O(\|d^k\|^2)}{\|d^k\|}=\nabla h(x^{\ast})^{\top}\widetilde{d^k}+\O(\|d^k\|),
%   &\hspace{2em}\longrightarrow  \nabla h(x^{\ast})^{\top}\widetilde{d^{\ast}}\ \ (k\in\L\to\infty)\\ 
\end{align*}
which along with driving $k(\in\K)\to\infty$ and multiplying $\balpha^{-1}$ 
implies $\nabla h(x^{\ast})^{\top}(\balpha^{-1} \hd^{\ast})=0$. 
%Recall again $(\Ya,\za)$ is unique by Proposition\,\ref{prop:analyticunique}.
Comparing this fact and \eqref{eq:0409-2}
to condition\,\eqref{eq:proofanalytic} with $v:=\balpha \hd^{\ast}$, we ensure that  
$(Y_{\ast},z^{\ast})$ is an analytic center of the NSDP at $x^{\ast}$, leading to 
$(\Ya,\za)=(Y_{\ast},z^{\ast})$ due to the uniqueness of analytic center by \cref{prop:analyticunique}.
Finally, recalling that $(Y_{\ast},z^{\ast})$ is an arbitrary accumulation point of
$\left\{(Y_k,z^k)\right\}$, we conclude that the whole sequence $\left\{(Y_k,z^k)\right\}$ converges to $(\Ya,\za)$.  The proof is complete.
\hfill\Halmos\endproof\vspace{0.5em}

Before moving on to the next subsection, we show that $\|\YkUN\|_{\rm F}$ and $\|\YkNN\|_{\rm F}$ are bounded by $\O(\muk)$.
\begin{proposition}\label{prop:ykmuk}
Suppose that Assumption\,\ref{assumA} holds. Then, we have 
\begin{equation*}
\|\YkUN\|_{\rm F}= \O(\mu_k),\  \|\YkNN\|_{\rm F}= \O(\mu_k).
\end{equation*}
\end{proposition}
\proof{\rm Proof.}
Note that $\{(\Yk,\zk)\}$ is convergent by \cref{thm:analytic}
\tred{and thus $\Yk=\O(1)$ and $\zk=\O(1)$.
Moreover, $\DeltaG(\xast;d^k)=\sum_{i=1}^nd^k_i\frac{\partial G(\xast)}{\partial x_i}=\O(\|d^k\|)$.}
Applying Taylor's expansion to $\Gk$ around $\xast$ and using $\Gk\Yk=\mu_kI$ and $\Gast\Yast=\Gast\Ya=O$ give 
\begin{align*}
\mu_k I &= \left(
\Gast + \DeltaG(\xast;d^k)+ \O(\|d^k\|^2)\right)\tred{\Yk}\\
&= \tred{\Gast(\Yk-\Ya) + \Gast\Ya+\DeltaG(\xast;d^k)\Yk+\O(\|d^k\|^2)}\\
&= \Gast(\Yk-\Ya) + \O(\|d^k\|)\\
            &=\Past^{\top}\begin{bmatrix}
O&O \\
            \GNN_{\ast}\YkNU& \GNN_{\ast}\YkNN 
            \end{bmatrix}\Past 
            +\O(\|d^k\|),
\end{align*}
\tred{where
  the third equality holds from $\Gast\Ya=O$, $\DeltaG(\xast;d^k)\Yk=\O(\|d^k\|)$,
  and $O(\|d^k\|^2)=\O(\|d^k\|)$.}
Recall that $\Past$ is an orthogonal matrix. Divide both the sides of the above by $\mu_k$ and drive $k\to\infty$.
From \cref{prop:mukdk}, we obtain
\begin{equation*}
 \frac{\|\GNN_{\ast}\YkNU\|_{\rm F}}{\mu_k} = \O(1),\  \frac{\|\GNN_{\ast}\YkNN\|_{\rm F}}{\mu_k} = \O(1),  
\end{equation*}
which together with $\GNN_{\ast}\in \mathbb{S}^{m-\rast}_{++}$ implies the desired assertions.
\hfill\Halmos\endproof\vspace{0.5em}

\subsection{{\bf Proof of Claim~(ii):}
  \tred{convergence of BKKT points along a specific direction}}
  \label{subsec:claim2}
%% In this subsection, we will prove that a smooth central path $w(\mu)=(x(\mu),Y(\mu),z(\mu))$ that converges to $w^{\ast}$
%% as $\mu\to 0$ exists {\it uniquely} in a sufficiently small vicinity of $w^{\ast}$. 
%% Furthermore, we will prove that
%% $(\dot{x}(\mu),\cdot{z}(\mu))$ is convergent as $\mu$ goes to $0$ and then charatatize its limit specifically.
Let $\{w^k=(\xk,\Yk,\zk)\}\subseteq \W_{++}$ be a sequence of BKKT triplets as described right after \cref{thm:barrier}. From \cref{thm:analytic}, $\{w^k\}$ converges to the KKT triplet $\wast=(\xast,\Ya,\za)$ with the analytic center $(\Ya,\za)$.
In this subsection, we study how $d^k/\muk$ behaves asymptotically, wherein $d^k$ is defined in \eqref{eq:dk}.

We begin by considering the following equation-system that comes from the BKKT conditions of the symmetric form:
\begin{equation}
\nabla_xL(w)=0,\ G(x)Y+YG(x)=2\mu I,\ h(x) = 0,\label{eq:BKKTsym}
\end{equation}
where $w=(x,Y,z)\in \W_{++}$. Suppose at this moment\footnote{In Theorem\,\ref{lastthm}, this assumption will be verified.}
that there exists a smooth function $w(\cdot): (0,\bar{\mu}]\to \W_{++}$
  with some $\bar{\mu}>0$ such that $w(\mu)$ is a BKKT triplet for each $\mu\in (0,\bar{\mu}]$ and we stand at $w=w(\mu)$.
 Differentiating equations\,\eqref{eq:BKKTsym} with respect to $\mu$ results in  
\begin{eqnarray}
&\nabla^2_{xx}L(w)\dot{x}- \mathcal{J}G(x)^{\ast}\dot{Y}+ \nabla h(x)\dot{z}=0,\label{eqn:diff1} \\
&\mathcal{L}_{G(x)}\dot{Y} + \mathcal{L}_Y\Delta G(x;\dot{x}) = 2I,\label{eqn:diff2}\\
&\nabla h(x)^{\top}\dot{x}= 0.\label{eqn:diff3}
\end{eqnarray}
 As for the definition of $\mathcal{L}_{(\cdot)}(\cdot)$, refer to the section of notations.
For later use, in terms of the matrix function
\begin{equation}
\JPHI(w):=\begin{bmatrix}
\nabla^2_{xx}L(w)&-\calJ G(x)^{\ast} &\nabla h(x)\\
\mathcal{L}_{Y}\calG_1(x)\cdots\mathcal{L}_{Y}\calG_n(x)& \mathcal{L}_{G(x)} &0\\
\nabla h(x)^{\top}&0 & 0 
\end{bmatrix},
\label{eq:jphi}
\end{equation}
we express the above equation-system \eqref{eqn:diff1}-\eqref{eqn:diff3} as
\begin{equation}
  \JPHI(w)\dot{w}=
  %= [0,2I,0]^{\top}.
\begin{bmatrix}0\\2I\\0\end{bmatrix}.
    \label{eq:JPHI}
 \end{equation}
{\begin{remark}\label{rem4}
 The Newton equation to the BKKT system\,\eqref{eq:BKKTsym} is expressed as
  \begin{equation}
  \JPHI(w)\Delta w= \begin{bmatrix}
    &-\nabla_xL(w)\\
    &2\mu I-\mathcal{L}_{G(x)}Y\\
    &-h(x)
    \end{bmatrix}.\label{eq:JPHINewton}
  \end{equation}
  This is often solved in the primal-dual interior point method for the NSDP\,\citep{yamashita2012primal,yamashita2021primal,yamakawa1,yamakawa2}.
% where $dw:=(dx,dY,dz)^{\top}\in \R^n\times \mathbb{S}^m\times \R^s$. 
\end{remark}}
Now, relevant to equation\,\eqref{eq:JPHI}, we consider the following equations defined at the KKT triplet $\wa=(\xast,\Ya,\za)$:
\begin{eqnarray}
&\UR^{\top}\left(\nabla^2_{xx}L(\wa)\Delta x- \mathcal{J}G(\xast)^{\ast}\Delta Y\right)=0,\label{eqn:diff4} \\
%&\nabla^2_{xx}L(\wa)\Delta x- \mathcal{J}G(\xast)^{\ast}\Delta Y+ \nabla h(\xast)\Delta z=0,\label{eqn:diff4} \\
&\mathcal{L}_{\Gast}\Delta Y + \mathcal{L}_{\Ya}\Delta G(\xast;\Delta x) = 2I,\label{eqn:diff5}\\
&\nabla h(\xast)^{\top}\Delta x= 0,\label{eqn:diff6}
\end{eqnarray}
where 
$\UR$ denotes an arbitrary matrix whose columns form an \tred{orthonormal} basis of 
the subspace 
\begin{equation}
\mathUast:=\{d\in \R^n\mid \DeltaGUU(\xast;d)=O, \nabla h(\xast)^{\top}d=0\}\label{eq:mathUast}
\end{equation}
and  we can write $\UR\in \R^{n\times \past}$ by letting $\past$ be the dimension of $\mathUast$.
Notice that the above equations \eqref{eqn:diff4}-\eqref{eqn:diff6} are derived by
changing the variables in \eqref{eqn:diff1}-\eqref{eqn:diff3}, pre-multiplying \eqref{eqn:diff1} by the matrix $\UR^{\top}$, and using the relation 
$\nabla h(\xast)^{\top}\UR=0$.
%Moreover, note that they are well-defined regardless of 
%For the sake of the subsequent analysis, we decompose a vector $\Delta x\in \R^n$ into orthogonal component vectors as follows: 
%\begin{equation}
%\Delta x = \UR \eta^1 + \UN\eta^2,
%\label{eq:decompose}
%\end{equation}
%where $(\eta^1,\eta^2)\in \R^{\past}\times \R^{n-\past}$ and $\UN\in \R^{n\times (n-\past)}$ is a matrix whose columns form an normal orthogonal basis of the orthogonal complement subspace for $\mathUast$.
The following proposition holds as to the solution set of equations\,\eqref{eqn:diff4}-\eqref{eqn:diff6}:
\begin{proposition}\label{lem:xi}
Suppose that Assumption\,\ref{assumA} holds.
Let $$
S:=\left\{(\Delta x,\Delta Y)\in \R^n\times \mathbb{S}^m \mbox{: solution to \eqref{eqn:diff4}-\eqref{eqn:diff6}
}\right\}.$$
If $S\neq \emptyset$, then the following properties hold:
\begin{enumerate}
%\item $S\neq \emptyset$; [[\tred{Proved Later}]]
\item\label{lem:xi-1} $\Delta x$-component in $S$ is unique, written as $\xiast\in \R^n$; 
\item\label{lem:xi-2} $\Delta\YNN= (\GNN_{\ast})^{-1}$, 
$\DeltaGUU(\xast;\xiast)=(\YUU_{\rm a})^{-1}$,\ 
%\YUU_{\rm a}\DeltaGUN(\xast;\xi)+\Delta\YUN\GNN_{\ast}=O$;
$\Delta\YUN = -\YUU_{\rm a}\DeltaGUN(\xast;\xiast)(\GNN_{\ast})^{-1}$.
%\item\label{lem:xi-2} $\DeltaYUU$ is free in $\rm Ker}$
%\todo{Proved Later}
\end{enumerate}
\end{proposition}
\proof{\rm Proof.}
See Appendix\,\ref{app:lem:xi}.
%Since $(\Delta Y,\Delta z)$ is uniquely determined for a given $\Delta x$ 
%through \eqref{al:0429-2} and \eqref{al:0429-3} together with the full column rank of $\nabla h(x^{\ast})$, 
%t suffices to show that $\Delta x$-component is not empty in $S$ so as to prove $S\neq \emptyset$.
\hfill\Halmos\endproof\vspace{0.5em}
%\begin{remark}
%In nonlinear optimization, a similar equation-system is considered 
%Unlike nonlinear equation, our analysis  
%\end{remark}

The following theorem shows that the limit of $d^k/\mu_k$ is actually equal to the direction $\xiast$, which is defined in the above proposition. 
\begin{theorem}\label{prop:xiast}
Suppose that Assumption\,\ref{assumA} holds.
Let $d^k$ be the vector defined in \eqref{eq:dk} and $\xiast$ be the one defined in \cref{lem:xi}. 
Then, we have
$$
\lim_{k\to\infty}\frac{d^k}{\mu_k}=\xiast.
$$
In particular, $\xiast\neq 0$ \tred{and $\nabla h(\xast)^{\top}{\xiast}=0$}.
\end{theorem}
\proof{\rm Proof.}
\tred{First, recall $\lim_{k\to\infty}\dk=0$.}
From \cref{prop:mukdk}, $\{{d^k}/{\mu_k}\}$ is bounded.
Let $\tilde{\xi}\in \R^n$ be an arbitrary accumulation point of $\{{d^k}/{\mu_k}\}$.
In order to prove the assertion, 
it suffices to show that $\tilde{\xi}$ is a $\Delta x$-component of the solution set of equations\,\eqref{eqn:diff4}-\eqref{eqn:diff6} because of item~\ref{lem:xi-1} of \cref{lem:xi}.
%Note that by recalling $\DeltaGUU(\xast;\UR^i)=O$ for each $i=1,2,\ldots,\past$, it holds that  
%\begin{align}
%\left(\UR^{\top}\JGast^{\ast}(\Yk-\Ya)\right)_i &= \begin{bmatrix}
%                                                             O & \DeltaGUN(\xast;\UR^i)\\
%                                                         \DeltaGNU(\xast;\UR^i)   & \DeltaGNN(\xast;\UR^i)
%                                                    \end{bmatrix}\bullet
%                                                     \begin{bmatrix}
%                                                             \YkUU-\YaUU & \YkUN\\
%                                                         \YkNU & \YkNN
%                                                    \end{bmatrix}\notag \\ 
%                                                    &={\rm Tr}\left(
%                                                              \DeltaGUN(\xast;\UR^i)\YkNU+
%                                                              \DeltaGNU(\xast;\UR^i)\YkUN+\DeltaGNN(\xast;\UR^i)\YkNN
%                                                    \right) \notag \\
%                                                    &=\O(\mu_k)    \label{al:0430-1}                                                 
%\end{align}
%where the third equality follows from Proposition\,\ref{prop:ykmuk}.
Altering $\Yk$ as 
$$\widehat{Y}_k:=\Past\begin{bmatrix}
\YaUU& \YkUN\\
\YkNU& \YkNN 
\end{bmatrix}
\Past^{\top}$$ for each $k$, we obtain  
%\begin{equation}
\begin{align}
  \|\widehat{Y}_k-\Ya\|_{\rm F}&=  \left\|
  \begin{bmatrix}
    \Ya^{\UU}-\Ya^{\UU}&\Yk^{\Un}-\Ya^{\Un}\\
    \Yk^{\NU}-\Ya^{\NU}&\Yk^{\NN}-\Ya^{\NN}
  \end{bmatrix}
  \right\|_{\rm F}\notag \\
 &=  \left\|
  \begin{bmatrix}
 O&\Yk^{\Un}\\
    \Yk^{\NU}&\Yk^{\NN}
  \end{bmatrix}
  \right\|_{\rm F}\notag \\
 &=\O(\|\Yk^{\Un}\|_{\rm F}+\|\Yk^{\NN}\|_{\rm F}) \notag \\
 &=\O(\mu_k),\notag 
\end{align}
%\label{eq:yktilde}
%\end{equation}
where the last equality follows from \cref{prop:ykmuk}, and thus 
$\left\{\frac{1}{\mu_k}(\widehat{Y}_k-\Ya)\right\}$ is bounded and has at least one accumulation point, say 
$\Delta Y_{\ast}$.
Without loss of generality, we assume that
\begin{equation}
\lim_{k\to\infty}\frac{1}{\mu_k}(\widehat{Y}_k-\Ya)=\Delta Y_{\ast}.\label{eq:240126}
\end{equation}
    {Note that 
      \begin{align}
        \nabla_xL(w^k)&=\nabla_xL(\xk,\Ya,\za)
        - \mathcal{J}G(\xk)^{\ast}(\Yk-\Ya)+\nabla h(\xk)(\zk-\za)\notag\\
        &=\left(\nabla_xL(\wa)+\nabla^2_{xx}L(\wa)d^k + \O(\|d^k\|^2)\right) - \mathcal{J}G(\xk)^{\ast}(\Yk-\Ya)+\nabla h(\xk)(\zk-\za)\notag\\
        &=\nabla^2_{xx}L(\wa)d^k - \mathcal{J}G(\xk)^{\ast}(\Yk-\Ya)+\nabla h(\xk)(\zk-\za) +\O(\|d^k\|^2)\notag\\
        &=\nabla^2_{xx}L(\wa)d^k - \mathcal{J}G(\xast)^{\ast}(\Yk-\Ya)+\nabla h(\xast)(\zk-\za) +\o(\|d^k\|),\label{al:231111-1}
      %  &=\nabla^2_{xx}L(\wa)d^k - \mathcal{J}G(\xast)^{\ast}(\Yk-\Ya)+\nabla h(\xast)(\zk-\za) +\O(\|d^k\|^2)
      \end{align}
      where the second equality follows from applying Taylor's expansion to
$\nabla_xL(\xk,\Ya,\za)$ around $\xast$ with respect to $x$ and the third one from $\nabla_xL(\wa)=0$.
      Moreover, the last one holds because $\o(\|d^k\|)+\O(\|d^k\|^2)=\o(\|d^k\|)$ and 
      $$
      - \mathcal{J}G(\xk)^{\ast}(\Yk-\Ya)+\nabla h(\xk)(\zk-\za)= - \mathcal{J}G(\xast)^{\ast}(\Yk-\Ya)+\nabla h(\xast)(\zk-\za) +\o(\|d^k\|)
      $$
      follows from Taylor's expansion of $\mathcal{J}G(\xk)$ and $\nabla h(\xk)$ around $\xast$ again and $\lim_{k\to\infty}(\Yk-\Ya,\zk-\za)=(O,0)$.
    }
%By the BKKT conditions and applying Taylor's expansion to $\nabla_xL(w^k)$ around $w^{\ast}$, it holds that
%By the BKKT conditions,
%and applying Taylor's expansion to $\nabla_xL(w^k)$ around \tred{$\wa=(\xast,\Ya,\za)$},
%  \footnote{By Taylor's expansion, we have
%  \tred{$\nabla f(\xk)+\mathcal{J}G(\xk)^{\ast}\Yk+\nabla h(x^k)\zk
%  =\nabla_xL(\wa)+\nabla^2_{xx}L(\wa)d^k+\mathcal{J}G(\xast)^{\ast}(\Yk-\Ya)+\nabla h_i(\xast)(\zk-\za)+\O(\|d^k\|)
%  $.}
%}
Then, it holds that 
\begin{align}
0&=\frac{1}{\mu_k}\UR^{\top}\nabla_xL(w^k)\notag \\
  &= \frac{1}{\mu_k} \UR^{\top}\left(\nabla_{xx}^2L(\wa)d^k - \JGast^{\ast}(\Yk-\Ya) \right)+ \tred{\frac{1}{\mu_k}\o(\|d^k\|)} \notag\\ 
  &=\UR^{\top}\left(\nabla_{xx}^2L(\wa)\frac{d^k}{\mu_k} - \JGast^{\ast}\frac{(\widehat{Y}_k-\Ya)}{\mu_k} \right)+ \tred{\frac{1}{\mu_k}\o(\|d^k\|)}, \notag 
\end{align}
where the first equality follows from 
$\nabla_xL(\wk)=0$, the second one does from
\tred{\eqref{al:231111-1}} and $\UR^{\top}\nabla h(\xast)=0$, and
the third one does from $\UR^{\top}\JGast^{\ast}\Yk=\UR^{\top}\JGast^{\ast}{\widehat{Y}_k}$.
Driving $k\to \infty$ above and \tred{using \cref{prop:mukdk} and \eqref{eq:240126}} imply
\begin{equation*}
\UR^{\top}\left(\nabla_{xx}^2L(\wa)\tilde{\xi} - \JGast^{\ast}\Delta Y_{\ast}\right)=0,%\label{eq:0430-1}
\end{equation*} 
which is nothing but \eqref{eqn:diff4} with $(\Delta x,\Delta Y)=(\widetilde{\xi},\Delta Y_{\ast})$.

Next, by $G_kY_k=\mu_kI$ from the BKKT conditions and also by noting $\Gast \widehat{Y}_k=\Gast\Yk$ together with $\Gast \Ya=O$, there holds that 
\begin{align*}
I&=\frac{1}{\mu_k}\left(\Gk\Yk - \Gast \Ya\right)\\
 &=\frac{1}{\mu_k}\left(\left(\Gast+\Delta G(\xast;d^k) + \O(\|d^k\|^2)\right)\Yk - \Gast \Ya\right)\\
  &=\frac{1}{\mu_k}\left(\Gast(\widehat{Y}_k-\Ya)+\Delta G(\xast;d^k)\Yk\right) + \O(\|d^k\|),
\end{align*}
wherein by driving $k\to \infty$, symmetrizing, and using $\lim_{k\to\infty}\Yk=\Ya$, we gain \eqref{eqn:diff5} with $(\Delta x,\Delta Y)=(\widetilde{\xi},\Delta Y_{\ast})$.
Finally, we can prove \eqref{eqn:diff6} with $(\Delta x,\Delta Y)=(\widetilde{\xi},\Delta Y_{\ast})$ by driving $k$ to $\infty$ 
in the relation
\begin{equation}
  0=\frac{1}{\mu_k}(h(x^k)-h(x^{\ast}))=\nabla h(\xast)^{\top}\frac{d^k}{\mu_k} + \O(\|d^k\|).\label{eq:shinjuku}
\end{equation}
Consequently, $(\widetilde{\xi},\Delta Y_{\ast})$ solves \eqref{eqn:diff4}-\eqref{eqn:diff6}.
Hence, $\widetilde{\xi}=\xiast$, namely, $\lim_{k\to\infty}d^k/\muk=\xiast$ is ensured by using item~\ref{lem:xi-1} of \cref{lem:xi}. The remaining assertions \tred{$\xiast\neq 0$ and $\nabla h(\xast)^{\top}\xiast=0$ follow} immediately
since $\|d^k\|=\Theta(\muk)$ from \cref{prop:mukdk} \tred{and \eqref{eq:shinjuku} holds.}
The proof is complete.
\hfill\Halmos\endproof\vspace{0.5em}

It is worth noting that we have multiple choices for $\{(d^k,\muk)\}$, while    
$\xiast$ is the constant vector that is uniquely determined as a $\Delta x$-component of the solution set to the equation-system\,\eqref{eqn:diff4}-\eqref{eqn:diff6}. Nevertheless, 
according to Theorem\,\ref{prop:xiast}, any $\{d^k/\muk\}$ converges to $\xiast$.

\cref{prop:xiast} yields the following corollary, a clear picture about how $x^k$ approaches $x^{\ast}$. 
\begin{corollary}\label{cor1}
Under Assumption\,\ref{assumA}, 
we obtain $\lim_{k\to\infty}\frac{x^k-\xast}{\|x^k-\xast\|}=\frac{\xiast}{\|\xiast\|}$. This indicates that
$x^k$ approaches $x^{\ast}$ along the direction $-\xiast$ asymptotically. 
\end{corollary}
\proof{\rm Proof.}
From \cref{prop:xiast}, we have
$$\lim_{k\to\infty}
\frac{x^k-\xast}{\|x^k-\xast\|}=\lim_{k\to\infty}
\frac{x^k-\xast}{\muk}\frac{\muk}{\|x^k-\xast\|}=\frac{\xiast}{\|\xiast\|}.$$
The proof is complete.
\hfill\Halmos\endproof\vspace{0.5em}
%\begin{remark}\cref{prop:xiast} along with \cref{prop:mukdk} yields $\lim_{k\to\infty}\frac{x^k-\xast}{\|x^k-\xast\|}=\frac{\xiast}{\|\xiast\|}$, which presents the clear picture that $x^k$ approaches $x^{\ast}$ along the direction $-\xiast$ asympotically.   
%\end{remark} 
%As will be proven, the central path, written $w(\mu)$, is smooth with respect to $\mu>0$ and leads to $\wast=(\xast,\Ya,\za)$ as $\mu\to 0$. 
%We can closely relate $\frac{1}{\muk}d^k$ to the tangential direction $\dot{x}(\mu)$ in the $x$-space.
%Indeed, if $\lim_{\mu\to 0}\dot{x}(\mu)$ exists, it is identical to $\lim_{\mu\to 0}\frac{1}{\mu}(x(\mu)-x(0))$ with $x(0)=\xast$, namely, $\xiast$.
%
%The above proposition leads us to the important geometric insight as to how $x^k$ moves:
%Let $\Ker\DeltaGUU(\xast;\cdot):=\Set{d\in \R^n|\DeltaGUU(\xast;d)=O}$. 
%Then, any vector in $\R^n$ can be represented as the summation of

%\begin{figure}
%\centering
%\includegraphics[scale=0.4,bb=0 0 993 1324]{fig1.png}
%\end{figure}
%where $\calF$ denotes the feasible region of the NSDP.
%For $\beta>0$ and $\mu\ge 0$, we denote 
%$$\calPc(\mu):=\Set{\xast + \mu\xiast + \eta | \eta\in \Ker\DeltaGUU(\xast;\cdot),\ \|\eta\|\le \beta \mu},$$

%\subsection{Existence and uniqueness of smooth central path}
For $\rho\in (0,1)$ and $\mu\ge 0$, define 
$$\calPc(\mu):=\{x\in \R^n \mid \|\xast + \mu\xiast-x\| < \rho \mu\|\xiast\|\}.$$
%Note that $\bigcup_{\mu\ge 0}\calPc(\mu)$ is a closed convex cone witn nonempty interior and admits $\xast$ and $\xiast$ as vertex and axis, respectively. 
%% From \cref{prop:xiast}, $x^k\in \calPc(\muk)$ holds for any $k$ sufficiently large. This fact implies that $\{\xk\}_{k\ge K}\subseteq \bigcup_{\mu\ge 0}\calPc(\mu)$ for a sufficiently large $K$.
%In the next subsection, we study properties of $\calPc(\mu)$ more precisely. 
\subsection{Some properties on $\calPc(\mu)$}\label{subsec:claim3}
%\subsection{Nonsingularity of the matrix $\mathcal{A}(w)$}\label{subsec:claim3}
%In this subsection, we will prove that the matrix $\mathcal{A}(w)$ defined in \eqref{eq:jphi} is nonsingular in a certain region around 
%$\xast$, 
%We begin by studying properties of region near the line emanuating from the KKT point $\xast$ in the direction $\xiast$.
%Let us study properties related to $\calPc(\mu)$ and $\bigcup_{\mu\ge 0}\calPc(\mu)$.
In this subsection, we will present some propositions about properties relevant to $\calPc(\mu)$. 
{The following proposition concerns the existence of a BKKT point in $\calPc(\mu)$.
  %  The following proposition is derived from 
\begin{proposition}\label{prop:bkktcalPc}
    Choose $\rho\in (0,1)$ arbitrarily. 
\begin{enumerate}
\item\label{prop:bkktcalPc:item1}
 There exists some $\bmurho>0$ such that, 
  %Let $\mu>0$ small enough. Then, there exists a BKKT point $x_{\mu}$ with barrier parameter $\mu>0$ in $\calPc(\mu)$.
  for any $0<\mu\le \bmurho$, a BKKT point $x_{\mu}$ with barrier parameter $\mu$ exists in $\calPc(\mu)$, but never  
  on the boundary
  $\bd\calPc(\mu):=\{x\in \R^n \mid \|\xast + \mu\xiast-x\|=\rho \mu\|\xiast\|\}$. 
        %Conversely, any BKKT points in  $\calPc(\mu)$ are ones with barrier parameter $\mu$.
\item\label{prop:bkktcalPc:item2}
  Let $\{\mu_k\}\subseteq \R_{++}$ and $\{x(\muk)\}$ be arbitrary sequences of barrier parameters converging to 0 and corresponding BKKT points converging to the KKT point $\xast$, respectively. Then, $x(\muk)\in\calPc(\muk)$ for any $k$ large enough.  
        %Moreover, there exists no BKKT point $x_{\mu}$ such that $x_{\mu}\in\calPc(\nu)\setminus \calPc(\mu)$ for some $\mu\neq \nu>0$.
  %        For any $\mu>0$ small enough, any BKKT points $x_{\mu}$ with barrier parameter $\mu>0$ exist in $\calPc(\mu)$.
  \end{enumerate}
  \end{proposition}
  \proof{\rm Proof.}
  See Appendix\,\ref{app:prop6}.
%% We show item\,\ref{prop:bkktcalPc:item1}.
%% We first consider the first-half claim.
%% For contradiction, assume that there exists an infinite sequence $\{\muk\}\subseteq \R_{++}$ converging to $0$ such that
%%  $\calPc(\muk)$ does not contain a BKKT point with barrier parameter $\muk$ for each $k$.
%% According to \cref{thm:barrier}, $\{\muk\}$ accompanies a sequence of BKKT points $\{x({\muk})\}$ which converges to the KKT point $\xast$.
%% By the above assumption, $x({\muk})\notin \calPc(\muk)$ for each $k$, implying
%% $\|x({\muk})-\xast-\muk\xiast\|/\muk\ge \rho\|\xiast\|>0$.
%% However, \cref{prop:xiast} implies 
%% \begin{equation}
%%   \|x({\muk})-\xast-\muk\xiast\|=\o(\muk).\label{eq:231225-1}
%% \end{equation}
%% This is a contradiction. Hence, the first-claim claim is obtained.
%% %The second-half one can be also established by deriving a contradiction and using \eqref{eq:231225-1} in a similar way.
%% The second-half one can be also established by deriving a contradiction. Suppose to the contrary that
%% there exist BKKT points $\{x(\muk)\}$ with $x(\muk)\in\bd \calPc(\muk)$. By the definition of $\bd\calPc(\muk)$, we see
%% $\lim_{k\to\infty}x(\muk)=\xast$, thus \eqref{eq:231225-1} holds again. However, this contradicts $(x(\muk)-\xast)/\muk=\rho,\forall k$ from $x(\muk)\in\bd \calPc(\muk)$.  
%% Item\,\ref{prop:bkktcalPc:item2} also follows readily from \eqref{eq:231225-1}.
%% %Next, we show the second-half one. As well as above, suppose that
%% %there exists an infinite sequence $\{\muk\}\subseteq \R_{++}$ converging to $0$ and 
%% %such that 
%% %to derive contradiction.
\hfill\Halmos\endproof\vspace{0.5em}}
%Since the proofs are long, we defer the proof to Appendix for readability. 
\tred{The next proposition shows existence and properties of $\mu G(x)^{-1}$ for $x\in \calPc(\mu)$.}
{\begin{proposition}\label{prop:0513-1}
    %The following properties hold.
    \tred{Choose $\brhoone\in (0,1]$ and $\bmuone>0$ sufficiently small. Then, the following properties hold:}
\begin{enumerate}
\item\label{prop:0513-1-item1-1}
  %There exist some $\bar{\rho}>0$ and $\bar{\mu}>0$ such that 
  $G(x)\in \mathbb{S}^m_{++}$  holds for any $\mu\in (0,\bmuone]$ and $x\in \tred{\cl\calPcbone(\mu)}$,
and
  $\Set{\mu G(x)^{-1}\in \mathbb{S}^m_{++}| x\in \cl\calPcbone(\mu),\mu\in (0,\bmuone]}$ is bounded.
\item\label{prop:0513-1-item2-2}
%For an arbitrary $\varepsilon>0$,
%there exists $\rho>0$ such that the distance between 
%the matrix $\Ya$ and arbitrary accumulation points of $\mu G(x)^{-1}\ (x\in\calPc(\mu))$ is smaller than $\varepsilon$.
  %When
  {For any $(\mu,\rho)\in (0,\bmuone]\times (0,\brhoone]$ and $x\in \tred{\cl\calPc(\mu)}$,
    there exists some $K_1>0$ such that
$\|\mu G(x)^{-1}-\Ya\|_{\rm F}\le K_1(\rho+\mu)$    
%% When we vary $\mu,\rho$ and $x$ with $x\in \calPc(\mu)$ and $(\mu,\rho)\in (0,\bmu]\times (0,\brho]$, 
%%For $\rho>0,\mu>0$, and $x\in\calPc(\mu)$, we have 
%% $x$, $\rho$, and $\mu$ are varied in such a way that $x\in\calPc(\mu)$, $\rho\to 0$, and $\mu\to 0$, we have 
%% $$
%% \lim_{\mu\to 0,\calPc(\mu)\ni x\to \xast}\mu G(x)^{-1}=\Ya,
%% $$
%% we have
%% $\|\mu G(x)^{-1}-\Ya\|_{\rm F}=\O(\rho+\mu)$.
  %and $\|\nabla_xL(w)\|=\O(\rho+\mu)$ hold.
}
\end{enumerate}
\end{proposition}
\proof{\rm Proof.}
See Appendix\,\ref{app:prop0513-1}.
\hfill\Halmos\endproof\vspace{0.5em}
%Let $\bar{\rho}$ and $\bar{\mu}$ be the constants defined in \cref{prop:0513-1}.
%% Let $\bar{\rho}$ be the constant defined in \cref{prop:0513-1}.
%% Recall that $x^k$ is a BKKT point with barrier parameter $\muk$ and converges to the KKT point $\xast$ as $k\to \infty$.
%% By the first-half assertion of item~\ref{prop:0513-1-item1-1} of \cref{prop:0513-1}, it holds that 
%% $$\left(\bigcup_{\mu>0}\calPcb(\mu)\cap \calB_{r}(x^{\ast})\right)\subseteq \Set{x|G(x)\in \mathbb{S}^m_{++}}$$ for $r>0$ sufficiently small, where $\calB_{r}(x^{\ast})$ denotes the closed Euclidean ball centered at $\xast$ with a radius $r$.
%Let $\brhoone$ be the constant defined in \cref{prop:0513-1}.
%Recall that $x^k$ is a BKKT point with barrier parameter $\muk$ and converges to the KKT point $\xast$ as $k\to \infty$
{As $\nabla h(\xast)$ is of full column rank since the MFCQ holds at $\xast$ and $\nabla h$ is continuous,
  there exists some closed ball $\hball\subseteq \R^n$ centered at $\xast$ such that
  \begin{equation}
\xast\in\hball\subseteq \Set{x\in\R^n|\nabla h(x)\mbox{ is of full-column rank}}.\label{eq:231216-3}
  \end{equation}  
%From \eqref{eq:231216-1} and
%  By the definition of $\calPc(\mu)$, we have $\calPc(\mu)\subseteq \mathcal{P}_{\rho^{\prime}}(\mu)$ for $\rho<\rho^{\prime}$.
%  As $\mu>0$ gets smaller, the distance of $\calPc(\mu)$ from $\xast$ gets smaller.
With sufficiently small $\mu$ and $\rho$, $\cl\calPc(\mu)\subseteq \mathcal{B}$ holds,
  %% \begin{equation}
  %%   \calPc(\mu)\subseteq \mathcal{P}_{\rho^{\prime}}(\mu),\ \forall\rho<\rho^{\prime}.
  %%   \label{eq:231216-1}
  %% \end{equation}
which together with the first-half assertion of item~\ref{prop:0513-1-item1-1} of \cref{prop:0513-1} yields 
that, by re-taking smaller $\bmuone$ and $\brhoone$ if necessary,
\begin{equation}
\cl\calPc(\mu)
\subseteq \Set{x\in\R^n|
%  G(x)\in \mathbb{S}^m_{++}, \nabla h(x)\mbox{ is of full-column rank}},\ \forall \rho\in (0,\brho]\label{eq:231202-1}
  G(x)\in \mathbb{S}^m_{++}}\cap\hball,\ \forall (\mu,\rho)\in (0,\bmuone]\times (0,\brhoone].\label{eq:231202-1}
\end{equation}
%with sufficiently small $\bmu>0$ and $\brho>0$.
%for $r>0$ sufficiently small, where $\calB_{r}(x^{\ast})$ denotes the closed Euclidean ball centered at $\xast$ with a radius $r$.
Next, we consider the following conditions: 
\begin{eqnarray}
&x\in \cl\calPc(\mu),\label{eqn:0524-0}\\
  &\|Y-\mu G(x)^{-1}\|_{\rm F}\le \gamma_1\mu,\label{eqn:0524-1}\\
  %&\|Y-\mu G(x)^{-1}\|_{\rm F}\tred{=\o(1)},\label{eqn:0524-1}\\
%\|\nabla_xL(w)\|\le \gamma_2\mu^{\tau}.\label{eqn:0524-2}
&\|z+(\nabla h(x)^{\top}\nabla h(x))^{-1}\nabla h(x)^{\top}
\left(\nabla f(x)
-\calJ G(x)^{\ast}Y\right)\| \le \gamma_2\mu.\label{eqn:0524-2}
%&\|\nabla_xL(w)\| \le \gamma_2\mu,\label{eqn:0524-2}
%&\|\nabla_xL(w)\|\le \gamma_2\mu^{\delta},\label{eqn:0524-2}
\end{eqnarray}
As discussed in subsection\,\ref{sec:unique}, the set of $(x,Y,z)$ which satisfies the above conditions is a neighborhood of the central path leading to the KKT triplet $\wa$.
%The following two propostions give crutial properties on the set of $w=(x,Y,z)$ which satisfies the above conditions.
The following two propositions give crucial properties on this set.
They will play important roles in proving \cref{prop:231219-1} and \cref{prop:0603-1} in subsection\,\ref{subsec:claim4}.
%These conditions are 
\begin{proposition}\label{prop:231211-1}
  Let $\gamma_1,\gamma_2>0$ and choose $(\bmutwo,\brhotwo)\in (0,\bmuone]\times (0,\brhoone]$ sufficiently small.
%When we vary $(\mu,\rho)\in (0,\bmu]\times (0,\brho]$ and $(x,Y,z)\in \R^n\times \mathbb{S}^m\times \R^s$ with
There exists some $K_2>0$ such that 
%% For the same $\mu,\rho$, and $x$ in item~\ref{prop:0513-1-item2-2} of \cref{prop:0513-1} and $Y\in \mathbb{S}^m_{++}$ and $z\in \R^s$
  %% such that $\|Y-\mu G(x)^{-1}\|_{\rm F}=\O(\mu)$ and
  %% $\|z+(\nabla h(x)\nabla h(x)^{\top})^{-1}\nabla h(x)^{\top}\left(\nabla f(x)-\mathcal{J}G(x)^{\ast}Y\right)\|=\O(\mu)$, we have
%%   $$
%% \|x-\xast\|=\O(\mu),\ \|Y-\Ya\|_{\rm F}=\O(\rho+\mu),\ \|z-\za\|=\O(\rho+\mu),\ \|\nabla_xL(w)\|=\O(\rho+\mu).
%%   $$
\begin{equation}
\|x-\xast\|\le K_2\mu,\ \ \max\left(\|Y-\Ya\|_{\rm F},\|z-\za\|,\|\nabla_xL(w)\|\right)\le K_2(\rho+\mu)\label{eq:231216-4}
\end{equation}
for all $(\mu,\rho)\in (0,\bmutwo]\times (0,\brhotwo]$ and $(x,Y,z)\in \R^n\times \mathbb{S}^m\times \R^s$ satisfying \eqref{eqn:0524-0}, \eqref{eqn:0524-1}, and \eqref{eqn:0524-2}.
  \end{proposition}
  \proof{\rm Proof.}
See Appendix\,\ref{app:prop:231211-1}.
\hfill\Halmos\endproof\vspace{0.5em}
  }

%The second proposition will play an important role for proving \cref{prop:0603-1} that will appear later.
\begin{proposition}\label{prop:0531-1}
Suppose that Assumption\,\ref{assumA} holds.
\tred{Let $\gamma_1,\gamma_2>0$. %Choose $\mu\in (0,\bmu]$ and $\rho\in (0,\brho]$ both sufficiently small.
  Choose $(\bmuthree,\brhothree)\in (0,\bmutwo]\times (0,\brhotwo]$ sufficiently small.
  %  Then, for any $w=(x,Y,z)\in \R^n\times \mathbb{S}^m_{++}\times \R^{\dimeq}$ satisfying \eqref{eqn:0524-0},
  For any $(\mu,\rho)\in (0,\bmuthree]\times (0,\brhothree]$ and $w=(x,Y,z)\in \R^n\times \mathbb{S}^m\times \R^s$ satisfying
  \eqref{eqn:0524-0}, \eqref{eqn:0524-1}, and \eqref{eqn:0524-2}, we have}
%\begin{enumerate}
%\item\label{prop7:item1} \tred{$\nabla h(x)$ is of full column rank for any $x\in \calPc(\mu)$.}
%
%\item\label{prop7:item2} 
%
%{\color{red}
%% \begin{eqnarray}
%% &x\in \calPc(\mu),\label{eqn:0524-0}\\
%%   &\|Y-\mu G(x)^{-1}\|_{\rm F}\le \gamma_1\mu,\label{eqn:0524-1}\\
%%   %&\|Y-\mu G(x)^{-1}\|_{\rm F}\tred{=\o(1)},\label{eqn:0524-1}\\
%% %\|\nabla_xL(w)\|\le \gamma_2\mu^{\tau}.\label{eqn:0524-2}
%% &\|z+(\nabla h(x)^{\top}\nabla h(x))^{-1}\nabla h(x)^{\top}
%% \left(\nabla f(x)
%% -\calJ G(x)^{\ast}Y\right)\| \le \gamma_2\mu,\label{eqn:0524-2}
%% %&\|\nabla_xL(w)\| \le \gamma_2\mu,\label{eqn:0524-2}
%% %&\|\nabla_xL(w)\|\le \gamma_2\mu^{\delta},\label{eqn:0524-2}
%% \end{eqnarray}
%}
\begin{equation}
  d^{\top}\nabla^2_{xx}L(w)d + \Delta G(x;d)\bullet \calL_{G(x)}^{-1}\calL_Y\left(\Delta G(x;d)\right) \tred{\ge \frac{\kappa}{2}\|d\|^2},\ \ \forall d\in \R^n\setminus \{0\}: \nabla h(x)^{\top}d=0,  \label{eq:prop:0531}
\end{equation}
\tred{where $\kappa>0$ is the constant defined in \eqref{eq:ssoc}.}
%\item\label{secondclaim}
  %% The matrix $\JPHI(w)$ defined in \eqref{eq:jphi}
  %% is nonsingular.
  %for any $w$.
  %satisfying conditions\,\eqref{eqn:0524-0}--\eqref{eqn:0524-2}.
%\end{enumerate} 
\end{proposition}
\proof{\rm Proof.}
See Appendix\,\ref{app:prop:0531-1}.
\hfill\Halmos\endproof\vspace{0.5em}

\subsection{{\bf Proof of Claim\,(iii):} uniqueness of BKKT point for each barrier parameter}\label{subsec:claim4}
\tred{In this section, in order to derive the smoothness of the central path by means of the classical implicit function theorem,
we first transform NSDP\,\eqref{al:nsdp} into an equivalent problem without equality constraints locally.}
Let $\calM:=\Set{x\in\R^n|h(x)=0}$.
Under the presence of the full column rank of $\nabla h(\xast)$ and twice continuous differentiability of $h$,
there exists an open set 
$U\subseteq \R^{n-\dimeq}$ together
with a $\calC^2$-{\it diffeomorphism} $\impfun:U\to \impfun(U)(\subseteq \R^n)$ such that $\xast\in \impfun(U)\subseteq \calM$.\footnote{
  \tred{More strictly speaking, there exists a $\calC^2$ mapping $\overline{{\impfun}}:U\to \R^n$ such that
    $\overline{{\impfun}}(U)$ and $U$ are diffeomorphic and furthermore $h((\overline{\impfun}(\uv)^{\top},\uv^{\top})^{\top})=0\ (\uv\in U)$ holds by re-ordering the variables in $x$ if necessary. Then, we let $\impfun(\uv):=(\overline{\impfun}(\uv)^{\top},\uv^{\top})^{\top}$.}}
\tred{Then, we can take an open {\it ball} $\UV$ such that $\cl\UV\subsetneq U$,
  $\xast\in \impfun(\UV)$, and $\nabla h(\impfun(\uv))$ is of full column rank for all $\uv\in\cl\UV$.
  %$\cl X$ denotes the closure of a set $X$.
  Thus, there exists $(\nabla h(x)^{\top}\nabla h(x))^{-1}$ with $x=\impfun(v)$ for all $v\in\cl V$.
}  
%% $\UV \subseteq \R^{n-\dimeq}$ together
%% with a $\calC^2$-{\it diffeomorphism} $\impfun:\UV\to \impfun(\UV)(\subseteq \R^n)$ such that $\xast\in \impfun(\UV)\subseteq \calM$.\footnote{
%%   \tred{More strictly speaking, there exists a $\calC^2$ mapping $\overline{{\impfun}}:\UV\to \R^n$ such that
%%     $\overline{{\impfun}}(\UV)$ and $\UV$ are diffeomorphic and furthermore $h(
%% (\overline{\impfun}(\uv)^{\top},\uv^{\top})^{\top}
%%     )=0\ (\uv\in \UV)$ holds by re-ordering the variables in $x$ if necessary. Then, we let $\impfun(\uv):=(\overline{\impfun}(\uv)^{\top},\uv^{\top})^{\top}$.}
%% }\tred{Here, we assume that $\nabla h(\impfun(\uv))$ is of full column rank for all $\uv$ in the closure of $V$ by taking $V$ smaller if necessary.}o

Let us give some relevant properties of $\Phi$ for later use. 
Since $\UV$ is bounded and $\impfun$ is smooth \tred{on the open set $U(\supsetneq \cl\UV)$}, there exists a \tred{Lipshitz} constant $M_1>0$ such that 
\begin{equation}
\|\impfun(u)-\impfun(v)\|\le M_1\|u-v\|,\ \ \forall u,v\in \tred{\cl\UV}. \label{eq:0608-1}
\end{equation}
Moreover, 
by noting that $\impfun$ is a diffeomorphism, there exists $M_2>0$ such that
\begin{equation*}
\|\impfun^{-1}(x)-\impfun^{-1}(y)\|\le M_2\|x-y\|,\ \ \forall x,y\in \impfun(\tred{\cl\UV}). \label{eq:0608-5}
\end{equation*}
    {\color{red}
      %% Since $h(\impfun(\uv))=0\ (\uv\in V)$, by differentiation with respect to $\uv$, we have
      %% \begin{eqnarray}
      %% &\nabla \impfun(\uv)\nabla_xh_i(\impfun(\uv))=0,&\label{eq:231210-2-1}\\
      %% &\sum_{j=1}^n\frac{\partial h_i(\impfun(\uv))}{\partial x_j}\nabla^2 \impfun_j(\uv)+\nabla \impfun(\uv)\nabla_{xx}^2h_i(\impfun(\uv))\nabla \impfun(\uv)^{\top}=O&\label{eq:231210-2-2}
      %% \end{eqnarray}
      %% for each $i=1,2,\ldots,s$, where $\impfun_j(\uv)$ is the $i$-th element of $\impfun(\uv)\in \R^n$.
      %% Given $z\in \R^s$, \eqref{eq:231210-2-2} yields 
      %% \begin{equation}
      %%   \sum_{j=1}^n\left(\sum_{i=1}^sz_i\frac{\partial h_i(\impfun(\uv))}{\partial x_j}\right)\nabla^2 \impfun_j(\uv)+\nabla \impfun(\uv)\left(\sum_{i=1}^sz_i\nabla_{xx}^2h_i(\impfun(\uv))\right)\nabla \impfun(\uv)^{\top}=O. 
      %%   \label{eq:231210-3}
      %% \end{equation}
    }
{Since $h(\impfun(\uv))=0\ (\uv\in V)$, by differentiation with respect to $\uv$, we have
      \begin{eqnarray}
      &\nabla \impfun(\uv)\nabla_xh_i(\impfun(\uv))=0,&\label{eqn:231213-2}\\
      &\sum_{j=1}^n\frac{\partial h_i(\impfun(\uv))}{\partial x_j}\nabla^2 \impfun_j(\uv)+\nabla \impfun(\uv)\nabla_{xx}^2h_i(\impfun(\uv))\nabla \impfun(\uv)^{\top}=O&\label{eq:231210-2-2}
      \end{eqnarray}
      for each $i=1,2,\ldots,s$, where $\impfun_i(\uv)$ stands for the $i$-th element of $\impfun(\uv)\in \R^n$.
      Note that $\nabla\impfun(\uv)^{\top}$ is of full column rank for any $\uv\in U$ because $\impfun$ is a diffeomorphism on $U$ and $\cl \UV$ is bounded by definition. From this fact, there exist some $M_3,M_4>0$ such that
\begin{equation}
  M_3\le \|\nabla\impfun(\uv)^{\top}y\|\le M_4,\ \forall (\uv,y)\in \cl\UV\times \R^{n-s}:\|y\|=1. \label{eq:231217-1}
\end{equation}
%for all $\uv\in \cl \UV$ and $y\in \R^{n-s}$ with $\|y\|=1$.
Since $\cl\UV$ is bounded and $\nabla^2\impfun$ is continuous on $\cl\UV$, there exists some $M_5>0$ such that 
\begin{equation}
\max_{i=1,2,\ldots,s}\|\nabla^2\impfun_i(\uv)\|_{\rm F}\le M_5,\ \forall v\in\cl\UV. \label{eq:231219-2}
\end{equation}
%Therefore, $\nabla\impfun(\uv)=0$
%\begin{equation}
}    
Let $\dist(x,\calM):=\min_{y\in \calM}\|x-y\|$
for $x\in \R^n$.
The following lemma holds.
\begin{lemma}\label{lem:0314-1}
It holds that 
\begin{equation}
\dist(\xmudel,\calM)=\O(\mu^2),\label{eq:0608-2} 
\end{equation}
where
$
\xmudel:=\xast+\mu\xiast\ \ \ (\mu\ge 0).
$
\end{lemma}
\proof{\rm Proof.}
See Appendix\,\ref{appA-2}.
\hfill\Halmos\endproof\vspace{0.5em}
In terms of $\impfun$, NSDP\,\eqref{al:nsdp} is reformulated as the following problem 
without equality constraints locally around $\uast:=\impfun^{-1}(\xast)\in \UV$: 
\begin{equation*}
\min_{\uv\in \UV}f(\impfun(\uv))\ \ \mbox{s.t. }\ \ G(\impfun(\uv))\in \mathbb{S}^m_{+}.
\end{equation*}
Accordingly, we obtain the following barrier penalized problem for each $\mu>0$:  
\begin{equation}
%\min_{\uv\in \UV}\Psi_{\mu}(\uv):=f(\impfun(\uv))-\mu\log\det G(\impfun(\uv))\ \ \mbox{s.t. }\ \ G(\impfun(\uv))\in \mathbb{S}^m_{++}. \label{Pmu}
\min_{\uv\in \UV}\Psi_{\mu}(\uv):=\tred{\psi_{\mu}(\impfun(\uv))}\ \ \mbox{s.t. }\ \ G(\impfun(\uv))\in \mathbb{S}^m_{++},\label{Pmu}
\end{equation}
where $\psi_{\mu}(x)=f(x)-\mu\log\det G(x)$ as defined in \eqref{eq:231213-4}. 
{\cred The gradient and Hessian of $\Psi_{\mu}$ are expressed as 
      \begin{eqnarray}
        &\nabla\Psi_{\mu}(\uv)=\nabla\impfun(\uv)\nabla_x\psi_{\mu}(\impfun(\uv)),&
        \label{al:231211-1}\\
        &\nabla^2\Psi_{\mu}(\uv)=
        \sum_{j=1}^n\frac{\partial \psi_{\mu}(\impfun(\uv))}{\partial x_j}\nabla^2 \impfun_j(\uv)+\nabla \impfun(\uv)\nabla_{xx}^2\psi_{\mu}(\impfun(\uv))\nabla \impfun(\uv)^{\top},&
        \label{al:231211-2}
      \end{eqnarray}
      respectively.
From \eqref{eqn:231213-2} and the fact that $\nabla h_1(x),\nabla h_2(x),\cdots,\nabla h_s(x)$ are linearly independent at $x=\impfun(v)$ and
the dimension of $\Ker\nabla\impfun(\uv)$ is $s$,
$\nabla h_1(x),\ldots,\nabla h_s(x)$ form a basis of $\Ker\nabla\impfun(\uv)$.
From this fact and equation\,\eqref{al:231211-1}, we see that 
%\begin{equation}
$\nabla\Psi_{\mu}(\uv)=0$ if and only if
$\nabla_x\psi_{\mu}(\impfun(\uv))=\nabla f(\impfun(\uv))-\mu\mathcal{J}G(\impfun(\uv))G(\impfun(\uv))^{-1}\in {\rm Im}\nabla h(x)$. Namely,
\begin{equation}
\nabla \Psi_{\mu}(\uv)=0\Leftrightarrow\mbox{$\impfun(\uv)$ is a BKKT point with barrier parameter $\mu$.}
  \label{eq:231221-1}
\end{equation}
}
%Here, note that
%%% \begin{align}
%%%   &\inte\left(\impfun^{-1}\left(\impfun(U)\cap \calPc(\mu)\right) \right)\neq \emptyset;\label{eq:0603-1}\\
%%%   &\inte\left(\impfun^{-1}\left(\impfun(U)\cap \calPc(\mu)\right) \right)\mbox{ is connected}\label{eq:0603-2}
%%% \end{align}
%\begin{align}
%\inte\left(\impfun^{-1}\left(\impfun(U)\cap \calPc(\mu)\right) \right)\mbox{ is nonempty and connected}\label{eq:0603-1}
%\end{align}
%for any $\mu>0$ small enough and a fixed $\rho>0$,
%where ``$\inte$'' denotes the topological interior.\footnote{We can show \eqref{eq:0603-1} with the facts that
%  \begin{itemize}
%  \item $\impfun$ is homeomorphic and hence $\impfun$ and $\impfun^{-1}$ preserve the openness and connectedness of a set;
%  \item for any $\mu>0$ small enough, there exists a BKKT point $x$, thus $x\in \calM$, such that $x\in \inte\calPc(\mu)$ by virtue of Proposition\,\ref{prop:xiast}.
%\end{itemize}}
%\tred{[[
%$\xast$に収束するBKKT点の列があるならば, それは
%十分小さい$\mu$に対しては, $\capcalPc$に入らなければならないことにどこかで言及すべし. 
%]]}
%he following lemma guarantees that such $\wk$ is uniquely determined for sufficiently large $k$.  
{\cred
  \begin{proposition}\label{prop:231219-1}
     Suppose that Assumption\,\ref{assumA} holds.
     Choose $(\brhofour,\bmufour)\in (0,\brhothree] \times (0,\bmuthree]$ sufficiently small, where $\brhothree$ and $\bmuthree$ are the constants
  defined in \cref{prop:0531-1}.
  Then, for any $(\rho,\mu)\in (0,\brhofour]\times (0,\bmufour]$ and $y\in \R^{n-s}$ with $\|y\|=1$, we have
  \begin{equation}
    y^{\top}\nabla^2\Psi_{\mu}(v)y\ge \frac{\kappa M_3^2}{4},\ \forall v\in V\cap \impfun^{-1}(\cl\calPc(\mu)),
    \label{eq:231219-3}
  \end{equation}
  where $\kappa$ and $M_3$ are defined in \eqref{eq:ssoc} and \eqref{eq:231217-1}, respectively.
\end{proposition}
\proof{\rm Proof.}
See Appendix\,\ref{app:prop:231219-1}.
$\hfill\Box$
}

  \begin{theorem}\label{prop:0603-1}
     Suppose that Assumption\,\ref{assumA} holds.
  %\tred{Take $\mu\in (0,\bmu]$ and $\rho\in [0,\brho]$ both sufficiently small.
     \tred{
%For $\brhofive\in (0,\brhofour]$ sufficiently small,
Choose $\rho\in (0,\brhofour]$ arbitrarily, where $\brhofour$ is the constant defined in \cref{prop:231219-1}. 
%     For any $\mu\in (0,\bmu]$,
 Then, there exists some $\bmurho\in (0,\bmufour]$
such that, for any $\mu\in (0,\bmurho]$
a unique BKKT point $\xmu$ with barrier parameter $\mu$ exists in $\impfun(\UV)\cap \calPc(\mu)$.
%% \tred{Choose $\bmu,\brho>0$ sufficiently small. For any $\mu\in (0,\bmu]$ and $\rho\in (0,\brho]$,
%%     there exists a unique BKKT point $\xmu$ with barrier parameter $\mu$ in
%% $\impfun(\UV)\cap \calPc(\mu)$}.
%$\impfun(\UV)\cap \capcalPcb$. 
  In particular, $\impfun^{-1}(\xmu)$ is a unique solution of the equation $\nabla\Psi_{\mu}(v)=0$ in the open set
  $\UV\cap \impfun^{-1}(\calPc(\mu))$.
  Moreover, $\nabla^2\Psi_{\mu}(\impfun^{-1}(\xmu))$ is positive definite, thus $\impfun^{-1}(\xmu)$ is a strict local optimum of problem\,\eqref{Pmu}.}
\end{theorem}
\proof{\rm Proof.}
%See Appendix\,\ref{app:prop:0603-1}.
      {\cred
       % We first fix $\rho$, for which \cref{}
        For a fixed $\rho\in (0,\brhofour]$, according to item\,\ref{prop:bkktcalPc:item1} of \cref{prop:bkktcalPc}, 
$x(\mu)\in \calPc(\mu)$ holds for any sufficiently small $\mu>0$.
Since $h(x(\mu))=0$ holds as $\xmu$ is a BKKT point,
we have $\xmu \in \impfun(\UV)$.
Then, noting $\lim_{\mu\to 0}\xmu=\xast\in \impfun(V)$, we ensure 
$\xmu\in \impfun(\UV)\cap \calPc(\mu)$ for any $\mu$ small enough.

Next, by deriving a contradiction,
we prove that such $x(\mu)$ is unique in $\impfun(V)\cap \calPc(\mu)$ for any sufficiently small $\mu$.}
Assume to the contrary that there exists an infinite sequence \tred{$\{\mul\}\subseteq (0,\bmufour]$} which 
converges to 0 and moreover accompanies two sequences $\{\xell\},\{\txell\}\subseteq
%\inte\left(\impfun(\UV)\cap \capcalPcbl\right)$
\tred{\impfun(\UV)\cap \calPc(\mul)}$
such that for each $\ell$,
  $\xell\neq \txell$, but $\xell$ and $\txell$ are both BKKT points with barrier parameter $\mul$.
Hence, 
$$\uell:=\impfun^{-1}(\xell),\ \tuell:=\impfun^{-1}(\txell)$$ exist in $\UV$, and in addition \eqref{eq:231221-1} yields 
  \begin{equation}
  \nabla\Psi_{\mul}(\uell)=\nabla\Psi_{\mul}(\tuell)=0,\label{eq:0604-1}
  \end{equation}
  where $\Psi_{\mu}$ is defined in \eqref{Pmu}.
%Choose $\bar{\rho}>0$ to be sufficiently small. 
  %Since the preceeding propositions in this subsection are true for both the sequences,
  %By taking the preceeding argument into consideration, 
Moreover, we have 
$$
\xell,\txell\in \impfun(\UV)\cap \calPc(\mul)
$$ for sufficiently large $\ell$. 
  Let $\bxmudell\in \argmin_{y\in \calM}\|\xmudell-y\|$, where $\chx(\cdot)$ is defined in Lemma\,\ref{lem:0314-1}, and let 
\begin{align}
%&\bxmudell\in \argmin_{y\in \calM}\|\xmudell-y\|, \notag \\
\umuell:=\impfun^{-1}(\bxmudell),\ \Uell:=\Set{\uv\in \UV|\|\uv-\umuell\|<\frac{\rho\mul}{2M_1}}\label{al:0508-1}
\end{align}
for each $\ell$.
  Then, by \eqref{eq:0608-2} in \cref{lem:0314-1},
    %Taylor's expansion and $\nabla h(\xast)^{\top}\xiast = 0 = h(\xast)$ we have $h(\xmudell)=\O(\mul^2)$ and thus by  
  \begin{equation}
  \dist\left(\xmudell,\calM\right)=\left\|\bxmudell-\xmudell\right\|=\O(\mul^2).\label{eq:0608-3}
  \end{equation}
For any $\uv\in \Uell$, 
\begin{align}
\|\impfun(\uv)-\xmudell\|&\le \|\impfun(\uv)-\bxmudell\| + \|\bxmudell-\xmudell\|\notag \\
                             & \le \|\impfun(\uv)-\impfun(\umuell)\| + \O(\mul^2) \notag \\ 
                             & \le  M_1\|\uv-\umuell\|  + \O(\mul^2) \notag \\ 
                             & \le  \frac{\rho\mul}{2}+ \O(\mul^2),\notag  
\end{align}
where the second inequality follows from \eqref{al:0508-1} and \eqref{eq:0608-3}, the third 
from \eqref{eq:0608-1}, and the fourth from \eqref{al:0508-1} and $\uv\in\Uell$.  
Hence, it holds that 
$
\|\impfun(\uv)-\xmudel\|\le \rho\mul,\ \forall \uv\in \Uell
$
for sufficiently large $\ell$, which yields that 
\begin{equation}
\impfun(\Uell)\subseteq \calPc(\mul).\label{eq:0508-2}
\end{equation}
Furthermore, since both $\xell$ and $\txell$ converge to $\xast$ as $\ell$ tends to $\infty$ because of
$\xell,\txell\in\calPc(\mul)$, \cref{prop:xiast} implies  
\begin{align}
\|\xmudell-\xell\|=\|\xast+\mul\xiast-\xell\|
                         =\mul{\left\|\xiast-\frac{\xell-\xast}{\mul}\right\|}
                         =\o(\mul),
\end{align}
and also $\|\xmudell-\txell\|=\o(\mul)$ in a similar way. 
These relations along with the triangle inequality 
and \eqref{eq:0608-3} yield 
$$\max\left(\|\bxmudell-\xell\|,\|\bxmudell-\txell\|\right)=\o(\mul),$$ which implies 
that for sufficiently large $\ell\ge 0$, 
$\max\left(
\|\bxmudell-\xell\|,
\|\bxmudell-\txell\|\right)\le \frac{\mul\rho}{4M_1M_2}$, thus
\begin{align*}
\|\umuell-\uell\|=\|\impfun^{-1}(\bxmudell)-\impfun^{-1}(\xell)\|
                                        \le M_2\|\bxmudell-\xell\| 
                                        \le  \frac{\rho\mul}{4M_1}. 
\end{align*}
Therefore, we gain $\tred{\uell\in\Uell}$ for $\ell$ large enough.  In a similar way, we can show $\tred{\tuell\in\Uell}$.
In short, from the above arguments we obtain that  
\begin{equation}
\tred{\{\uell,\tuell\}\subseteq \Uell}\label{eq:0608-4}
\end{equation}
for sufficiently large $\ell$.

%\eqref{al:231213-1} from \cref{lem:231213-5}
%Note \eqref{eq:0508-2}.

{\cred
%We next show that for $\ell$ large enough \eqref{Pmu} with $V:=\Uell$ is strongly convex by showing that $\nabla^2\Psi_{\mu}$ is positive definite
%on $\Uell$.
%Pick $\uv\in \Uell$ arbitrarily.
From \eqref{eq:0508-2} together with $\Uell\subseteq \UV$,
$\Uell\subseteq \UV\cap \impfun^{-1}(\calPc(\mul))$ holds. 
Then, according to \cref{prop:231219-1}, when $\rho\le \brhofour$ and $\ell$ is so large that $\mul\le \bmufour$,  
$\nabla^2\Psi_{\mul}(v)$ is positive definite for all $\uv \in \Uell\subseteq \UV\cap \impfun^{-1}(\calPc(\mul))$.
%\eqref{eq:231202-1} implies $G(x)\in \mathbb{S}^m_{++}$
%$\nabla^2\Psi_{\mul}(v)$ is positive definite for any $\uv\in\Uell$
%when $\ell$ sufficiently large,
Thus,
%by re-taking $\brho(>\rho)$ smaller if necessary,
problem\,\eqref{Pmu} with $V$ replaced by $\Uell$ can be viewed as a strongly convex problem for any large $\ell$.
}
 Therefore, a point \tred{$\uv\in \Uell$} which fulfills $\nabla \Psi_{\mul}(\uv)=0$ must be unique, which together with \eqref{eq:0604-1} and \eqref{eq:0608-4} implies $\tuell=\uell$. 
This gives $\xell=\txell$, a contradiction. 
With this, we ensure that by setting $\mu$ to be small enough,
a BKKT point $\xmu$ exists uniquely in $\impfun(\UV)\cap \calPc(\mu)$.
%\tred{Since $\impfun$ is a bijective mapping from $V$ to $\impfun(V)$, }
Moreover, we also see $\impfun(\xmu)$ is a unique local optimum of \eqref{Pmu} in $\UV\cap \impfun^{-1}(\calPc(\mu))$.
\tred{The positive definiteness of $\nabla^2\Psi_{\mu}(\impfun^{-1}(\xmu))$ is clear from \cref{prop:231219-1} along with
$\impfun^{-1}(\xmu)\in\UV\cap\impfun^{-1}(\calPc(\mu))$.}
We thus obtain the desired assertion.
\hfill\Halmos\endproof\vspace{0.5em}

\subsection{Main claim I: existence and uniqueness of central path}\label{sec:claim5}
%% Recall that $\wa = (\xast,\Ya,\za)$ as defined in \eqref{eq:analyticKKT}, and also recall the definitions of 
%% the sequence $\{x^k\}$ and the barrier parameter sequence $\{\muk\}$ organized in subsections\,\ref{subsec:claim1} and \ref{subsec:claim2}. 
%% For each $k$, $x^k$ is a BKKT point with barrier parameter $\mu_k$ and $\lim_{k\to\infty}x^k=\xast$.
%% Moreover, $\mu_k G(x^k)^{-1},z^k$
\tred{
By \cref{prop:0603-1}, we have ensured the uniqueness and existence of BKKT points around the KKT point $\xast$. 
In the following theorem, we prove that these BKKT points together with the corresponding Lagrange multiplier vectors and matrices form a smooth central path leading to the KKT triplet $\wa=(\xast,\Ya,\za)$, and moreover show such a path is uniquely determined.}
%The following main theorem claims that
%such sequences are on a smooth path leading to $\wa$. 
%Given a barrier parameter sequence $\{\muk\}$ such that $\lim_{k\to\infty}\muk=0$, $\{x^k\}$ is a sequence of BKKT points converging to the KKT point $\xast$. 
%\cref{prop:0603-1} implies that, given $\{\mu_k\}$, such $x^k$ is uniquely determined for each $k$ large enough. 
%Taking this fact into consideration, we establish the following main result.
\begin{theorem}\label{lastthm}
 % Suppose that NC, SSOC, and SC hold at a local opt $x^{\ast}$.
 Suppose that Assumption\,\ref{assumA} holds.
%There exists $(x(\mu),Y(\mu),z(\mu))$ for $\mu>0$ and continuous at $\mu=0$. 
\tred{For a sufficiently small $\bmu>0$, there exists a unique central path $w:(0,\bmu)\to \W_{++}$ such that}
% ``unique'' central path $\wmu:=(\xmu,\Ymu,\zmu)\in \W_{++}$ such that
  \begin{enumerate}
  \item\label{item:thm1}
    \tred{it is smooth} and, 
    %$w(\mu)=(\xmu,\Ymu,\zmu)$ is smooth at any $\mu\in (0,\bmu)$.
    for each $\mu\in (0,\bmu)$, $w(\mu)$ is a BKKT triplet of the NSDP with barrier parameter $\mu$. Moreover,  
      %and $x(\mu)$ is a strict local optimum of problem\,\eqref{Pmu};
  \item\label{item:thm2} $\displaystyle \lim_{\mu\to 0}w(\mu)=\wa$.
%  \item\label{item:thm3} $\displaystyle \lim_{\mu\to 0}\dot{x}(\mu)=\xiast$. 
  \end{enumerate}
%  \tred{In particular, such a path $w(\mu)$ is unique.}
%In particular,  $(Y(0),z(0))$ is the analytic center for NSDP\,\eqref{al:nsdp} at $x^{\ast}$. 
\end{theorem}
\proof{\rm Proof.}
      {\cred
%        The proof is three-fold.
       % By means of the implicit function theorem,
        We first show that there exists a unique smooth path $\uv(\cdot)$ such that $\nabla\Psi_{\mu}(\uv(\mu))=0$ and $\uv(\mu)\in V\cap \impfun^{-1}(\calPc(\mu))$ for each $\mu\in (0,\bmu)$ with some $\bmu$.
 %       Next, we prove the assertion of the theorem by mapping $\uv(\cdot)$ into the set $\mathcal{M}$ via the smooth function $\impfun$.  
        Choose $\rho\in (0,\brho_4]$ arbitrarily and consider $\bmurho>0$ defined as in \cref{prop:0603-1}.
%Choose $\brho,\bmu>0$ small enough and fix $\rho\in (0,\brho]$.
              From \cref{prop:0603-1}, for each $\mu\in (0,\bmurho]$,
                there exists a unique $\overline{\uv}_{\mu}\in V\cap \impfun^{-1}(\calPc(\mu))$ such that $\nabla \Psi_{\mu}(\overline{\uv}_{\mu})=0$.
                Moreover, 
              $\nabla^2\Psi_{\mu}(\overline{\uv}_{\mu})$ is positive definite, thus nonsingular. 
              By applying the implicit function theorem to the equation $\nabla \Psi_{\mu}(\uv)=0$,
%% there exists some $0<l_{\mu}< \mu$, and a unique smooth path
%% $\uv:[l_{\mu},\mu]\to \impfun^{-1}(\calPc(\mu))$
%%   satisfying $\nabla \Psi_{t}(\uv(t))=0$  for each $t\in [l_{\mu},\mu]$
              %%   and $\uv(\mu)=\overline\overline{\uv}_{\mu}$.
there exist some $l_{\mu}\in (0,\mu)$ and $u_{\mu}\in \left(\mu,\min\left(2\mu,\bmurho\right)\right)$ together with a smooth path
$\uvmu:(l_{\mu},u_{\mu})\to\UV\cap\impfun^{-1}(\calPc(\mu))$
satisfying
$\uvmu(\mu)=\overline{\uv}_{\mu}$ and $\nabla \Psi_{t}(\uvmu(t))=0$ for each $t\in (l_{\mu},u_{\mu})$.
  {\cred
    In fact, $\impfun(\uvmu(t))\in\calPc(t)$ holds for each $t\in  (l_{\mu},u_{\mu})$ by re-taking smaller $\bmurho$ if necessary.
For the proof, see the footnote\footnote{\cred Suppose to the contrary that there exists a sequence $\{\mul\}$ converging to $0$
  along with $\{t_{\ell}\}$ such that $t_{\ell}\in  (l_{\mu_{\ell}},u_{\ell})$ and $\impfun(\uv_{\mul}(t_{\ell}))\in \calPc(\mul)\setminus
  \calPc(t_{\ell})$. By noting $l_{\mul}<t_{\ell}<u_{\mul}\le 2\mul$, \eqref{eq:231221-1}, and the fact of $\impfun(\uv_{\mul}(t_{\ell}))\in \calPc(\mul)$, it follows that
  $\lim_{\ell\to\infty}t_{\ell}=0$,
$\impfun(\uv_{\mul}(t_{\ell}))$ is a BKKT point with barrier parameter $t_{\ell}$ and $\lim_{\ell\to\infty}\impfun(\uv_{\mul}(t_{\ell}))=\xast$.
  However, according to item\,\ref{prop:bkktcalPc:item2} of \cref{prop:bkktcalPc}, $\impfun(\uv_{\mul}(t_{\ell}))\in \calPc(t_{\ell})$ must hold for any $\ell$ large enough, a contradiction to the assumption of $\impfun(\uv_{\mul}(t_{\ell}))\in \calPc(\mul)\setminus
  \calPc(t_{\ell})$.}.
%$l_{\mu}$
%Actually, we can set $l_{\bmu}=0$.
%because $\bd\calPc(\mu)$ never contains a BKKT point with $\mu>0$ when $\bmu$ is small due to item\,\ref{prop:bkktcalPc:item1} of \cref{prop:bkktcalPc}.
%Notice that
Thus, due to \cref{prop:0603-1} again,
for each $t\in (l_{\mu},u_{\mu})$, $\uvmu(t)$ is the unique solution to $\nabla\Psi_t(\uv)=0$ in the open set $\UV\cap\impfun^{-1}(\calPc(t))$
and $\nabla^2\Psi_{t}(\uvmu(t))$ is positive definite.
%We can consider the set of paths $\uvmu(\cdot)$ for each $\mu\in (0,\bmurho]$ in the same manner.
Taking this fact into account and connecting the paths $\uvmu(\cdot)$ constructed in the above way for each $\mu\in (0,\bmu)$ with $\bmu:=\bmurho$, we can ensure there exists a unique smooth path $\uv:(0,\bmu)\to\UV$ such that $\nabla\Psi_{\mu}(\uv(\mu))=0$ and  $\uv(\mu)\in V\cap \impfun^{-1}(\calPc(\mu))$ for each $\mu\in (0,\bmu)$.
  }
%This path is defined on some open interval containing $\mu$.
%% there exists some $\delta_{\mu}\in (0,\mu)$ and a smooth path $\uv(\cdot):[\mu-\delta_{\mu},\mu+\delta_{\mu}]\to \UV$
%% such that $\nabla \Psi_{t}(\uv(t))=0$ for $t\in [\mu-\delta_{\mu},\mu+\delta_{\mu}]$.
%% From \eqref{eq:231221-1}, $\impfun(\uv(\mu-\delta_{\mu}))$ is a BKKT point with barrier parameter $\mu-\delta_{\mu}$ and thus
%% $\impfun(\uv(\mu-\delta_{\mu}))\in\calPc(\mu-\delta_{\mu})$
%% by \cref{prop:0603-1} again. 
%and $\impfun(\uv(\mu))=x_{\mu}\in\calPc(\mu)$ for each $\mu\in (0,\bmu)$.

By letting $x(\mu):=\impfun(\uv(\mu))$ for each $\mu\in (0,\bmu)$,
$x(\cdot)$ is smooth on $(0,\bmu)$ as $\uv(\cdot)$ is smooth and  $\impfun:\UV\to \mathcal{M}$ is a diffeomorphism.
Furthermore, $x(\mu)$ is a BKKT point with barrier parameter $\mu$ because of $\nabla\Psi_{\mu}(\impfun(\uv(\mu)))=0$ and \eqref{eq:231221-1}.
Since $\xmu=\impfun(\uv(\mu))\in\calPc(\mu)$, it is clear from the definition of $\calPc(\mu)$ that 
\begin{equation}
  \lim_{\mu\to 0}\xmu=\xast.\label{eq:231220-1}
\end{equation}
From \eqref{eq:231216-3} and \eqref{eq:231202-1},
$G(\xmu)\in \mathbb{S}^m_{++}$ and $\nabla h(\xmu)$ is of full column rank, and therefore
$G(\xmu)^{-1}$ and $(\nabla h(\xmu)^{\top}\nabla h(\xmu))^{-1}$ exist.
  Since $\xmu$ is a BKKT point as shown above, there exists $z_{\mu}\in \R^s$ such that
  $$\nabla_xL(\xmu,\mu G(\xmu)^{-1},z_{\mu})=\nabla f(\xmu)-\mu\mathcal{J}G(\xmu)^{\ast}G(\xmu)^{-1}+\nabla h(\xmu)z_{\mu}=0.$$
  Therefore, by premultiplying this equation with $(\nabla h(\xmu)^{\top}\nabla h(\xmu))^{-1}\nabla h(\xmu)^{\top}$, we obtain
  $z_{\mu}=-(\nabla h(\xmu)^{\top}\nabla h(\xmu))^{-1}\nabla h(\xmu)^{\top}(f(\xmu)-\calJ G(\xmu)^{\ast}Y(\mu))$ for each $\mu\in (0,\bmu)$.
%$\lim_{\mu\to 0}\impfun(\uv(\mu))=\xast$.
  Thus, we can define $w(\cdot):(0,\bmu)\to \W_{++}$ by $Y(\mu):=\mu G(x(\mu))^{-1}$
  and $z(\mu):=z_{\mu}$.
  Since $x(\cdot)$ is smooth, so is $w(\cdot)$, and $w(\mu)$ is a BKKT triplet with \eqref{eq:231220-1}. 
   Recall that \cref{thm:analytic} implies any sequence of BKKT triplets $\{w^k=(\xk,\Yk,\zk)\}$ such that $\lim_{k\to\infty}\xk=\xast$ converges to $\wa$.
   Therefore, $\lim_{\mu\to 0}\wmu=\wa$ follows.
   Finally, since $x(\cdot)$ is uniquely determined as shown above and $Y(\cdot)$ and $z(\cdot)$ are uniquely constructed from $x(\cdot)$,
    we can conclude the uniqueness of $w(\cdot)$. The proof is complete.     \hfill\Halmos\endproof\vspace{0.5em}}
\subsection{\tred{Main claim I\hspace{-0.1em}I: unique solvability of the Newton equation in the primal-dual interior point method}}\label{sec:unique}
\tred{As remarked in \cref{rem4}, $\JPHI(w)$ is the coefficient matrix of the Newton equation to the BKKT system\,\eqref{eq:BKKTsym}.
In the following theorem, $\JPHI(w)$ is shown to be nonsingular near the central path.}
\begin{theorem}\label{lastthm2}
      Let the same assumptions as in \cref{prop:0531-1} hold.
For any $w=(x,Y,z)\in \W$ satisfying \eqref{eqn:0524-0}, \eqref{eqn:0524-1}, and \eqref{eqn:0524-2}, 
the matrix $\JPHI(w)$ defined in \eqref{eq:jphi} is nonsingular.
\end{theorem}
\proof{\rm Proof.}
%The second assertion is easily obtained from the first one. Indeed, for proving this,
We have only to show that $\calA(w)dw=0$ when $x\in\cl\calPc(\mu)$ for $dw:=(dx,dY,dz)^{\top}\in \W$ implies $dw=0$. From $\calA(w)dw=0$, it holds that
\begin{eqnarray}
  &\nabla_{xx}^2L(w)dx-\calJ G(x)^{\ast}dY+\nabla h(x)dz=0,\label{eqn:0528-1}\\
  &\calL_{G(x)}dY+\calL_{Y}\DeltaG(x;dx)=O,\label{eqn:0528-2}\\
  &\nabla h(x)^{\top}dx=0.\label{eqn:0528-3}
\end{eqnarray}
Note that $(G(x),Y)\in \mathbb{S}^m_{++}\times \mathbb{S}^m_{++}$ and $\nabla h(x)$ is of full column rank \tred{due to $x\in\cl\calPc(\mu)$, \eqref{eq:231216-3}, and \eqref{eq:231202-1}.}
Pre-multiplying \eqref{eqn:0528-1} with $dx^{\top}$ and substituting \eqref{eqn:0528-3} and $dY=-\calL_{G(x)}^{-1}\calL_{Y}\DeltaG(x;dx)$ from \eqref{eqn:0528-2} into it, we have
$dx^{\top}\nabla^2_{xx}L(w)dx + \Delta G(x;dx)\bullet \calL_{G(x)}^{-1}\calL_Y\left(\Delta G(x;dx)\right)=0$.
Then, because of \cref{prop:0531-1}, $dx=0$ must hold, which together with \eqref{eqn:0528-2} implies $dY=O$. Moreover, \eqref{eqn:0528-1} and the full column rank of
$\nabla h(x)$ give $dz=0$. Hence, we obtain $dw=0$ and thus the second assertion is obtained.  
 \hfill\Halmos\endproof\vspace{0.5em}
%By virtue of item\,\ref{secondclaim} of \cref{prop:0531-1},
{\cred
  %  An equation-system that admits $\mathcal{A}(w)$ as a coefficient-matrix, such as the one\,\eqref{eqn:diff1}-\eqref{eqn:diff3} and Newton equations relevant to the KKT and BKKT conditions, is ensured to be uniquely solvable as long as $w$ satisfies \eqref{eqn:0524-0}-\eqref{eqn:0524-2}.
  %Near $\wa$,
  The set of $w=(x,Y,z)\in \W_{++}$ which fulfills \eqref{eqn:0524-0}-\eqref{eqn:0524-2}, write $\mathcal{N}\subseteq \W_{++}$, contains
  any BKKT triplets with sufficiently small barrier parameters. Indeed, from item\,\ref{prop:bkktcalPc:item2} of \cref{prop:bkktcalPc},
  \eqref{eqn:0524-0} holds at any BKKT points $\xmu$ with sufficiently small barrier parameter $\mu$, and $\left(\nabla h(\xmu)\nabla h(\xmu)^{\top}\right)^{-1}$ exists from \eqref{eq:231216-3} and \eqref{eq:231202-1}. Moreover, the expressions on the left-hand sides of \eqref{eqn:0524-1} and \eqref{eqn:0524-2} are both equal to 0
  because $\nabla_xL(w)=0$ and $G(x)Y=\mu I$ hold on the central path, and hence \eqref{eqn:0524-1} and \eqref{eqn:0524-2} hold true.
  %  Therefore, $\mathcal{N}$ can be viewed as a neighborhood of the central path.
    Therefore, $\mathcal{N}$ contains the central path.
  Thus, \cref{lastthm2} indicates that the Newton equation\,\eqref{eq:JPHINewton} is uniquely solvable when $w$ is close to the central path.
% Therefore, the set of $w\in \W_{++}$ satisfying \eqref{eqn:0524-0}-\eqref{eqn:0524-2} indicates a region near the central path such that
  This fact would be useful particularly when applying the Newton method in the primal-dual interior point method. 
}

%% {\color{blue}
%%   $$
%%   \Set{w\in\W_{++}| \|\nabla_xL(w)\|\le \gamma_3\mu^2,\|G(x)Y+YG(x)-2\mu I\|_{\rm F}\le\gamma_4\mu^2,
%% \|h(x)\|\le \gamma_5\mu^2.}
%%   $$
%% }
\section{Concluding remarks and future work}\label{sec:con}
In this paper, we have studied properties of a central path for nonlinear semidefinite optimization problems (NSDPs). 
Specifically, we have proven that, under the strict complementarity condition, strong second-order sufficient condition, and Mangasarian-Fromovitz constraint qualification, there exists a smooth central path 
which converges to a KKT triplet with an analytic center. In particular, given a KKT triplet, a central path leading to that KKT triplet is uniquely determined. Unlike the past results concerning the central path for the NSDP, the nondegeneracy condition is not assumed.
The author believes that the results obtained in this paper will play a substantial role for further development of the primal-dual interior point method for the NSDP. 

There exist two directions for future works. 
The first one is concerned with limiting behavior of the tangential direction $\dot{x}(\mu)$ in the $x$-space as $\mu\to 0$. 
We have the following conjecture: 
$$
\lim_{\mu\to 0}\dot{x}(\mu)=\xiast,
$$
where $\xiast$ was defined in \cref{lem:xi}.  
For nonlinear optimization, the corresponding result was proven by \citet[Theorem~12]{wright2002properties}.
The second direction of future works is to mitigate the strict complementarity (SC) condition from our assumptions.  
The SC actually plays a key role in our analysis, particularly when establishing \cref{prop:mukdk}, a base for proving the subsequent theorems.
{\cred
  Indeed, the following example, which is obtained from \citep[Section~4]{wright2002properties} with slight modification, shows that \cref{prop:mukdk} does not hold when the SC condition fails:
  $$
  \min_{x_1,x_2,x_3}\ \frac{1}{2}(x_1^2+x_2^2+x_3^2)\ \mbox{s.t. }
   \begin{bmatrix}
    x_1-1&x_3&&\\
    x_3&x_2&&\\
       & &x_1-1&x_3\\
       & &x_3&x_2 
  \end{bmatrix}\in \mathbb{S}^4_+.
   $$
   Its optimum $\xast$ is only $(1,0,0)^{\top}$, and the set of corresponding dual matrices is
   $$
\Set{   \begin{bmatrix}
    \lambda_1&0 &&\\
    0 &0 &&\\
       & &\lambda_2&0\\
       & &0 &0 
  \end{bmatrix}|
\lambda_1+\lambda_2=1,\lambda_1\ge 0,\lambda_2\ge 0}.
$$
It is easy to find that neither the SC nor the nondegeneracy condition holds, while both the MFCQ and {\SSOSC} hold at $\xast$.
For barrier parameter $\mu>0$, the BKKT point is $\xmu=(\frac{1+\sqrt{8\mu+1}}{2},\sqrt{2\mu},0)^{\top}$, which converges to the optimum $\xast$ as $\mu\to 0$. However, $\|\xmu-\xast\|=\sqrt{\frac{8\mu+1-\sqrt{8\mu+1}}{2}}\le 2\sqrt{\mu}\neq {\rm \Theta}(\mu)$, and therefore \cref{prop:mukdk} cannot hold without the SC. 
}
%We can find several works which does not request that condition.
%For example, 

\vspace{0.5em}\noindent{\bf Acknowledgments}: The author thanks Professor Yoshiko Ikebe for much advice.
He is also grateful to anonymous referees for many valuable comments and suggestions.

\def\thesection{A}
\renewcommand{\theequation}{A.\arabic{equation}}
\renewcommand{\thetheorem}{A.\arabic{theorem}}
\renewcommand{\thelemma}{A.\arabic{lemma}}
\setcounter{equation}{0}
\setcounter{theorem}{0}
\setcounter{lemma}{0}
\section{Omitted Proofs}\label{sec:appendix1}
In this appendix, we give the proofs which are not shown in the main part of this paper.
\subsection{Proof of \cref{lem:sigma}}
\label{subsec:lemma-1}
%\begin{proof}
The second equality follows from the direct calculation along with the fact of $\DeltaGUN=(\DeltaGNU)^{\top}$,
and the first one is derived from the following transformation:
\begin{align*}
d^{\top}\Omega(x^{\ast},Y)d=& 2\Tr\left(\sum_{j=1}^n\sum_{i=1}^nd_id_jY\Gi^{\ast} G_{\ast}^{\dag}\G_j^{\ast}\right)\\
=& 2\Tr\left(\left(P_{\ast}^{\top}YP_{\ast}\right)\left(\sum_{i=1}^nd_iP_{\ast}^{\top}\Gi^{\ast}P_{\ast}\right)\left(P_{\ast}^{\top}G_{\ast}^{\dag}P_{\ast}\right)\left(\sum_{j=1}^nd_jP_{\ast}^{\top}\G_j^{\ast}P_{\ast}\right)\right)\\
&
2
\Tr\left(
\begin{bmatrix}
\YUU&O\\
O&O
\end{bmatrix}
%\left(P_{\ast}^{\top}\sum_{i=1}^nd_i\Gi^{\ast}P_{\ast}\right)
\begin{bmatrix}
\DeltahGUU(x^{\ast};d)&\DeltahGUN(x^{\ast};d)\\
\DeltahGNU(x^{\ast};d)&\DeltahGNN(x^{\ast};d)
\end{bmatrix}
\begin{bmatrix}
O&O\\
O&(\Gast^{\NN})^{-1}
\end{bmatrix}
\begin{bmatrix}
\DeltahGUU(x^{\ast};d)&\DeltahGUN(x^{\ast};d)\\
\DeltahGNU(x^{\ast};d)&\DeltahGNN(x^{\ast};d)
\end{bmatrix}\right)
\notag\\
=&
2
\Tr\left(
\YUU\DeltahGNU(x^{\ast};d)
(\Gast^{\NN})^{-1}\DeltahGUN(x^{\ast};d)
\right).
\end{align*}
The proof is complete.
\hfill$\Box$

\subsection{Proof of \cref{prop:analyticunique}}\label{app:prop:analyticunique}
%{\it Proof of \cref{prop:analyticunique}: }
%We begin by proving that a global optimum of \eqref{eq:analy2} {\it exists uniquely}.
We first show the first assertion: the unique existence of the analytic center at $\xast$.
Note that because of relation\,\eqref{eq:past} for $(Y,z)\in\Lambda(\xast)$, 
\eqref{eq:analy} is equivalent to the following problem with respect to only $Y$:
\begin{align}
\begin{array}{cc}
\displaystyle{\min_{Y\in \mathbb{S}^m}}&-\log\det \YUU \\
\mbox{ ${\rm s.t}$. }&\nabla f(x^{\ast})-\mathcal{J}\GUU(x^{\ast})^{\ast}\YUU \in {\rm Im} \nabla h(x^{\ast}),\\ 
      &\YUN=\YNU=O,\ \YNN=O,\\
      &\YUU\in \mathbb{S}^{m-\rast}_{+},\\
\end{array}
\label{eq:analy3}
\end{align}
where $\mathcal{J}\GUU(x^{\ast})^{\ast}Z:=\left[\left(\PU^{\top}\Gi(x^{\ast})\PU\right)\bullet Z\right]_{i=1}^{n}\in \R^{n}$ for $Z\in \mathbb{S}^{m-\rast}$. 
%Moreover, the first consrtaint in the above is equivalent to the one in \eqref{eq:analy2} un%der the presence of the other constrains, because $\mathcal{J}\GUU(x^{\ast})^{\ast}\YUU=\mathcal{J}G(x^{\ast})^{\ast}Y$ holds. 

We establish existence of optima of \eqref{eq:analy}. 
By the strict complementarity condition as for the NSDP, there exists $(Y,z)\in \Lambda(\xast)$ such that $Y+\Gast\in \mathbb{S}^m_{++}$, which implies $\YUU\in \mathbb{S}^{m-\rast}_{++}$.
This means that a finite objective value of \eqref{eq:analy3} is attained at such a matrix $Y$.  Moreover, as $\Lambda(\xast)$ is convex and bounded from the MFCQ at $\xast$ for the NSDP, so is the feasible region of \eqref{eq:analy3}. 
By combining these facts, 
\eqref{eq:analy3} is ensured to have an optimum, say $\Ya\in \mathbb{S}^m_+$.
%After all, by the strict convexity of $-\log\det\YUU$ with respect to $Y$, \eqref{eq:analy3} is ensured to have a unique optimum, say $\Ya\in \mathbb{S}^m_+$.
From the full column rankness of $\nabla h(\xast)$, we see that 
the linear equation $\nabla f(x^{\ast})-\mathcal{J}\GUU(x^{\ast})^{\ast}\Ya^{\UU} +\nabla h(x^{\ast})z=0$
has a unique solution $z\in \R^s$, written $\za$. This $(\Ya,\za)$ is nothing but an optimum of \eqref{eq:analy}. 

Next, consider the following problem:
\begin{align}
\begin{array}{cc}
\displaystyle{\min_{Z}}&-\log\det Z\\
\mbox{ ${\rm s.t}$. }&\nabla f(x^{\ast})-\mathcal{J}\GUU(x^{\ast})^{\ast}Z \in {\rm Im} \nabla h(x^{\ast}),\\ 
      &Z\in \mathbb{S}^{m-\rast}_{+}.\\
\end{array}
\label{eq:analy4}
\end{align}
For a feasible point $Y$ of \eqref{eq:analy3}, $\YUU$ is clearly feasible to \eqref{eq:analy4}, and hence 
so is $\Ya^{\UU}$ to \eqref{eq:analy4}. Furthermore, we can ensure that $\Ya^{\UU}$ is optimal to \eqref{eq:analy4}.
Indeed, if not, there exists $Z$ such that $Z$ is feasible to \eqref{eq:analy4} and $-\log\det Z < -\log\det\Ya^{\UU}$.
Since $Y:=\PU Z\PU^{\top}\in \mathbb{S}^m_+$ is feasible to \eqref{eq:analy3} and $\log\det \YUU=\log\det Z$, we gain 
$-\log\det \YUU < -\log\det \Ya^{\UU}$, a contradiction to the optimality of $\Ya^{\UU}$ for \eqref{eq:analy3}.
Lastly, since \eqref{eq:analy4} is a strictly convex problem, we see that $\Ya^{\UU}$ is a {\it unique} optimum of \eqref{eq:analy4}.

In turn, we establish the uniqueness of $(\Ya,\za)\in \mathbb{S}^m_+\times \R^{\ell}$ as optimum of 
\eqref{eq:analy}. 
To derive a contradiction, assume that there exist two distinct optima $(\Ya,\za)$ and $(\hYa,\hza)$ at $x^{\ast}$, which 
yields that $\Ya$ and $\hYa$ are both optima of \eqref{eq:analy3} by the preceding argument. 
Thus, so are $\Ya^{\UU}$ and $\hYa^{\UU}$ to \eqref{eq:analy4}, in particular $\Ya^{\UU}=\hYa^{\UU}$, according to the preceding argument again. 
%Since a global optimum of \eqref{eq:analy} is unique according to the preceeding argument, 
%we obtain $(\YUUa,\za)=(\hYUUa,\hza)$.
Hence, we have  
\begin{equation*}
P_{\ast}^{\top}(\Ya-\hYa)P_{\ast}=\begin{bmatrix}
\Ya^{\UU}-\hYa^{\UU}
& O\\
O& O 
\end{bmatrix} = O.
\end{equation*}
%where the first equality follows from the relation in \eqref{eq:past} together with $(\Ya,\za),(\hYa,\hza)\in\Lambda(\xast)$. 
Since $P_{\ast}$ is nonsingular, we obtain $\Ya=\hYa$, which together with 
the full column rankness of $\nabla h(\xast)$ implies $\za=\hza$. 
Hence we ensure $(\Ya,\za)=(\hYa,\hza)$, which is a contradiction.  Consequently, 
\eqref{eq:analy} has a unique optimum,  and thus we obtain the first claim.

There remains to verify the second claim as for \eqref{eq:proofanalytic}. 
If $(\Ya,\za)$ is the analytic center at $\xast$, $(\Ya,\za)\in \Lambda(\xast)$ holds by definition, and from the above proof, 
$\Ya^{\UU}$ is the unique optimum of \eqref{eq:analy4}. Hence, by the KKT conditions of \eqref{eq:analy4}, 
there exists $v\in \R^s$ such that \eqref{eq:proofanalytic} holds. 
Conversely, if such $v$ exists and $(\Ya,\za)\in \Lambda(\xast)$,
$\Ya^{\UU}$ solves \eqref{eq:analy4}, and hence
$\PU\Ya^{\UU}\PU^{\top}=\Ya$ does \eqref{eq:analy3}. 
This means that $(\Ya,\za)$ is the analytic center. The whole proof is complete.   
$\hfill\Box$
\subsection{Proof of \cref{thm:barrier}}\label{app:thm:barrier}
%\begin{proof}
\tred{Since $\xast$ is a strict local optimum because of the {\SSOSC}, we can take} a compact set $B\subseteq \R^n$ with nonempty interior such that $x^{\ast}\in {\rm int}\,B$ and it is a unique optimum of the problem
  \begin{equation}
\min\ f(x)\ \mbox{s.t. }h(x)=0,\ G(x)\in \mathbb{S}^m_{+},\ x\in B.\label{eq:probB}
  \end{equation}
  Consider the sequence of the relevant barrier problems parameterized with $\mu_k$ as in the following:  
  \begin{equation}
\min\ f(x)-\mu_k\log\det G(x)\ \mbox{s.t. }h(x)=0,\ G(x)\in \mathbb{S}^m_{++},\ x\in B,\label{eq:logB}
  \end{equation}
  and let $x^k$ be an optimum of problem\,\eqref{eq:logB} for each $k$.
  
 We will prove the theorem by showing that the above-defined sequence $\{x^k\}$ is nothing but the desired one.
To this end, it suffices to prove that $\{x^k\}$ converges to $x^{\ast}$. 
Indeed, because $x^{\ast}\in {\rm int}\,B$, the constraint $x\in B$ for problem\,\eqref{eq:logB} is inactive at $x^k$ for sufficiently large $k$, and thus
$x^k$ eventually becomes a local optimum of \eqref{eq:log}.

We write $f_k:=f(\xk)$ for each $k$ and $f_{\ast}:=f(\xast)$ for the sake of simplicity.
Recall $\Gk=G(\xk)$ and $\Gast=G(\xast)$.
We first consider the case~(i) where $\Gast\in \mathbb{S}^m_{+}\setminus \mathbb{S}^m_{++}$, i.e., $\Gast$ is on the boundary of $\mathbb{S}^m_{+}$ and thus $\det\,\Gast=0$.
The proof for the other case~(ii) where $\Gast\in \mathbb{S}^m_{++}$ will be given later.
Letting $\phi_k:=f_k-\mu_k\log\det \Gk$ for each $k$, the first goal is to prove  
\begin{equation}
\lim_{k\to\infty}\phi_k=f_{\ast}. \label{eq:limphik0}
\end{equation}
Without loss of generality, by re-taking a smaller $B$ with ${\rm int}\,B\ni x^{\ast}$ if necessary, we can suppose that $\det G(x)<1$ for all $x\in B$ because of $\det \Gast=0$, yielding 
\begin{equation}
  -\log\det G(x)>0,\ \ \forall x\in B,\label{eq:logdet1}
\end{equation}
which together with the feasibility of $x^k$ for \eqref{eq:probB} implies 
\begin{equation}
  -\mu_k\log\det \Gk>0>f_{\ast}-f_k.
  \label{eq:2}
\end{equation}
Using the two inequalities in \eqref{eq:2} yields
\begin{align}
  f_{\ast}&<f_k\notag\\
  &<f_k-\mu_k\log\det \Gk\ (=\phi_k) \notag \\
  &\le f_{k-1}-\mu_{k}\log\det G_{k-1}\notag\\
  &\le f_{k-1}-\mu_{k-1}\log\det G_{k-1}\ (=\phi_{k-1}), 
\end{align}
where the third inequality follows from the optimality of $x^k$ for problem\,\eqref{eq:logB} and the fourth one is due to $\mu_k\le \mu_{k-1}$ and
$-\log\det G_{k-1}>0$ from \eqref{eq:logdet1} and $x^{k-1}\in B$. 
From the above inequalities, we find that $\{\phi_k\}$ is a nonincreasing sequence such that it is bounded by $f_{\ast}$ from below. Therefore, we 
ensure the existence of $\lim_{k\to\infty}\phi_k$ and moreover obtain 
\begin{equation}
f_{\ast}\le \lim_{k\to\infty}\phi_k.\label{eq:limitphi}
\end{equation}
To verify \eqref{eq:limphik0}, there remains to prove the converse inequality.
Related to $\{x^k\}$, under the MFCQ at $\xast$,
we can construct another sequence $\{x^{\ell(k)}\}$ feasible to problem\,\eqref{eq:probB} such that it converges to $x^{\ast}$ and also satisfies $\det G_{\ell(k)}=\mu_k$ for each $k\ge K$ with sufficiently large $K>0$.\footnote{This fact is verified as follows:
From the MFCQ at $\xast$, there exists $d\in \R^n$ such that $\Gast+\DeltaG(x;d)\in \mathbb{S}^m_{++}$ and $\nabla h(\xast)^{\top}d=0$. By the full column rankness of $\nabla h(\xast)$, 
we can ensure existence of a smooth curve 
$x(\cdot):[0,\bar{t}]\to \R^s$ with some $\bar{t}>0$ such that 
$x(0) = \xast$, $\dot{x}(0) = d$, and 
$x(t)$ is feasible to \eqref{eq:logB}, $\forall t\in (0,\bar{t}]$.
Particularly, $G(x(t))\in \mathbb{S}^m_{++}$ holds for all $t\in (0,\bar{t}]$.
Therefore, as $\det G(x(t))$ is continuous w.r.t. $t\ge 0$ and takes $0$ at $t=0$ by the assumption $\Gast\in \mathbb{S}^m_+\setminus \mathbb{S}^m_{++}$, we conclude that for any sufficiently small $\alpha>0$, 
$\det G(x(t))=\alpha$ is attained by some $t\in (0,\bar{t}]$. The proof is complete.
} 
We then obtain $\lim_{k\to\infty}\phi_k\le f_{\ast}$ since
$a\log a\to 0$  as $a\to 0+$ and $\phi_k\le f_{\ell(k)}-\mu_k\log\det G_{\ell(k)}$ holds by the definition of $x^k$.
Together with \eqref{eq:limitphi}, it derives the target equation\,\eqref{eq:limphik0}.

The convergence of $\{x^k\}$ to $x^{\ast}$ is not difficult to derive from \eqref{eq:limphik0}. 
%Recall that $\{x^k\}$ is bounded since $\{x^k\}\subseteq B$. 
Letting $\bar{x}$ be an arbitrary accumulation point of $\{x^k\}$ and taking into consideration $-\mu_k\log\det \Gk>0$ in $\phi_k$, 
we get $\limsup_{k\to\infty}\phi_k\ge f(\bar{x})$, which combined with \eqref{eq:limphik0} implies $f_{\ast}\ge f(\bar{x})$.
By the feasibility of $\bar{x}$ and the unique optimality of $x^{\ast}$ for \eqref{eq:probB}, 
we gain $x^{\ast}=\bar{x}$.
Finally, since $\bar{x}$ was an arbitrary accumulation point of $\{x^k\}$, we conclude that $\lim_{k\to\infty}x^k=x^{\ast}$.

We next consider case~(ii) where $\Gast\in \mathbb{S}^m_{++}$.
Note that $\log\det \Gast$ is finite in this case.
Without loss of generality, we may assume that $\det G(x)>0$ for all $x\in B$, by taking a smaller $B(\ni x^{\ast})$ if necessary.
Let $\bar{x}$ be an arbitrary accumulation point of $\{x^k\}$ and note that
$\log\det G(\bar{x})$ is also finite since $\det G(\bar{x})>0$ by virtue of $\bar{x}\in B$.
By the optimality of $x^{k}$ and feasibility of $x^{\ast}$ to \eqref{eq:logB}, it follows that 
%\begin{align}
$$
f_k-\mu_k\log\det \Gk\le f_{\ast} - \mu_k\log\det \Gast,
$$
where driving $k\to \infty$ and taking a subsequence if necessary imply $f(\bar{x})\le f_{\ast}$.
Then, in virtue of feasibility of $\bar{x}$ and unique optimality of $x^{\ast}$ for \eqref{eq:probB}, we have $x^{\ast}=\bar{x}$.
We hence conclude that $\lim_{k\to\infty}x^k=x^{\ast}$ as for case(ii).
Consequently, the desired result is obtained and the proof is complete.
\hfill$\Box$

\subsection{Proof of \cref{prop:mukdk}}\label{app:prop:mukdk}
We use the notations described before \cref{prop:mukdk}.
In particular, recall \eqref{eqn:1212-1} and \eqref{eqn:1212-2}.
\vspace{0.5em}
 
{\it Proof of \cref{prop:mukdk}:} To begin with, for each $k\ge 0$, let $$\hd^k:=\frac{d^k}{\|d^k\|}.$$
Since $\{\hd^k\}$ is bounded, it has at least one accumulation point, say $\hd^{\ast}$. 
Choose an arbitrary subsequence $\{\hd^k\}_{k\in \K}$ which converges to $\hd^{\ast}$. 
From \cref{rem:20211217}, $\{\wk\}_{k\in\K}$ has an accumulation point, say $\wast:=(\xast,\Yast,\zast)$. Without loss of generality, we assume $\lim_{\K\ni k\to \infty}\wk=\wast$.

We prove the assertion by two steps.
As the first step, we prove
\begin{equation}
\liminf_{k\to\infty}\frac{\mu_k}{\|d^k\|}>0. \label{eq:20211215-1}
\end{equation}
  In order to derive a contradiction, suppose to the contrary that there exists a subsequence of $\left\{\frac{\mu_k}{\|d^k\|}\right\}_{k\in \K}$ such that it converges to 0. We may assume $\lim_{k(\in\K)\to\infty}\frac{\mu_k}{\|d^k\|}=0$ by retaking $\K$ if necessary. Since $\lim_{k(\in\K)\to\infty}\hd^k=\hd^{\ast}$, $\hd^{\ast}$ satisfies  
\begin{equation}
\left(\PU^{\top}\DeltaG(x^{\ast};\hd^{\ast})\PU=\right)\DeltahGUU(x^{\ast};\hd^{\ast})\in \mathbb{S}^{\rast}_{+},\ \nabla h(x^{\ast})^{\top}\hd^{\ast}=0,
\label{eq:sm++}
\end{equation}
where these relations are derived from dividing the following equations by $\|d^k\|$ and passing to the limit:
\begin{eqnarray}
&\mathbb{S}^r_{++}\ni \PU^{\top}\Gk\PU=\PU^{\top}\left(\Gk-\Gast\right)\PU
=\DeltaGUU(\xast;d^k)+\O(\|d^k\|^2),&\notag \\
&0=h(\xk)=h(\xast)+\nabla h(\xast)^{\top}d^k+\O(\|d^k\|^2).&\notag
\end{eqnarray}

As $w^k=(\xk,\Yk,\zk)$ satisfies the BKKT conditions and
$\Past=[\PU,\PN]$ is an orthogonal matrix, we obtain, for each $k\in\K$,
\begin{align}
\frac{\mu_kI_{r_{\ast}}}{\|d^k\|}&=\frac{\PU^{\top}G_kY_k\PU}{\|d^k\|}\notag \\
                            &=\frac{\PU^{\top}\left(\Gast+\DeltaG(x^{\ast};d^k)+\O(\|d^k\|^2)\right)\begin{bmatrix}\PU\ \PN
                            \end{bmatrix}
                            \begin{bmatrix}
                            \PU^{\top}\\
                            \PN^{\top}
                            \end{bmatrix}
                            Y_k\PU}{\|d^k\|},\notag
\end{align}
which together with driving $k(\in\K)\to \infty$ yields 
\begin{equation}
\DeltahGUU(x^{\ast};\hd^{\ast})\YUU_{\ast}=O,\label{eq:deltaguu0}
\end{equation}
where we have used the relations $\Gast^{\UU}=O$ and $\Yast^{\NU}=O$ from \eqref{eq:anotherform}.
%% $$
%% \PU^{\top}\Gast\PU=O,\ 
%% \begin{bmatrix}
%%                             \PU^{\top}\Yast\PU\\
%%                             \PN^{\top}\Yast\PU
%%                             \end{bmatrix}=\begin{bmatrix}
%%                             {\YUU_{\ast}}\\
%%                             O
%%                             \end{bmatrix}$$
%%                             from \eqref{eq:anotherform}.

As $w^{\ast}=(x^{\ast},\Yast,\zast)$ and $(x^{\ast},\Ya,\za)$ satisfy the KKT conditions, it follows that 
\begin{align}
\nabla f(x^{\ast})&=\mathcal{J}G(x^{\ast})^{\ast}Y_{\ast}-\nabla h(x^{\ast})z^{\ast},
\label{eq:kktmathj1}\\
%\end{equation}
%and  
%\begin{equation}
&=\mathcal{J}G(x^{\ast})^{\ast}\Ya-\nabla h(x^{\ast})\za. \label{eq:kktmathj2}
\end{align} 
Pre-multiplying both \eqref{eq:kktmathj1} and \eqref{eq:kktmathj2} by $(\hd^{\ast})^{\top}$ and noting \eqref{eq:deltaguu0} lead to 
\begin{equation}
\nabla f(x^{\ast})^{\top}\hd^{\ast}={\rm Tr}\left(
\DeltahGUU(x^{\ast};\hd^{\ast})\YUU_{\rm a}\right)
={\rm Tr}\left(
\DeltahGUU(x^{\ast};\hd^{\ast})\YUU_{\ast}\right)
=0.\label{eq:threeeq}
\end{equation}
From \eqref{eq:threeeq} and \eqref{eq:sm++}, we ensure  
\begin{equation}
\hd^{\ast}\in C(x^{\ast}),\label{eq:hdcritical}
\end{equation}
where $C(x^{\ast})$ is defined in \eqref{eq:critical}.
As $\YUU_{\rm a}$ is positive definite by definition and 
$\DeltahGUU(x^{\ast};\hd^{\ast})\in \mathbb{S}^m_{+}$ follows from \eqref{eq:sm++} again, 
$
{\rm Tr}\left(
\DeltahGUU(x^{\ast};\hd^{\ast})\YUU_{\rm a}\right)=0
$
in \eqref{eq:threeeq} yields
\begin{equation}
\DeltahGUU(x^{\ast};\hd^{\ast})=O.
\label{deltaguu0}
\end{equation}

Next, we transform $(\hd^{\ast})^{\top}\mathcal{J}G(x^{\ast})^{\ast}(Y_k-\olY)$ as 
\begin{align}
(\hd^{\ast})^{\top}\mathcal{J}G(x^{\ast})^{\ast}(Y_k-\Yast)&=\DeltaG(x^{\ast};\hd^{\ast})\bullet (Y_k-\Yast)\notag\\
                                                               &=\Tr\left(\begin{bmatrix}\PU^{\top}\\
                                                               \PN^{\top}
                                                               \end{bmatrix}
                                                               \DeltaG(x^{\ast};\hd^{\ast})
                                                               \begin{bmatrix}\PU\ \PN
                                                               \end{bmatrix}
                                                               \begin{bmatrix}\PU^{\top}\\
                                                               \PN^{\top}
                                                               \end{bmatrix}
                                                               (Y_k-\Yast)
                                                               \begin{bmatrix}\PU\ \PN
                                                               \end{bmatrix}
                                                               \right)
                                                               \notag\\
                                                               %%%
                                                                &=\begin{bmatrix}
                                                                \DeltahGUU(x^{\ast};\hd^{\ast})&\DeltahGUN(x^{\ast};\hd^{\ast})\\
                                                                \DeltahGNU(x^{\ast};\hd^{\ast})&\DeltahGNN(x^{\ast};\hd^{\ast})
                                                               \end{bmatrix}
                                                               \bullet
                                                               \begin{bmatrix}
                                                                \DeltaYUU_k-\YNN_{\ast}&\DeltaYUN_k\\
                                                                \DeltaYNU_k           &\DeltaYNN_k
                                                               \end{bmatrix}
  \notag\\
                                                                &=\begin{bmatrix}
                                                               O&\DeltahGUN(x^{\ast};\hd^{\ast})\\
                                                                \DeltahGNU(x^{\ast};\hd^{\ast})&\DeltahGNN(x^{\ast};\hd^{\ast})
                                                               \end{bmatrix}
                                                               \bullet
                                                               \begin{bmatrix}
                                                                \DeltaYUU_k-\YNN_{\ast}&\DeltaYUN_k\\
                                                                \DeltaYNU_k           &\DeltaYNN_k
                                                               \end{bmatrix}
  \notag\\    
                                                              &=2{\rm Tr}\left(\DeltahGUN(x^{\ast};\hd^{\ast})\DeltaYNU_k\right)
                                                              +{\rm Tr}\left(\DeltahGNN(x^{\ast};\hd^{\ast})\DeltaYNN_k\right),   
                                                                                                                        \label{al:2237}
\end{align}
where 
the second equality follows from the fact that $[\PU,\PN](=P_{\ast})$ is an orthogonal matrix and 
the fourth one is due to \eqref{deltaguu0}.

Since $w^k$ and $\wast$ satisfy the BKKT and KKT conditions, respectively, we have
$\nabla_{x}L(\wast)=\nabla_xL(w^k)=0$ for each $k\in\K$, yielding 
\begin{align}
0&=(\hd^{\ast})^{\top}\frac{\left(\nabla_xL(w^k)-\nabla_{x}L(\wast)\right)}{\|d^k\|}\notag\\
  &=(\hd^{\ast})^{\top}
 \frac{
  \nabla^2_{xx}L(\wast)(d^k)
  -\mathcal{J}G(x^{\ast})(Y_k-\Yast)
  +\nabla h(x^{\ast})^{\top}(z^k-\zast)
  +\O(\|d^k\|^2)
  }{
  \|d^k\|
  } 
  \notag\\ 
  &=
 \frac{
  (\hd^{\ast})^{\top}\nabla^2_{xx}L(\wast)(d^k)
  -\Delta G(x^{\ast};\hd^{\ast})\bullet (Y_k-\Yast)
  +\O(\|d^k\|^2)
  }{
  \|d^k\|
  } 
  \notag\\
  &=
  (\hd^{\ast})^{\top}\nabla^2_{xx}L(\wast)\frac{d^k}{\|d^k\|}
  -
  \frac{2{\rm Tr}\left(\DeltahGUN(x^{\ast};\hd^{\ast})\DeltaYNU_k\right)
                                                              +{\rm Tr}\left(\DeltahGNN(x^{\ast};\hd^{\ast})\DeltaYNN_k\right)}{\|d^k\|}+\O(\|d^k\|),\label{al:2330}
\end{align}
where the last equality follows from \eqref{al:2237}.
Notice that the off-diagonal elements of $P_{\ast}^{\top}G_kY_kP_{\ast}(=\mu_k I)$ are zeros for all $k$.
Hence, for each $k(\in \K)\ge 0$, we have 
\begin{align}
O&=\PN^{\top}G_kY_k\PU\notag \\
   &=\PN^{\top}G_k
\begin{bmatrix}\PU,\PN\end{bmatrix}
\begin{bmatrix}\PU^{\top}\\ \PN^{\top}\end{bmatrix}
Y_k\PU\notag\\
&=
\PN^{\top}G_k\PU\YUU_k+\PN^{\top}G_k\PN\YNU_k\notag
\end{align}
for each $k\in\K$.
Substituting Taylor's expansion $G_k=\Gast+\DeltaG(x^{\ast};d^k)+\O(\|d^k\|^2)$ into the last equation yields 
\begin{equation*}
(\GNU_{\ast}+\DeltaGNU(x^{\ast};d^k))\YUU_k
+\left(\GNN_{\ast}+\DeltaGNN(x^{\ast};d^k)\right)\YNU_k=\O(\|d^k\|^2),
\end{equation*}
where $\|\Yk\|_{\rm F}=\O(1)$ was used for the last equality.
Noting $\GNU_{\ast}= O$ and dividing both the sides of the above by $\|d^k\|$ give
\begin{align}
&\frac{\DeltaGNU(x^{\ast};d^k)}{\|d^k\|}\YUU_k
+
\left(\GNN_{\ast}\frac{\YNU_k}{\|d^k\|}+\frac{\DeltaGNN(x^{\ast};d^k)}{\|d^k\|}\YNU_k\right)=\O(\|d^k\|)\label{al:0408}
\end{align}
for $k\in\K$.
Note that
$\lim_{k\to\infty}\YNU_k=O$ holds, which implies $\lim_{k\to\infty}\frac{\DeltaGNN(x^{\ast};d^k)}{\|d^k\|}\YNU_k=O$. Moreover, together with letting $k(\in\K)\to \infty$, equation\,\eqref{al:0408} implies  
\begin{align}
&\lim_{k\to\infty}\frac{\YNU_k}{\|d^k\|}=
-(\Gast^{\NN})^{-1}\DeltaGNU(x^{\ast};\hd^{\ast})\YUU_{\ast}.\label{al:ynukdk}
\end{align}
On the other hand, the $(2,2)$-block matrix of $P_{\ast}^{\top}G_kY_kP_{\ast}/\|d^k\|(=\mu_k I/\|d^k\|)$ is calculated as 
\begin{align}
  &\frac{1}{\|d^k\|}\left(\PN^{\top}\left(G_{\ast}+\Delta G(x^k;d^k)\right)\PU\YUN_k+
\PN^{\top}\left(G_{\ast}+\Delta G_k(x^k;d^k)\right)\PN\YNN_k
\right)
+\frac{\O(\|d^k\|^2)}{\|d^k\|}\notag\\
=&
\frac{1}{\|d^k\|}\left(\PN^{\top}\Delta G(x^k;d^k)\PU\YUN_k+
\left(\Gast^{\NN}+\PN^{\top}\Delta G_k(x^k;d^k)\PN\right)\YNN_k
\right)
+\O(\|d^k\|),\label{al:three}
\end{align}
where we have used
\begin{equation}
\Gast^{\NU}=O,\ \Gast^{\NN}=O\label{eq:0328-2317}
\end{equation}
from \eqref{eq:anotherform}.
In particular, from $\lim_{k\to\infty}\Delta G(x^k;d^k)/\|d^k\|=\Delta G(x^{\ast};\hd^{\ast})$ and $\lim_{k\to\infty}(\YUN_k,\YNN_k)=(O,O)$, we see 
$$\lim_{k\to\infty}\frac{\PN^{\top}\Delta G(x^k;d^k)\PU\YUN_k+\PN^{\top}\Delta G(x^k;d^k)\PN\YNN_k}{\|d^k\|}=O.$$
Since the limit of \eqref{al:three} is zero by the assumption $\lim_{k\to\infty}\mu_k/\|d^k\|=0$ again and
recalling that \eqref{al:three} is the $(2,2)$-block of $\mu_kI/\|d^k\|$, the above equation yields 
$
\lim_{k\to\infty}{\Gast^{\NN}\YNN_k}/{\|d^k\|}=O,
$
which combined with the nonsingularity of $\Gast^{\NN}$ induces 
\begin{equation}
\lim_{k\to\infty}\frac{\YNN_k}{\|d^k\|}=O.\label{al:ynnkdk}
\end{equation}

By taking \eqref{al:ynukdk} and \eqref{al:ynnkdk} into consideration and driving $k\to\infty$ in the equation in \eqref{al:2330},
it holds that 
$$
  (\hd^{\ast})^{\top}\nabla^2_{xx}L(\wast)\hd^{\ast}=
  -2{\rm Tr}\left(\YUU_{\ast}\DeltahGUN(x^{\ast};\hd^{\ast})
  (\Gast^{\NN})^{-1}\DeltaGNU(x^{\ast};\hd^{\ast})
  \right).$$
Combined with \cref{lem:sigma}, this equation further implies
$$(\hd^{\ast})^{\top}\nabla^2_{xx}L(\wast)\hd^{\ast}+(\hd^{\ast})^{\top}\Omega(x^{\ast},Y_{\ast})\hd^{\ast}=0.$$
However, in view of $\hd^{\ast}\in C(x^{\ast})$ by \eqref{eq:hdcritical} and $\hd^{\ast}\neq 0$ and by noting $(\Yast,\zast)\in \Lambda(\xast)$,
the above equation contradicts the \SSOSC. 
Therefore, we conclude \eqref{eq:20211215-1}.

In turn, we show $\frac{\mu_k}{\|d^k\|}=\O(1)$ as the second step. 
As $w^k$ satisfies the BKKT conditions, we obtain, for each $k$,
\begin{align}
\frac{\mu_kI_{r_{\ast}}}{\|d^k\|}=\frac{\PU^{\top}G_kY_k\PU}{\|d^k\|}
                            =\frac{\PU^{\top}\left(\Gast+\mathcal{J}G(x^{\ast})(d^k)+\O(\|d^k\|^2)\right)\begin{bmatrix}\PU\ \PN
                            \end{bmatrix}
                            \begin{bmatrix}
                            \PU^{\top}\\
                            \PN^{\top}
                            \end{bmatrix}
                            Y_k\PU}{\|d^k\|},\notag
\end{align}
which together with \eqref{eq:0328-2317}
%% the fact from \eqref{eq:anotherform} that 
%% $$
%% \Gast^{\UU}=O,\ 
%% \begin{bmatrix}
%%                             \PU^{\top}\Yast\PU\\
%%                             \PN^{\top}\Yast\PU
%%                             \end{bmatrix}=\begin{bmatrix}
%%                             {\YUU_{\ast}}\\
%%                             O
%%                             \end{bmatrix}$$
 yields 
\begin{equation*}
\lim_{k\to\infty}\frac{\mu_kI_{r_{\ast}}}{\|d^k\|}=\DeltahGUU(x^{\ast};\hd^{\ast})\YUU_{\ast}.
\end{equation*}
This means that the sequence $\left\{\frac{\mu_k}{\|d^k\|}\right\}$ is bounded, and we thus obtain the desired consequence. 
By combining \eqref{eq:20211215-1} and this fact, 
the proof is complete.
$\hfill\Box$
\subsection{Proof of \cref{lem:xi}}\label{app:lem:xi}
To start with, decompose a vector $\Delta x\in \R^n$ into orthogonal component vectors as follows: 
\begin{equation}
\Delta x = \UR \eta^1 + \UN\eta^2,
\label{eq:decompose}
\end{equation}
where $(\eta^1,\eta^2)\in \R^{\past}\times \R^{n-\past}$ and $\UN\in \R^{n\times (n-\past)}$ is a matrix whose columns form an orthonormal basis of the orthogonal complement subspace of $\mathUast$, where $\mathUast$ is defined in \eqref{eq:mathUast}.

Let $\UR^i$ be the $i$-th column of $\UR$ for each $i=1,2,\ldots,\past$.
%Decompose $\Delta x$ as in \eqref{eq:decompose}.
From \eqref{eqn:diff5}, we have
$$
{\rm Sym}\left(
\begin{bmatrix}
\YUU_{\rm a}\DeltaGUU(\xast;\Delta x)& \YUU_{\rm a}\DeltaGUN(\xast;\Delta x) \\
\GNN_{\ast}\Delta\YNU & \GNN_{\ast}\Delta\YNN
\end{bmatrix}\right)= I,
$$
of which the block components together with $\GNN_{\ast}\in \mathbb{S}^{m-\rast}_{++}$ and $\YUU_{\rm a}\in \mathbb{S}^{\rast}_{++}$
yield
\begin{eqnarray}
&\Delta\YNN= (\GNN_{\ast})^{-1},\label{al:0429-1}\\
&\DeltaGUU(\xast;\Delta x)=(\YUU_{\rm a})^{-1},\label{al:0429-2}\\ 
&\YUU_{\rm a}\DeltaGUN(\xast;\Delta x)+\Delta\YUN\GNN_{\ast}=O.\label{al:0429-3}
\end{eqnarray}
%By premutiplying \eqref{eqn:diff4} by $\UR^{\top}$, we obtain 
By \eqref{eqn:diff4}, we obtain
$$
\UR^{\top}\nabla^2_{xx}L(\wa)\Delta x-  \left(\begin{bmatrix}
O &\DeltaGUN(\xast;\UR^i)\\
\DeltaGNU(\xast;\UR^i) &\DeltaGNN(\xast;\UR^i)
\end{bmatrix}
\bullet \Past^{\top}\Delta Y\Past \right)_{i=1}^{\past}=0, 
$$
leading to
$$
\UR^{\top}\nabla^2_{xx}L(\wa)\Delta x
-  \left(
\DeltaGNU(\xast;\UR^i)\bullet \Delta \YUN
+\DeltaGUN(\xast;\UR^i)\bullet \Delta \YNU
+ \DeltaGNN(\xast;\UR^i)\bullet \Delta \YNN
\right)_{i=1}^{\past}
= 0, 
$$
which is moreover rephrased as 
$$
\UR^{\top}\nabla^2_{xx}L(\wa)\Delta x
-  \left(
2\DeltaGNU(\xast;\UR^i)\bullet \Delta \YUN
+ \DeltaGNN(\xast;\UR^i)\bullet \Delta \YNN
\right)_{i=1}^{\past}
= 0.
$$
Combined with \eqref{al:0429-1} and \eqref{al:0429-3}, this equation yields
\begin{align*}
&\UR^{\top}\nabla^2_{xx}L(\wa)\Delta x
+
2\left({\rm Tr}\left(\DeltaGNU(\xast;\UR^i)\YUU_{\rm a}\DeltaGUN(\xast;\Delta x) (\GNN_{\ast})^{-1}\right)\right)_{i=1}^{\past}=
\left(\DeltaGNN(\xast;\UR^i)\bullet(\GNN_{\ast})^{-1}\right)_{i=1}^{\past}.
%\left(-\DeltaGUN(\xast;\UR^i)\bullet (\YUU_{\rm a})^{-1}+\left( \DeltaGUN(\xast;\UR^i)-\DeltaGNN(\xast;\UR^i)\right)\bullet(\GNN_{\ast})^{-1}\right)_{i=1}^{\past}.
\end{align*}
Decomposing $\Delta x$ as in \eqref{eq:decompose}, we obtain from the above equation that 
\begin{align}
&\UR^{\top}\nabla^2_{xx}L(\wa)\UR\eta^1
+
2\left(\sum_{j=1}^{\past}\eta_j^1{\rm Tr}\left(\DeltaGNU(\xast;\UR^i)\YUU_{\rm a}\DeltaGUN(\xast;\UR^j) (\GNN_{\ast})^{-1}\right)\right)_{i=1}^{\past}\notag\\
=&  
-\UR^{\top}\nabla^2_{xx}L(\wa)\UN\eta^2
-2\left({\rm Tr}\left(\DeltaGNU(\xast;\UR^i)\YUU_{\rm a}\DeltaGUN(\xast;\UN\eta^2) (\GNN_{\ast})^{-1}\right)\right)_{i=1}^{\past}\notag \\
&\hspace{2em}+
\left(
-\DeltaGUN(\xast;\UR^i)\bullet (\YUU_{\rm a})^{-1}+
\left( \DeltaGUN(\xast;\UR^i)-\DeltaGNN(\xast;\UR^i)\right)\bullet(\GNN_{\ast})^{-1}\right)_{i=1}^{\past},
\label{al:eta1eta2}
\end{align}
where $\eta^1:=(\eta^1_1,\eta^1_2,\ldots,\eta^1_{\past})^{\top}$.

Next, we prove that $\eta^1$ and $\eta^2$ are uniquely determined. 
To this end, note that 
$\UN \eta^2\in \mathUast^{\perp}$ by definition, and 
$\DeltaGUU(x^{\ast};\UN \eta^2)=(\YUU_{\rm a})^{-1}$ follows
from \eqref{al:0429-2} and $\DeltaGUU(x^{\ast};\UR \eta^1)=O$.
From this, $\UN \eta^2$ turns out to be unique,\footnote{
More precisely speaking, to derive the uniqueness of $\UN \eta^2$, we have made use of the following fundamental result from linear algebra: 
given $A\in \R^{q_1\times q_2}$ and $b\in\R^{q_1}$, 
assume that the linear equation $A\theta=b$ has a nonempty solution set. 
Pick a solution $u$ arbitrarily and 
decompose it as $u=u^1+u^2$ with 
$u^1\in {\rm Ker}\,A$ and $u^2\in \left({\rm Ker}\,A\right)^{\perp}$, where ${\rm Ker}\,A$ denotes the kernel or null space of the matrix $A$. Then, $u^2$ is uniquely determined regardless of choice for $u$, whereas $u^1$ is free in ${\rm Ker}\,A$. 
In the proof, 
there exist correspondences between
$Au=b$ and $\DeltaGUU(x^{\ast};\UR \eta^1+\UN \eta^2)=(\YUU_{\rm a})^{-1}$, 
$u^1$ and $\UR \eta^1$,
$u^2$ and  $\UN \eta^2$, and 
${\rm Ker}\,A$ and $\mathUast$, respectively.
}
which together with the full column rank of $\UN$ yields the uniqueness of $\eta^2$.
In view of this fact
and \eqref{al:eta1eta2}, $\eta_1$ is also uniquely determined, because the matrix $$\UR^{\top}\nabla^2_{xx}L(\wa)\UR+
2\left({\rm Tr}\left(\DeltaGNU(\xast;\UR^i)\YUU_{\rm a}\DeltaGUN(\xast;\UR^j) (\GNN_{\ast})^{-1}\right)\right)_{1\le i\le j\le \past}$$ is actually positive definite by virtue of the {\SSOSC}.
%Together this fact, note that  
%\eqref{al:eta1eta2} is equivalent to \eqref{al:0429-1}
%under the presence of \eqref{al:0429-2} and \eqref{al:0429-3}. 

As a result, $\Delta x = \UR\eta^1 + \UN\eta^2$ is the unique $\Delta x$-component of solutions to equations\,\eqref{eqn:diff4}-\eqref{eqn:diff6}. 
Therefore, we ensure Item\,\ref{lem:xi-1}. 
Item\,\ref{lem:xi-2} follows immediately from \eqref{al:0429-1}-\eqref{al:0429-3} with $\Delta x = \xiast$.
$\hfill\Box$
{\cred
\subsection{Proof of \cref{prop:bkktcalPc}}\label{app:prop6}
We show item\,\ref{prop:bkktcalPc:item1}.
We first consider the first-half claim.
For contradiction, assume that there exists an infinite sequence $\{\muk\}\subseteq \R_{++}$ converging to $0$ such that
 $\calPc(\muk)$ does not contain a BKKT point with barrier parameter $\muk$ for each $k$.
According to \cref{thm:barrier}, $\{\muk\}$ accompanies a sequence of BKKT points $\{x({\muk})\}$ which converges to the KKT point $\xast$.
By the above assumption, $x({\muk})\notin \calPc(\muk)$ for each $k$, implying
$\|x({\muk})-\xast-\muk\xiast\|/\muk\ge \rho\|\xiast\|>0$.
However, \cref{prop:xiast} implies 
\begin{equation}
  \|x({\muk})-\xast-\muk\xiast\|=\o(\muk).\label{eq:231225-1}
\end{equation}
This is a contradiction. Hence, the first-claim claim is obtained.
%The second-half one can be also established by deriving a contradiction and using \eqref{eq:231225-1} in a similar way.
The second-half one can be also established by deriving a contradiction. Suppose to the contrary that
there exist BKKT points $\{x(\muk)\}$ with $x(\muk)\in\bd \calPc(\muk)$. By the definition of $\bd\calPc(\muk)$, we see
$\lim_{k\to\infty}x(\muk)=\xast$, thus \eqref{eq:231225-1} holds again. However, this contradicts $(x(\muk)-\xast)/\muk=\rho,\forall k$ from $x(\muk)\in\bd \calPc(\muk)$.  
Item\,\ref{prop:bkktcalPc:item2} also follows readily since the same relation as \eqref{eq:231225-1} is obtained from \cref{prop:xiast} again.
\hfill\Halmos\endproof
}
\subsection{Proof of \cref{prop:0513-1}}\label{app:prop0513-1}
%\tred{[[Read the whole proof again!]]}
We prove the first assertion in item\,\ref{prop:0513-1-item1-1}.  
Write $\widetilde{X}:=\Past^{\top}X\Past$ for any $X\in \mathbb{S}^m$ \tred{by convension}.
%In what follows, we set 
%$G(\cdot)$ and $\DeltaG(\cdot;\xiast)$ as $X$ here to express
%$\widetilde{G}(\cdot)$ and $\widetilde{\DeltaG}(\cdot\,;\,\cdot)$.
In particular, we set $G(\cdot)$ and $\DeltaG(\xast;\cdot)$ to $X$.  
Denote
$$
R(x,\mu) := \widetilde{G}(x) - \wtGast  - \mu\widetilde{\DeltaG}(\xast;\xiast).
$$
We next consider to bound the magnitude of $R(x,\mu)$ when $(x,\mu)$ is varied.
Recall $\GUU_{\ast}=O$, $\GUN_{\ast}=\GNU_{\ast}=O$, and 
$\GNN_{\ast}\in \mathbb{S}^{\tred{\rast}}_{++}$. 
It follows that  
\begin{align}
 \frac{1}{\mu}\widetilde{G}(x) & =  \frac{1}{\mu}\wtGast+ \widetilde{\DeltaG}(\xast;\xiast) +  \frac{1}{\mu}R(x,\mu)\notag \\
 & =
 %% \begin{bmatrix}
 %%    O&O\\
 %%    O&\GNN_{\ast}
 %% \end{bmatrix} +
 \begin{bmatrix}
  (\YUU_{\rm a})^{-1} + \frac{1}{\mu}R_1(x,\mu)&  \DeltaGUN(\xast;\xiast)+ \frac{1}{\mu}R_2(x,\mu)\\
   \DeltaGNU(\xast;\xiast)+ \frac{1}{\mu}R_2(x,\mu)^{\top} &  \DeltaGNN(\xast;\xiast)+ \frac{1}{\mu}\GNN_{\ast}+ \frac{1}{\mu}R_3(x,\mu)   
\end{bmatrix},\label{al:05010-1}
\end{align}
where $R_i(x,\mu)\ (i=1,2,3)$ represent block submatrices of $R(x,\mu)$ with appropriate sizes and the second equality follows from $\DeltaGUU(\xast;\xiast)=(\YUU_{\rm a})^{-1}$ by item-\ref{lem:xi-2} of \cref{lem:xi}.
Taylor's expansion of $\widetilde{G}$ at $\xast$ gives  
\begin{align}
R(x,\mu)&=\widetilde{\DeltaG}\left(\xast,x-\xast-\mu\xiast\right)+\O(\|x-\xast\|^2)\notag\\
&=\O(\mu\rho\|\xiast\|+\|x-\xast\|^2)\notag\\
&=\O(\mu\rho+\mu^2)
\label{eq:0510-1}
\end{align}
for $x\in \tred{\cl\calPc(\mu)}$, where the last equality follows since
$\|x-\xast\|\le \mu(\rho+1)\|\xiast\|$ by $x\in \tred{\cl\calPc(\mu)}$.
By \eqref{eq:0510-1}, the fact that $(\YUU_{\rm a})^{-1}\in \mathbb{S}^{m-\rast}_{++}$, and taking $\mu_1>0$ and $\rho_1>0$ so small, $\frac{1}{\mu}
\|R(x,\mu)\|$ can be so small that 
the $(1,1)$-block matrix of $\frac{1}{\mu} \widetilde{G}(x)$ is symmetric positive definite for any $(\rho,\mu)\in (0,\rho_1]\times (0,\mu_1]$, that is to say,
\begin{equation}
Q(x,\mu):=(\YUU_{\rm a})^{-1} +\frac{1}{\mu}R_1(x,\mu)\in \mathbb{S}^{m-\rast}_{++},\ \forall (\rho,\mu)\in (0,\rho_1]\times (0,\mu_1],\ x\in\tred{\cl}\calPc(\mu),
\label{eq:0519-1}
\end{equation}
which along with \eqref{eq:0510-1} implies 
\begin{equation}
Q(x,\mu)^{-1}= \YUU_{\rm a}\left(I+\O\left(\rho+\mu\right)\right)^{-1}.\label{eq:limitQ}
\end{equation}
Meanwhile, the Schur complement of $\frac{1}{\mu} \widetilde{G}(x)$ is expressed as 
$$
\Sc(x,\mu):=\DeltaGNN(\xast;\xiast)+\frac{1}{\mu}\GNN_{\ast} + \frac{1}{\mu}R_3- \left(\DeltaGNU(\xast;\xiast)+\frac{1}{\mu}R_2^{\top}\right)
Q^{-1}\left(\DeltaGUN(\xast;\xiast) + 
\frac{1}{\mu}R_2\right),
$$
where 
we have dropped the arguments $(x,\mu)$ from the functions $R_1$, $R_2$, $R_3$, and $Q$ for simplicity.
From \eqref{eq:0510-1}, 
by re-taking $(\mu_1,\rho_1)$ sufficiently small if necessary, we find that 
the above $\Sc$ is symmetric positive definite for any $\mu\in (0,\mu_1]$ because
$\frac{1}{\mu}\GNN_{\ast}\in \mathbb{S}^{\rast}_{++}$ is eventually dominant therein as $\mu>0$ gets smaller 
and $x\in \tred{\cl}\mathcal{P}_{\rho}(\mu)$ holds by assumption.
%In particular, the least eigenvalue of $\Sc(x,\mu)$ is bounded from below with some positive number, and thus
\tred{Hence,}
\begin{equation}
\Sc(x,\mu)^{-1} = \O(\mu)\ \ \ (x\in \cl\calPc(\mu)).\label{eq:1216-boundedsc}
\end{equation}
\tred{Moreover,
in view of \eqref{al:05010-1}, from \eqref{eq:0519-1} and $\Sc\in \mathbb{S}^{\rast}_{++}$ shown above}, we conclude $\frac{1}{\mu} \widetilde{G}(x)\in \mathbb{S}^m_{++}$ for any $(\rho,\mu)\in (0,\rho_1]\times (0,\mu_1]$ and $x\in \cl\calPc(\mu)$, implying $G(x)\in \mathbb{S}^m_{++}$. 
Setting $(\bmuone,\brhoone):=(\mu_1,\rho_1)$, we ensure the first assertion. 

We next prove the second assertion in item~\ref{prop:0513-1-item1-1}.
Taking the inverse of $\mu^{-1}\widetilde{G}(x)$ by applying the formula of the inverse of a partitioned matrix (e.g., \citet[Section~0.7.3]{horn2012matrix})
to \eqref{al:05010-1}, we obtain
\begin{align}
\mu \wt{G}(x)^{-1}=\begin{bmatrix}
M_{11} &M_{12}\\
M_{12}^{\top} &M_{22}
\end{bmatrix},
\label{al:muginv} 
\end{align}
where each block component is defined as 
\begin{align*}
&M_{11}:=Q(x,\mu)^{-1}+Q(x,\mu)^{-1}\left(
\DeltaGUN(\xast;\xiast)+ \frac{1}{\mu}R_2(x,\mu)
\right)\Scxu^{-1}\left(
\DeltaGUN(\xast;\xiast)+ \frac{1}{\mu}R_2(x,\mu)
\right)^{\top}Q(x,\mu)^{-1},\\
&M_{12}:=-Q(x,\mu)^{-1}\left(\DeltaGUN(\xast;\xiast)+ \frac{1}{\mu}R_2(x,\mu)\right)\Scxu^{-1},\\
&M_{22}:=\Scxu^{-1}.
\end{align*}
%In view of \eqref{eq:0510-1} and \eqref{eq:0519-1}, we find that,
%by taking $\rho_1$ again sufficiently small if necessary, the least eigenvalue of $Q(x,\mu)$ is bounded from below by a positive number as $\mu$ approaches 0. Hence, $Q^{-1}(x,\mu)=\O(1)$.
Moreover, we have $\frac{1}{\mu}R_i(x,\mu)=\O(\rho+\mu)\ (i=1,2,3)$
from \eqref{eq:0510-1}.
These facts together with \eqref{eq:limitQ}, \eqref{eq:1216-boundedsc}, and \eqref{al:muginv} yield 
\begin{equation}
M_{11}=\YUU_{\rm a}\left(I+\O\left(\rho+\mu\right)\right)^{-1}+\O(\mu), M_{12}=\O(\mu),\ M_{22}=\O(\mu),
\label{eq:m11m12m22}
\end{equation}
which together with $\mu\widetilde{G}(x)^{-1}=\mu \Past^{\top}G(x)^{-1}\Past$ implies that $\Set{\mu G(x)^{-1}| x\in \tred{\cl\calPcbone(\mu)},
\mu\in (0,\tred{\bmuone}]}$ is bounded.

In turn, we prove item~\ref{prop:0513-1-item2-2}.   
{\cred
First, in view of \eqref{eq:m11m12m22}, we obtain 
\begin{align}
M_{11}-\YUU_{\rm a}%&=\YUU_{\rm a}\left(I+\O\left(\rho+\mu\right)\right)^{-1}+\O(\mu)-\YUU_{\rm a}\notag\\                                           
      	      	   &=\YUU_{\rm a}\left(I+\O\left(\rho+\mu\right)\right)^{-1}\left(I-\left(I+\O\left(\rho+\mu\right)\right)\right)+\O(\mu)\notag\\
                   &=\O(\rho+\mu).\label{eq:231122}
\end{align}
%we obtain $\|\mu G(x)^{-1}-\Ya\|=\O(\rho+\mu)$.
}
%Let us consider the first assertion in this item. 
%In what follows,  
We drive $(x,\mu)\to (\xast,0)$ along with satisfying $x\in \cl\calPc(\mu)$. 
%and consider limit points and accumulation points of $M_{11}$, $M_{12}$, and $M_{22}$.
%From \eqref{eq:m11m12m22}, $M_{12}$ and $M_{22}$ converge to $O$.
%Meanwhile,
{\cred
From \eqref{eq:m11m12m22} and \eqref{eq:231122}, it follows that 
$$
\mu \wt{G}(x)^{-1}-\wt{\Ya}=
\begin{bmatrix}
M_{11}-\YUU_{\rm a}&M_{12}\\
M_{12}^{\top}&M_{22}
\end{bmatrix}
=\begin{bmatrix}
\O(\rho+\mu)&\O(\mu)\\
\O(\mu)&\O(\mu)
\end{bmatrix},
$$
where $\wt{\Ya}=\Past^{\top}\Ya\Past$.
%according to \eqref{eq:231122},
%Therefore, 
%given $\varepsilon>0$, by taking $\rho>0$ sufficiently small, 
%\tred{all accumulation points of $\mu\wt{G}(x)^{-1}$ lie from $\wt{\Ya}$ within $\varepsilon$.
Thus, we have $\|\mu G(x)^{-1}-\Ya\|_{\rm F}=\|\mu\wt{G}(x)^{-1}-\wt{\Ya}\|_{\rm F}=\O(\rho+\mu)$.
This means that there exists some $K_1>0$ such that 
$\|\mu G(x)^{-1}-\Ya\|_{\rm F}\le K_1(\rho+\mu)$ as claimed.}
%the first assertion is gained.}
%all accumulation points of $M_{11}$ lie from $\YUU_{\rm a}$ within $\varepsilon$. Hence, 
%by recalling the relationship of $\mu \wt{G}(x)^{-1}$ and $\mu G(x)^{-1}$, the first assertion is obtained.
%Next, we consider the second assertion in item~\ref{prop:0513-1-item2-2}.
%Using \eqref{eq:limitQ} and the assumption \tred{$\rho\to 0+$},
%the convergence of $Q(x,\mu)^{-1}$ to $\YaUU$ is ensured.
%Under the assumption that $\rho>0$ is driven to 0,
%$M_{11}$ converges to $\YUU_{\rm a}$.
%Therefore, we obtain $\mu$
%we obtain $\|\mu G(x)^{-1}-\Ya\|_{\rm F}=\|\mu \wt{G}(x)^{-1}-\wt{\Ya}\|_{\rm F}=\O(\rho+\mu)$.
%Then, in a similar manner to the first assertion in item~\ref{prop:0513-1-item2-2}, the second assertion is established with \eqref{eq:m11m12m22}.
%The whole proof is complete.
$\hfill\Box$

{\cred
\subsection{Proof of \cref{prop:231211-1}}\label{app:prop:231211-1}
First,
from $x\in\cl\calPc(\mu)$,
it follows that $\|x-\xast\|\le \|x-\xast-\mu\xiast\|+\mu\|\xiast\|\le (\rho+1)\|\xiast\|\mu$.
Second, it follows that $\|Y-\Ya\|_{\rm F}\le \|Y-\mu G(x)^{-1}\|_{\rm F}+\|\Ya-\mu G(x)^{-1}\|_{\rm F}
\le (\gamma_1+K_1)\mu+K_1\rho$
from \eqref{eqn:0524-1} and item~\ref{prop:0513-1-item2-2} of \cref{prop:0513-1}.
Moreover, since these inequalities yield
\begin{align*}
&\|\za+(\nabla h(x)^{\top}\nabla h(x))^{-1}\nabla h(x)^{\top}\left(\nabla f(x)-\mathcal{J}G(x)^{\ast}Y\right)\|\\
=&\|(\nabla h(x)^{\top}\nabla h(x))^{-1}
\nabla h(x)^{\top}\left(\nabla f(x)-\mathcal{J}G(x)^{\ast}Y+\nabla h(x)\za\right)\|\\
=&\|(\nabla h(x)^{\top}\nabla h(x))^{-1}
\nabla h(x)^{\top}\left(\nabla_xL(\wa)+\O(\|Y-\Ya\|_{\rm F}+\|x-\xast\|)\right)\|\\
%=&\O(\|Y-\Ya\|_{\rm F}+\|x-\xast\|)\\
=&\O(\rho+\mu),
\end{align*}
where the second and third equalities are derived from applying Taylor's expansion to
$\nabla_xL(x,Y,\za)$ at $\wa$ and the facts that $\nabla_x L(\wa)=0$ and
$\|(\nabla h(x)^{\top}\nabla h(x))^{-1}\nabla h(x)^{\top}\|_{\rm F}=\O(1)$ for $x\in \hball$, where $\hball$ is the ball defined in \eqref{eq:231216-3}, 
we obtain
\begin{align*}
\|z-\za\|&\le \|z+(\nabla h(x)^{\top}\nabla h(x))^{-1}\nabla h(x)^{\top}\left(\nabla f(x)-\mathcal{J}G(x)^{\ast}Y\right)\|+\\
&\hspace{2em}\|\za+(\nabla h(x)^{\top}\nabla h(x))^{-1}\nabla h(x)^{\top}\left(\nabla f(x)-\mathcal{J}G(x)^{\ast}Y\right)\|\\
&=\O(\rho+\mu),
\end{align*}
where
we have used the assumption 
$\|z+(\nabla h(x)^{\top}\nabla h(x))^{-1}\nabla h(x)^{\top}\left(\nabla f(x)-\mathcal{J}G(x)^{\ast}Y\right)\|\le \gamma_2\mu$.
Finally, using the above facts together with $\nabla_xL(\wa)=0$, we have
$$\nabla_xL(w)=\nabla_xL(\wa)+\O\left(\|x-\xast\|+\|Y-\Ya\|_{\rm F}+\|z-\za\|\right)
            =\O(\rho+\mu)$$
%            &=\O(\|x-\xast\|)+\O(\|\mu G(x)^{-1}-\Ya\|_{\rm F})+\O(\|Y-\mu G(x)^{-1}\|_{\rm F})\\
%            &+\O\left(\|\za+(\nabla h(x)\nabla h(x)^{\top})^{-1}\nabla h(x)^{\top}\left(\nabla f(x)-\mathcal{J}G(x)^{\ast}Y\right)\|\right)\\
%            &+\O\left(\|z+(\nabla h(x)\nabla h(x)^{\top})^{-1}\nabla h(x)^{\top}\left(\nabla f(x)-\mathcal{J}G(x)^{\ast}Y\right)\|\right)
%\end{align*}
where the first equality follows from Taylor's expansion of $\nabla_xL$ at $\wa$.
As a consequence, by taking $K_2>0$ sufficiently small, we ensure the desired inequalities.
$\hfill\Box$}

\subsection{Proof of \cref{prop:0531-1}}\label{app:prop:0531-1}
{\cred
To start with,
choose $\rho\le \brho_2$ and consider an arbitrary sequence $\{\wl=(\xl,\Yell,\zell)\}$
and $\{\mul\}\subseteq (0,\bmu_2]$ such that $\mul\to 0$ as $\ell\to \infty$
and \eqref{eqn:0524-0}, \eqref{eqn:0524-1}, and \eqref{eqn:0524-2} are fulfilled for each $\ell$.
Write $\Gell:=G(\xell)$ for each $\ell$. From \cref{prop:231211-1}, we see that $\{\wl\}$ is bounded and $\lim_{\ell\to \infty}\xell=\xast$.
Note that $\Gell^{-1}$ exists by virtue of $\xell\in \cl\calPc(\mul)$ and \eqref{eq:231202-1} along with $\bmu_2\le \bmu_1$ and $\brho_2\le \brho_1$, and also note that 
$\{\mul \Gell^{-1}\}$ are bounded from item\,\ref{prop:0513-1-item1-1} of \cref{prop:0513-1}.
Moreover, \eqref{eqn:0524-1} implies that $\{\Yell\}$ and $\{\mul \Gell^{-1}\}$ accumulate at identical points in $\mathbb{S}^m_{+}$.
Denote an arbitrary accumulation point of $\{(\Yell,\zell)\}$ by $(\Yast,\zast)$.
From the above argument and the fact that $\|\Gell\Yell-\mul I\|_{\rm F}=\|\Gell(\Yell-\mul \Gell^{-1})\|_{\rm F}\le
\|\Gell\|_{\rm F}\|\Yell-\mul \Gell^{-1}\|_{\rm F}\le \gamma_1\mul\|\Gell\|$,
we obtain
\begin{equation}
\Gast\Yast=O,\ \Yast\in \mathbb{S}^m_+.\label{eq:scGastYast}
\end{equation}
Moreover, from \cref{prop:231211-1}, for any $\ell$, we have
$
\max\left(\|\Yell-\Ya\|_{\rm F},\|\zell-\za\|\right)\le K_2(\rho+\mul),
$
where $K_2>0$ is the constant defined in \cref{prop:231211-1}.
Then by driving $\ell\to \infty$, we obtain 
\begin{equation}
\max\left(\|\Yast-\Ya\|_{\rm F},\|\zast-\za\|\right)\le K_2\rho.\label{eq:231216-5}
\end{equation}
For $X\in S^m$, define $\lambda_{\min}(X)$ as the least eigenvalue of $X$.  
From $\Gast+\Ya\in \mathbb{S}^m_{++}$ and \eqref{eq:231216-5}, it follows that
$
\lambda_{\rm min}(\Gast+\Yast) \ge \lambda_{\rm min}(\Gast+\Ya)+\lambda_{\rm min}(\Yast-\Ya)
\ge
\lambda_{\rm min}\left(\Gast+\Ya\right)-\|\Yast-\Ya\|_{\rm F}
\ge 
\lambda_{\rm min}\left(\Gast+\Ya\right)-K_2\rho,
$
which together with $\lambda_{\rm min}\left(\Gast+\Ya\right)>0$ from $\Gast+\Ya\in\mathbb{S}^m_{++}$ implies 
\begin{equation}
\Gast+\Yast\in\mathbb{S}^m_{++}\ \mbox{
when $0<\rho\le \frac{\lambda_{\rm min}\left(\Gast+\Ya\right)}{2K_2}$}.\label{eq:231202-3}
\end{equation}
Next,
let $\widetilde{K}:=\sum_{i=1}^s\|\nabla^2h_i(\xast)\|_{\rm F}+
\frac{n(n+1)}{2}\max_{1\le i,j\le n}
\left\|\frac{\partial^2G(\xast)}{\partial x_i\partial x_j}\right\|_{\rm F}
+n\|(\Gast^{\NN})^{-1}\|_{\rm F}\max_{1\le i\le n}\|\PU\calG_i(\xast)\PN\|^2_{\rm F}>0$ and 
recall that $C(\xast)$ is a critical cone defined in \eqref{eq:critical}.
%recall that $\nabla^2_{xx}L$ is locally lipshdtz and, by Lemma\,\ref{lem:sigma}, also so is $\Omega(\xast,Y)$ with respect to $Y\in\mathbb{S}^m_+$.
For any $d\in C(\xast)$, we have
\begin{align}
&d^{\top}\left(\nabla^2_{xx}L(\xast,\Yast,\zast)+\Omega(\xast,\Yast)\right)d\notag \\
=&d^{\top}\left(\nabla^2_{xx}L(\wa)+\left(\frac{\partial^2G(\xast)}{\partial x_i\partial x_j}\bullet (\Yast-\Ya)\right)_{1\le i,j\le n}+\sum_{i=1}^s\nabla^2h_i(\xast)(\zast_i-\za_i)\right)d\notag\\
&\hspace{2em}+d^{\top}\Omega(\xast,\Ya)d+
2
\Tr\left(
\left(\YUU_{\ast}-\YUU_{\rm a}\right)\DeltahGNU(x^{\ast};d)
(\Gast^{\NN})^{-1}\DeltahGUN(x^{\ast};d)
\right)\notag\\
\ge&d^{\top}\left(\nabla^2_{xx}L(\wa)+\Omega(\xast,\Ya)\right)d
-\|\zast-\za\|\|d\|^2\sum_{i=1}^s\|\nabla^2h_i(\xast)\|_{\rm F}
-\frac{n(n+1)}{2}\|\Yast-\Ya\|_{\rm F}\|d\|^2
\max_{1\le i,j\le n}
\left\|\frac{\partial^2G(\xast)}{\partial x_i\partial x_j}\right\|_{\rm F}\notag\\
&\hspace{2em}-n\|d\|^2\|\YUU_{\ast}-\YUU_{\rm a}\|_{\rm F}\|(\Gast^{\NN})^{-1}\|_{\rm F}
\max_{1\le i\le n}\|\PU\calG_i(\xast)\PN\|^2_{\rm F}\notag
\\
\ge &d^{\top}\left(\nabla^2_{xx}L(\wa)+\Omega(\xast,\Ya)\right)d-\rho \widetilde{K}K_2\|d\|^2\notag\\
\ge&(\kappa-\rho \widetilde{K}K_2) \|d\|^2, \notag
\end{align}
where the first equality follows from
\cref{lem:sigma} and \eqref{eq:scGastYast}, and
the third inequality follows from \eqref{eq:231216-5} and
$\|\YUU_{\ast}-\YUU_{\rm a}\|_{\rm F}\le \|\Yast-\Ya\|_{\rm F}$ and the last inequality from \tred{{\SSOSC}}\,\eqref{eq:ssoc} with $(Y,z)=(\Ya,\za)$.
Thus, the last inequality implies
\begin{equation}
d^{\top}\left(\nabla^2_{xx}L(\xast,\Yast,\zast)+\Omega(\xast,\Yast)\right)d\ge \frac{\kappa\|d\|^2}{2}\ \mbox{when }0<\rho\le \frac{\kappa}{2\widetilde{K}K_2}.
\label{eq:231202-2}
\end{equation}}
\tred{Hereafter, we set $0<\rho\le \min\left(\frac{\lambda_{\rm min}(\Gast+\Ya)}{2K_2},\frac{\kappa}{2\widetilde{K}K_2},\brho_2\right)$ so that
\eqref{eq:231202-3} and \eqref{eq:231202-2} hold.
Notice that this choice of $\rho$ is independent from the sequence $\{\wl\}$.}

To prove the desired claim, we derive a contradiction by assuming to the contrary, that is, there exists infinite sequences\footnote{
To abuse notation, we use $\well$ and $\mul$ again to denote a sequence.
}
{\cred
$$
\{\mul\}\subseteq \R_{++},\
\{\well:=(\xell,\Yell,\zell)\}\subseteq  \W_{++},\
\{\vell\}\subseteq \R^n
$$}
such that \tred{$\lim_{\ell\to \infty}\mul=0$} and, for each $\ell$, it holds that
$\|\vell\|=1$, $\nabla h(\xell)^{\top}\vell=0$, \eqref{eqn:0524-0}, \eqref{eqn:0524-1}, and \eqref{eqn:0524-2} are fulfilled with $(\mu,w):=(\mul,\well)$, and
\begin{equation}
\calHl:=(\vell)^{\top}\nabla^2_{xx}L(\well)\vell + \Delta G(\xell;\vell)\bullet \calL_{\Gell}^{-1}\calL_{\Yell}\left(\Delta G(\xell;\vell)\right)
\tred{<\frac{\kappa}{2}}.\label{eq:calHl}
\end{equation}

%From \eqref{eq:calHl}, it follows that
\tred{
%Note that the above argument also holds true for this $\{\well\}$
By calculation, we have}
\begin{align}
% \begin{autobreak}
   \calHl&=(\vell)^{\top}\nabla^2_{xx}L(\well)\vell+\DeltaG(\xell;\vell)\bullet \calL_{\Gell}^{-1}\calL_{\mul \Gell^{-1}}\left(\DeltaG(\xell;\vell)\right)
   +\DeltaG(\xell;\vell)\bullet \calL_{\Gell}^{-1}\calL_{\Yell-\mul \Gell^{-1}}\left(\DeltaG(\xell;\vell)\right)
    \notag
    \\
    &=(\vell)^{\top}\nabla^2_{xx}L(\well)\vell+\mul\DeltaG(\xell;\vell)\bullet \Gell^{-1}\DeltaG(\xell;\vell)\Gell^{-1}
   +\DeltaG(\xell;\vell)\bullet \calL_{\Gell}^{-1}\calL_{\Yell-\mul \Gell^{-1}}\left(\DeltaG(\xell;\vell)\right)
    \notag\\
   &= (\vell)^{\top}\nabla^2_{xx}L(\well)\vell+\mul\left\|\Gell^{-\fr}\DeltaG(\xell;\vell) \Gell^{-\fr}\right\|_{\rm F}^2
   +\DeltaG(\xell;\vell)\bullet \calL_{\Gell}^{-1}\calL_{\Yell-\mul \Gell^{-1}}\left(\DeltaG(\xell;\vell)\right).\label{eq:hessipsimudd-4}
 % \end{autobreak}
  \end{align}
As will be verified later on, we actually have the following relationships:
\begin{align}
&\lim_{\ell\to \infty}\DeltaG(\xell;\vell)\bullet \calL_{\Gell}^{-1}\calL_{\Yell-\mul \Gell^{-1}}\left(\DeltaG(\xell;\vell)\right)=0,\label{eq:0528-3-0}\\
&\liminf_{\ell\to\infty}\left((\vell)^{\top}\nabla^2_{xx}L(\well)\vell+\mul\left\|\Gell^{-\fr}\DeltaG(\xell;\vell) \Gell^{-\fr}\right\|_{\rm F}^2\right)\tred{\ge \frac{\kappa}{2}}.\label{eq:0528-3}
\end{align}
From these results 
and \eqref{eq:hessipsimudd-4}, $\liminf_{\ell\to\infty}\calHl\ge \tred{\frac{\kappa}{2}}$ holds.
However, this contradicts hypothesis \eqref{eq:calHl}. Therefore, we have reached the first assertion.
$\hfill\Box$
\subsubsection*{Proofs of \eqref{eq:0528-3-0} and \eqref{eq:0528-3}}
For making the above proof \tred{complete}, it remains to prove \eqref{eq:0528-3-0} and \eqref{eq:0528-3}.
We also suppose the same assumptions as those made for contradiction at the beginning of the above proof. In particular, we will use the same notations and symbols, such as $\{\well\}$, $\{v^{\ell}\}$, and $\{\mu_{\ell}\}$.

Before starting the proofs of \eqref{eq:0528-3-0} and \eqref{eq:0528-3}, we shall give some preliminary results.
\tred{First, note that $\{\well\}$ is bounded as described at the beginning of this section.}
Let $\wrho$ denote an accumulation point of $\{\well\}$.
\tred{Then, notice that the $x$-component of $\wrho$ is the KKT point $\xast$ and denote the $(Y,z)$-component of $\wrho$ by $(\Yrho,\zrho)$.}
Moreover, let $\vast$ be an accumulation point of $\{\vell\}$.  
Choose an orthogonal matrix $\Pell$ for each $\ell$ so that $\Gell$ is diagonalized with $\Pl$ and the eigenvalues of the resultant diagonal matrix is aligned in the ascending order.
By re-choosing $\Past$ and taking 
a subsequence of $\{(x^{\ell},\vell,\Pl)\}$ if necessary, we can suppose, w.l.o.g\footnote{Recall that 
 $\Past$ was an arbitrary orthogonal matrix such that \eqref{eq:past} holds. Even if $\Past$ is reset as the limit of $\{\Pell\}$ here, it satisfies \eqref{eq:past} again, and thus never affects the theoretical results established so far.},  
\begin{equation}
\lim_{\ell\to \infty}(\well,\vell,\Pell)=(\wrho,\vast,P_{\ast}).\label{pelllimit}
\end{equation}
Note that, \tred{as $\nabla h(\xell)^{\top}\vell=0,\|\vell\|=1$ for each $\ell$,}
it follows that 
\begin{equation}
\|\vast\|=1,\ \nabla h(\xast)^{\top}\vast=0.\label{vast}
\end{equation}

Next, so as to match $\Past=[\PU,\PN]$, we partition $\Pell$ as  
$$\Pell=[\PUell,\PNell],$$ which along with \eqref{pelllimit} implies $\lim_{\ell\to \infty}(\PUell,\PNell)=(\PU,\PN)$.
Let the resultant diagonal matrix obtained from $\Gell\in \mathbb{S}^m_{++}$ using $\Pell$
be $D_{\Gell}$, and also let
$D_{\Gell}^0$ and $D_{\Gell}^{++}$ be the block diagonal matrices 
of $D_{\Gell}$ that converge to the $(m-\rast)\times(m-\rast)$ zero matrix and the positive diagonal matrix $\Gast^{\NN}$, respectively. 
Moreover, we often write simply
$$\wtGell:=\Pl\Gell\Pl^{\top}$$ 
for each $\ell$.
In summary, it holds that 
%Finally, we define $D_{\Yell}^0$ and $D_{\Yell}^{++}$ in the same manner. Therefore, 
$$
\wtGell=\Pell\Gell\Pell^{\top}=\begin{bmatrix}
D_{\Gell}^0&O\\
O&D_{\Gell}^{++}
\end{bmatrix},\ \ \lim_{\ell\to\infty}(D_{\Gell}^0,D_{\Gell}^{++})=(O,\Gast^{\NN}).
$$
Accordingly, we denote
$$
\wtYell:=\Pell \Yell\Pell^{\top},\ \DeltawtGell :=\Pell \DeltaGell \Pell^{\top}
$$
with $\DeltaGell:=\DeltaG(\xell;\vell)$.
%Let $\wast=(\xast,\Yast,\zast)$ be an arbitrary accumulation point of $\{\well\}$ and, without loss of generarity, assume that $\well$ converges to $\wast$ by taking a subsequence further if necessary.
%as a symbolic rule, we often writ $\wt{Z}:=\Px Z\Px^{\top}$ for $Z\in \mathbb{S}^m$ for the sake of simplicity. 
Furthermore, so as to match the partition pattern of 
 $\left[\begin{smallmatrix}
 O&O\\
 O&\Gast^{\NN}
 \end{smallmatrix}\right]$, partition a given matrix $Z\in \mathbb{S}^m$ as 
$$Z=\left[\begin{smallmatrix}
 \ZAA&\ZAB\\
 \ZAB^{\top}&\ZBB
 \end{smallmatrix}\right],\ \ZAA\in \mathbb{S}^{m-\rast},\ \ZAB\in \R^{(m-\rast)\times \rast},\ \ZBB\in \mathbb{S}^{\rast}.
$$ 
Now, we start proving \eqref{eq:0528-3-0} and \eqref{eq:0528-3}.
\subsubsection*{Proof of \eqref{eq:0528-3-0}}
\tred{First, recall that $\mathcal{L}_XY=XY+YX$ for $X,Y\in \mathbb{S}^m$.
If $X\in\mathbb{S}^m_{++}$, the linear operator $\mathcal{L}_X$ is invertible, namely, $\mathcal{L}_X^{-1}$ exists.
%Hence, if $X\in \mathbb{S}^m_{++}$,.
}
\tred{Next,} note that, given $W\in \mathbb{S}^m$, a solution $Z\in \mathbb{S}^m$ to $\calL_{\wtGell}Z=W$
 satisfies 
 \begin{align}
  &\VAA = \calL_{\DGxzero}^{-1}\WAA,\ \VBB = \calL_{\DGxpp}^{-1}\WBB,\label{al:0527-1}\\
   &\VAB(i,j)= \frac{1}{\DGxzero(i,i) + \DGxpp(j,j)}\WAB(i,j)\ \ \ (1\le i\le m-\rast, 1\le j \le \rast),\label{al:0527-2} 
 \end{align}
 which are verified by representing $\calL_{\wtGell}Z=W$ as 
 \begin{align}
 %&\BlockDiagGx \blockV+ \blockV \BlockDiagGx= \blockW  \\
  \begin{bmatrix}
      \DGxzero\VAA+\VAA\DGxzero - \WAA&\DGxzero \VAB+\VAB\DGxpp  - \WAB \\
      \DGxpp\VAB^{\top}+\VAB^{\top}\DGxzero- \WAB^{\top}& \DGxpp \VBB+ \VBB\DGxpp - \WBB
    \end{bmatrix}=O.\notag
    \end{align}
\tred{We have that
$\|\Yell-\mul \Gell^{-1}\|_{\rm F}\le \gamma_1\mul$ from \eqref{eqn:0524-1}} and $\{\Pell\}$ is bounded since $\Pell$ is an orthogonal matrix.
These facts yield
\begin{eqnarray}
 &\wtYell-\mul \wtGell^{-1}=\O(\mul).\label{eq:0527-4}
% &\DeltaGell \bullet \calL_{\Gell}^{-1}\calL_{\Yell-\mul \Gell^{-1}} \DeltaGell = \DeltawtGell \bullet \calL_{\wtGell}^{-1}\calL_{\wtYell-\mul \wtGell^{-1}}\DeltawtGell.
% \label{eq:0527-5}
 \end{eqnarray}
In view of \eqref{al:0527-1} and \eqref{al:0527-2} with $W:=\calL_{\wtYell-\mul \wtGell^{-1}} \DeltawtGell$, the solution $Z$ satisfies  
 \begin{align}
 \VAA = \O(\|\DeltawtGell\|_{\rm F}),\ \VBB = \O(\mul\|\DeltawtGell\|_{\rm F}),\ \VAB(i,j)=\O(\mul\|\DeltawtGell\|_{\rm F})\ \ \ (1\le i\le m-\rast, 1\le j \le \rast), \label{al:0528-1}
 \end{align} 
 where the first equation in \eqref{al:0528-1}
 is derived from the fact  
 $$Z=\mul \calL_{\wtGell}^{-1}\calL_{\frac{1}{\mul}{(\wtYell-\mul \wtGell^{-1})}} \DeltawtGell=\O(\|\DeltawtGell\|_{\rm F}),$$
 which is ensured by \eqref{eq:0527-4} and  the boundedness of $\mul \calL_{\wtGell}^{-1}$.
 (For the proof of the boundedness of $\mul \calL_{\wtGell}^{-1}$, see the footnote\,\footnote{
 Note that, for any $X\in \mathbb{S}^m$ having $m$ eigenvalues $\alpha_1\le \alpha_2\le\cdots,\le \alpha_m$,
the linear operator $\mathcal{L}_X$ is symmetric and has $m(m+1)/2$ eigenvalues
$
\alpha_1,\alpha_2,\ldots,\alpha_m,\{(\alpha_i+\alpha_j)/2\}_{i\neq j}.
$
\tred{Letting $(0<)\lambda_1^{(\ell)}\le \lambda_2^{(\ell)}\le \cdots \lambda_m^{(\ell)}$ be the eigenvalues of $\widetilde{\Gell}\in \mathbb{S}^m_{++}$ for each $\ell$,
the eigenvalues of the symmetric linear operator $\mul\mathcal{L}_{\widetilde{\Gell}}^{-1}$ are
$\mul/\lambda_1^{(\ell)},\mul/\lambda_2^{(\ell)},\cdots \mul/\lambda_m^{(\ell)}$, and $\left\{(\mul/\lambda^{(\ell)}_i+\mul/\lambda^{(\ell)}_j)/2\right\}_{i\neq j}$.
Since 
$\Set{\mu G(x)^{-1}| x\in \cl\mathcal{P}_{\brho_1}(\mu), \mu\in (0,\bmu_1]}$ is bounded from item~\ref{prop:0513-1-item1-1} of \cref{prop:0513-1}, 
so are $\{\mul \Gell^{-1}\}$ and $\{\mul \widetilde{\Gell}^{-1}\}$. Hence, $\{\mul/\lambda_i^{(\ell)}\}$ is also bounded for each $i$, which
together with $\|\mul\mathcal{L}_{\widetilde{\Gell}}^{-1}\|_2:=\displaystyle{\max_{Z\in \mathbb{S}^m:\|Z\|_{\rm F}=1}}\|\mul\mathcal{L}_{\widetilde{\Gell}}^{-1}Z\|_{\rm F}\le \mul/\lambda^{(\ell)}_1$
yields the boundedness of $\{\mul\mathcal{L}_{\widetilde{\Gell}}^{-1}\}$.
}}.)
Moreover, the second and third equations 
 in \eqref{al:0528-1}
 are implied by \eqref{al:0527-2} and the right equation in \eqref{al:0527-1}.
 Using \eqref{al:0528-1} again \tred{and noting that $\{\widetilde{\DeltaGell}\}$ is bounded}, we obtain
 \begin{align}
 \DeltaGell \bullet \calL_{\Gell}^{-1}\calL_{\Yell-\mul \Gell^{-1}} \DeltaGell 
 =&\DeltawtGell \bullet \calL_{\wtGell}^{-1}\calL_{\wtYell-\mul \wtGell^{-1}} \DeltawtGell\notag\\
 =&
 \DeltawtGell \bullet Z\notag\\
 =&{\rm Tr}\left(
 (\DeltawtGell)_{11}Z_{11}+ 2(\DeltawtGell)_{12}Z_{12}^{\top}
+ (\DeltawtGell)_{22}Z_{22} 
 \right)\notag\\
 =&O\left(\|(\DeltawtGell)_{11}\|_{\rm F}+\mul\right).\label{al:1221-1}
 \end{align}
\tred{In order to prove the desired equation\,\eqref{eq:0528-3-0}, it suffice to verify}  
\begin{equation}
\DeltaGUU(\xast;\vast)=O,\label{eq:1221-2}
\end{equation}
%where $C(\xast)$ is the critical cone and expressed as in \eqref{eq:critical}.
\tred{because
\eqref{eq:0528-3-0} is verified by using \eqref{al:1221-1} together with the fact that
$\lim_{\ell\to\infty}(\DeltawtGell)_{11}=\DeltaGUU(\xast;\vast)=O$ following from \eqref{eq:1221-2}.} 
To this end,
we evaluate $\mul\left\|\Gell^{-\fr}\DeltaGell\Gell^{-\fr}\right\|_{\rm F}^2$ in \eqref{eq:hessipsimudd-4} as follows:
\begin{align*}
 \mul\left\|\Gell^{-\fr}\DeltaGell\Gell^{-\fr}\right\|_{\rm F}^2=&\mul {\rm Tr}\left(\Gell^{-1}\DeltaGell\Gell^{-1}\DeltaGell\right)\notag\\
 %                                                                         &={\rm Tr}\left(D_{G_k}^{-1}\Delta\widetilde{G}(x^k;v^k)D_{Y_k}\Delta\widetilde{G}(x^k;v^k)\right)\notag\\
                                                                          =&\mul\Tr\left(\left\{
                                                                          \begin{bmatrix}
                                                                          (D_{\Gell}^0)^{-1}&O\\
                                                                          O&(D_{\Gell}^{++})^{-1}
                                                                          \end{bmatrix}
                                                                          \begin{bmatrix}
                                                                          (\DeltawtGell)_{11}&(\DeltawtGell)_{12}\\
                                                                          (\DeltawtGell)_{21}&(\DeltawtGell)_{22}
                                                                          \end{bmatrix}\right\}^2
                                                                          \right)\notag \\
                                                                            =&\mul\Tr\left(
                                                                          \begin{bmatrix}
                                                                          (D_{\Gell}^0)^{-1}(\DeltawtGell)_{11}&&(D_{\Gell}^0)^{-1}(\DeltawtGell)_{12}\\
                                                                          (D_{\Gell}^{++})^{-1}(\DeltawtGell)_{21}&&(D_{\Gell}^{++})^{-1}(\DeltawtGell)_{22}
                                                                          \end{bmatrix}^2                                                                                                      
                                                                          \right)\notag\\
                                                                          =&\underbrace{\mul\Tr\left((D_{\Gell}^0)^{-1}(\DeltawtGell)_{11}(D_{\Gell}^{0}
                                                                          )^{-1}(\DeltawtGell)_{11}\right)}_{(\aell)}
                                                                          +\underbrace{2\mul\Tr\left((D_{\Gell}^0)^{-1}(\DeltawtGell)_{12}(D_{\Gell}^{++})^{-1}(\DeltawtGell)_{12}^{\top}
                                                                          \right)}_{(\bell)}\notag\\
                                                                          &\hspace{0.5em}
                                                                          +\underbrace{\mul\Tr\left((D_{\Gell}^{++})^{-1}(\DeltawtGell)_{22}(D_{\Gell}^{++})^{-1}(\DeltawtGell)_{22}\right)                                                                    }_{(\cell)}, \notag 	
 \end{align*}
 and herein we obtain 
 \begin{equation}
   \lim_{\ell\to\infty}(\cell)=0\label{eq:dell}
 \end{equation}
because $\lim_{\ell\to\infty}\mul=0$ and the matrices in
the trace-part of $(\cell)$ are convergent. 
It holds that $\lim_{\ell\to\infty}\mul\Gell^{-1}=\Yrho$ 
since $\mul \Gell^{-1}$ and $\Yell$ accumulate at identical points due to \eqref{eqn:0524-1}, thus
which together with \eqref{eq:231202-3} yields
\begin{equation}
\lim_{\ell\to\infty}\mul (D_{\Gell}^{0})^{-1}=\PU^{\top}{\Yrho}\PU=(\wt{Y}_{\ast})_{11}\in \mathbb{S}^m_{++}, \label{eq:231230-1}
\end{equation}
which yields
\begin{align}
\lim_{\ell\to\infty}(\bell)=2\Tr\left((\Gast^{\NN})^{-1}\DeltaGNU(x^{\ast};v^{\ast})(\wt{Y}_{\ast})_{11}\DeltaGUN(x^{\ast};v^{\ast})\right)=(v^{\ast})^{\top}\Omega(x^{\ast},\Yrho)v^{\ast},
%=&2\Tr\left((\Gast^{\NN})^{-1}\DeltaGNU(x^{\ast};v^{\ast})(\wt{Y}_{\brho})_{11}-D_{\Ya}^{++})\DeltaGUN(x^{\ast};v^{\ast})\right)
%+2\Tr\left((\Gast^{\NN})^{-1}\DeltaGNU(x^{\ast};v^{\ast})D_{\Ya}^{++}\DeltaGUN(x^{\ast};v^{\ast})\right)\notag \\
%=&2\Tr\left((\Gast^{\NN})^{-1}\DeltaGNU(x^{\ast};v^{\ast})\left((\wt{Y}_{\brho})_{11}-D_{\Ya}^{++}\right)\DeltaGUN(x^{\ast};v^{\ast})\right)+(v^{\ast})^{\top}\Omega(x^{\ast},\Ya)v^{\ast},
\label{al:0526-2}
\end{align}
where the last equality follows from \cref{lem:sigma}. 
Therefore, in view of \eqref{eq:calHl}, the last equality in \eqref{eq:hessipsimudd-4}, and \eqref{al:1221-1}, by noting that $(\vell)^{\top}\nabla^2_{xx}L(\well)\vell$ is bounded, we find that $\mul\left\|\Gell^{-\fr}\DeltaGell\Gell^{-\fr}\right\|_{\rm F}^2$ is bounded,
which together with \eqref{eq:dell} and \eqref{al:0526-2} implies that $\{(\aell)\}$ is also bounded. 
\tred{From this fact together with \eqref{eq:231230-1}, we ensure
$\DeltaGUU(\xast;\vast)=\lim_{\ell\to\infty}(\DeltawtGell)_{11}=O$}
thus \eqref{eq:1221-2}. 
%\tred{because $\mul (D^0_{\Gell})^{-1}$ converges to a positive definite matrix from \eqref{eq:231202-3}.}
The proof of \eqref{eq:0528-3-0} is complete.
%$\hfill\Box$
%which combined with \eqref{eq:0526-2} and the convergence of $\mul (D^0_{\Gell})^{-1}$ shown above 
%$\Delta \GUU(\xast;\vast)=O$, leading to $v^{\ast}\in C(x^{\ast})$.

%% Since $\left\{\mul \Gell^{-1}\right\}$ is bounded from \cref{prop:0513-1}, 
%% \cref{lem:0528-1} implies that the sequence of operators 
%% $\left\{\mul\calL_{\Gell}^{-1}\right\}$ is also bounded. 
%% Moreover, by assumption\,\eqref{eqn:0524-1}, the sequence   
%% $\{\frac{1}{\mul}\calL_{\Yell-\mul \Gell^{-1}}\}$ is bounded, too.  Using these boundedness results and noting that 
%% $\{\Pl\DeltaG(\xell;\vell)\Pl^{\top}\}$ is convergent, 
%% we obtain 
%% \begin{align}
%%   \DeltaG(\xell;\vell)\bullet \calL_{\Gell}^{-1}\calL_{\Yell-\mul \Gell^{-1}}\left(\DeltaG(\xell;\vell)\right)=
%%   %% =&O\left(\|\Pl\DeltaG(\xell;\vell)\Pl^{\top}\|_{\rm F}\left(\|\DeltaGUU(\xell;\vell)\|_{\rm F}+\mul\left(\|\DeltaGUN(\xell;\vell)\|_{\rm F}+\|\DeltaGNN(\xell;\vell)\|_{\rm F}\right)\right)\right)\label{eq:deltag-0} \\	
%%   %% =&
%%   \O(1).\label{eq:deltag}
%% \end{align}

\subsubsection*{Proof of \eqref{eq:0528-3}}
From \eqref{eq:1221-2} and \eqref{vast}, $\vast\in C(\xast)$ holds, where $C(\xast)$ is the critical cone and expressed as in \eqref{eq:critical}.
%% Note that $(\Yrho,\zrho)\in \Lambda(\xast)$. Hence, by the SSOSC
% (see \eqref{eq:ssoc} again for the definition), 
%% \begin{equation}
%% (v^{\ast})^{\top}\left(\nabla^2_{xx}L(\wrho)+\Omega(x^{\ast},\Yrho)\right)v^{\ast}>0.\label{eq:1222-2}
%% \end{equation}
%Additioanlly, we have 
% \begin{align}
% \eqref{al:0526-2}
% \ge & -2\|\DeltaGNU(x^{\ast};v^{\ast})\|_{\rm F}^2\|
% \YUU_{\brho}-D_{\Ya}^{++}
% \|_{\rm F}\|(\Gast^{\NN})^{-1}\|_{\rm F}
% +(v^{\ast})^{\top}\Omega(x^{\ast},\Ya)v^{\ast}\notag \\
% \ge&	-2\|\DeltaGNU(x^{\ast};v^{\ast})\|_{\rm F}^2\|\PU\|^2_{\rm F}
% \|\Yrho-\Ya\|_{\rm F}\|(\Gast^{\NN})^{-1}\|_{\rm F}
% +(v^{\ast})^{\top}\Omega(x^{\ast},\Ya)v^{\ast}. \label{al:0526-4}
% \end{align}
Let $$
\Xi_{\ell}:=(\vell)^{\top}\nabla^2_{xx}L(\well)\vell+\mul\left\|\Gell^{-\fr}\DeltaG(\xell;\vell) \Gell^{-\fr}\right\|_{\rm F}^2$$
for each $\ell$.
Since $(\aell)\ge 0$ for each $\ell$ and $\{(\aell)\}$ is bounded as shown in the proof of \eqref{eq:0528-3-0},
any accumulation point of $\{(\aell)\}$ is nonnegative.
Combining this fact with \eqref{eq:dell}, \eqref{al:0526-2}, and \tred{\eqref{eq:231202-2} with $d:=\vast$} yields
$$
\liminf_{\ell\to\infty}\Xi_{\ell}=(v^{\ast})^{\top}\left(\nabla^2_{xx}L(\wrho)+\Omega(x^{\ast},\Yrho)\right)v^{\ast}
+\liminf_{\ell\to\infty}(\aell)\tred{\ge \frac{\kappa}{2}\|v^{\ast}\|^2=\frac{\kappa}{2}},
$$
which implies \eqref{eq:0528-3}. The proof is complete.

% Combined with \eqref{al:0526-2}, \eqref{al:0526-3}, \eqref{al:0526-4}, and 
% condition\,\eqref{al:0526-1} yields 
% \begin{align}
% (v^{\ast})^{\top}\nabla^2_{xx}L(\wrho)v^{\ast}+\lim_{\ell\to\infty}\left((\cell)+(\bell)\right)>0,
% \end{align}
% by which along with $(\aell)\ge 0$ for all $\ell$, \eqref{eq:0528-3}, and \eqref{eq:dell}, we gain that $\calHl$ is bounded from below with a positive number for any $\ell$ sufficiently large.  However, this contradicts the hypothesis $\calHl\le 0$ for each $\ell$.

$\hfill\Box$

\subsection{Proof of \cref{lem:0314-1}}\label{appA-2}
{\cred
First, we show that 
\begin{equation}
\dist(\xmudel,\calM)=\O(\|h(\xmudel)\|)\label{eq:distxM}
\end{equation}
for $\xmudel \in \hball$.
Note that \eqref{eq:231216-3} holds as for $\hball$.}
%% Recall that $\nabla h(\xast)$ is of full column rank.
%% First, we show that 
%% \begin{equation}
%% \dist(x,\calM)=\O(\|h(x)\|)\label{eq:distxM}
%% \end{equation}
%% for any $x\in \R^n$ locally around $\xast$.
%% Let $C\subseteq \R^n$ be a closed ball centered at $\xast$, such that 
%% $\nabla h(x)$ is of full column rank for all $x\in C$, namely, 
%% \begin{equation}
%% {\rm rank}\nabla h(x) = {\rm rank}\nabla h(\xast)=s,\ \forall x\in C.\label{eq:rankC}
%% \end{equation}
To prove \eqref{eq:distxM}, assume to the contrary: there exists a sequence 
\tred{$\{\mul\}$ converging to $0$ such that
\begin{equation}
\dist(\xmudell,\calM)\neq 0,\ \forall \ell,\ \ \lim_{\ell\to\infty}\frac{\|h(\xmudell)\|}{\dist(\xmudell,\calM)}=0.\label{eq:240122-1}
\end{equation}
}
For each $l$, let $\yell\in \argmin_{y\in \calM\cap \tred{\hball}}\|\xmudell-y\|$. Then, we have that, for any $\ell$ large enough, 
\begin{equation}
\dist(\xmudell,\calM)=\|\xmudell-\yell\| \label{eq:240122-2}
\end{equation}
\tred{since $\lim_{\ell\to \infty}\xmudell=\xast\in {\rm int}\,\hball$, $\lim_{\ell\to\infty}\xmudell-\yell=0$, and thus $\yell\in {\rm int}\,\hball$. In what follows, we consider only $\ell$ large enough and assume $\yell\in {\rm int}\,\hball$.}
%Notice $\yell \to \tred{\xast}$ as $\xmudell$ approaches $\xast$ and thus $\tred{\xmudell}-\yell \to 0$. 
%Since $\yell$ solves $\min_{y\in \calM\cap C}\|\xmudell-y\|^2$ by definition 
%and $\yell$ stays in the interior of $C$ for any $l$ large enough, 
\tred{Since $\yell$ solves $\min_{y\in \calM} \frac{1}{2}\|\xmudel-y\|^2$ and the LICQ holds at $\yell$ from \eqref{eq:231216-3},
the KKT conditions holds, namely} there exists $\etal$ such that $$\xmudell-\yell=\nabla h(\yell)\etal.$$ 
\tred{From \eqref{eq:240122-1} and \eqref{eq:240122-2}, $\xmudell-\yell=\nabla h(\yell)\etal\neq 0$.}
Letting $\bar{v}$ be an accumulation point of $\left\{(\xmudell-\yell)/\|\xmudell-\yell\|\right\}$, we may assume that
\begin{equation}
\frac{(\xmudell-\yell)}{\|\xmudell-\yell\|}=\frac{\nabla h(\yell)\etal}{\|\nabla h(\yell)\etal\|}\to \bar{v}\ \ (l\to \infty)
\label{eq:linf}
\end{equation}
without loss of generality. It follows that $\|\bar{v}\|=1$.
As
\begin{align*}
\frac{\|h(\xmudell)\|}{\dist(\xmudell,\calM)}=\frac{\|h(\xmudell)-h(\yell)\|}{\|\xmudell-\yell\|}
                           =\frac{\|\nabla h(\yell)^{\top}(\xmudell-\yell)+\O(\|\xmudell-\yell\|^2)\|}{\|\xmudell-\yell\|},
\end{align*}
driving $\ell\to \infty$ herein yields  
\begin{equation}
\nabla h(\xast)^{\top}\bar{v}=0.
\label{eq:linf-2}
\end{equation}
%Note $\yell\to \xast$ holds.        
\tred{Since $\nabla h(\xast)$ is of full column rank, $\nabla h$ is continuous, and $\lim_{\ell\to \infty}\yell=\xast$,}
%because of $\xast\in C$ and \eqref{eq:rankC},
\eqref{eq:linf-2} actually implies that there exists $\{\vell\}$ such that
it converges to $\bar{v}$ and, for each $\ell$, 
\begin{equation}
\nabla h(\yell)^{\top}\vell=0.
\label{eq:linf-3}
\end{equation}
%% Thus, 
%% \begin{equation}
%% \frac{(v^l)^{\top}\nabla h(\yell)\etal}{\|\nabla h(\yell)\etal\|}=0\label{eq:linf-3}
%% \end{equation}
%% for each $l$.
Meanwhile, with \eqref{eq:linf} and $\lim_{l\to\infty}\vell=\bar{v}$, we gain 
$$
\lim_{\ell\to\infty}\frac{(\vell)^{\top}\nabla h(\yell)\etal}{\|\nabla h(\yell)\etal\|}= \|\bar{v}\|^2,
$$
which combined with \eqref{eq:linf-3} gives $\|\bar{v}\|^2=0$. However, this contradicts $\|\bar{v}\|=1$.  We thus obtain \eqref{eq:distxM}.

We next prove the desired relation\,\eqref{eq:0608-2}. 
Note that $h(\xmudel)=h(\xast+\mu\xiast)=\O(\mu^2)$ holds by Taylor's expansion together with
$h(\xast)=0=\nabla h(\xast)^{\top}\xiast$, where the last equality follows from \cref{prop:xiast}.
With this fact along with \eqref{eq:distxM}, \eqref{eq:0608-2} is ensured. The proof is complete. 
$\hfill\Box$

{\cred
\subsection{Proof of \cref{prop:231219-1}}\label{app:prop:231219-1}
In order to prove \cref{prop:231219-1}, we start by showing the following lemma:
\begin{lemma}\label{lem:231213-5}
  For $v\in V$ and $y\in \R^{n-s}$, let
  %\begin{align}
  $(x,Y,d):=(\impfun(\uv),\mu G(\impfun(\uv))^{-1},\nabla\impfun(\uv)^{\top}y)$
  %\end{align*}
  and suppose $G(x)\in \mathbb{S}^m_{++}$. Moreover, let $z\in \R^s$ and $\mu>0$.
Then, we have
  \begin{equation}
y^{\top}\nabla^2\Psi_{\mu}(\uv)y=d^{\top}\nabla_{xx}^2L(x,Y,z)d+
  \sum_{j=1}^n\left(\frac{\partial L(x,Y,z)}{\partial x_j}
  y^{\top}\nabla^2 \impfun_j(\uv)y\right)+\Delta G(x;d)\bullet \calL_{G(x)}^{-1}\calL_{Y}\left(\Delta G(x;d)\right).\label{al:231213-1}
  \end{equation}
\end{lemma}
\proof{\rm Proof.}
By calculation, as regards the function $\psi_{\mu}$ defined in \eqref{eq:231213-4}, 
we have
\begin{align}
  %      \frac{\partial \psi_{\mu}(x)}{\partial x_j}=&\frac{\partial f(x)}{\partial x_j}-\mu \Gj(x)\bullet G(x)^{-1}\ (j=1,2,\ldots,n),\label{eqn:231211-1}\\
      \nabla \psi_{\mu}(x)=&\nabla f(x)-\mu \mathcal{J}G(x)^{\ast}G(x)^{-1}=\nabla_xL(x,\mu G(x)^{-1},z)-\nabla h(x)z,\label{eqn:231211-1}\\
      d^{\top}\nabla^2\psi_{\mu}(x)d=&d^{\top}\nabla^2f(x)d-\mu
      \left(\sum_{i=1}^n\sum_{j=1}^n
      d_id_j\frac{\partial^2 G(x)}{\partial x_i\partial x_j}\right)\bullet G(x)^{-1}+\mu\|G(x)^{-\fr}\Delta G(x;d)G(x)^{-\fr}\|_{\rm F}^2\notag\\
      =&d^{\top}\nabla_{xx}^2L(x,\mu G(x)^{-1},z)d-d^{\top}
\left(\sum_{i=1}^sz_i\nabla^2h_i(x)\right)d+\mu\|G(x)^{-\fr}\Delta G(x;d)G(x)^{-\fr}\|_{\rm F}^2.
      \label{eqn:231211-2}
\end{align}
Moreover,  it follows from \eqref{eq:231210-2-2} that  
      \begin{equation}
        \sum_{j=1}^n\left(\sum_{i=1}^sz_i\frac{\partial h_i(\impfun(\uv))}{\partial x_j}\right)\nabla^2 \impfun_j(\uv)+\nabla \impfun(\uv)\left(\sum_{i=1}^sz_i\nabla_{xx}^2h_i(\impfun(\uv))\right)\nabla \impfun(\uv)^{\top}=O. 
        \label{eq:231210-3}
      \end{equation}
      Then, with $(x,Y,d)=\left(\impfun(\uv),\mu G(\impfun(\uv))^{-1},\nabla\impfun(\uv)^{\top}y\right)$, equation\,\eqref{al:231211-2} yields
\begin{align}
  y^{\top}\nabla^2\Psi_{\mu}(\uv)y=&
  d^{\top}\nabla_{xx}^2L(x,Y,z)d+
  \sum_{j=1}^n\frac{\partial \psi_{\mu}(\impfun(\uv))}{\partial x_j}
  y^{\top}\nabla^2 \impfun_j(\uv)y-
  d^{\top}\left(\sum_{i=1}^sz_i\nabla^2h_i(x)\right)d
  %\frac{\partial \psi_{\mu}(x)}{\partial x_j}\right)y^{\top}\nabla^2 \impfun_j(\uv)y
  +\mu\|G(x)^{-\fr}\Delta G(x;d)G(x)^{-\fr}\|_{\rm F}^2\notag\\
=&
  d^{\top}\nabla_{xx}^2L(x,Y,z)d+
  \sum_{j=1}^n\left(\frac{\partial \psi_{\mu}(x)}{\partial x_j}
  +\sum_{i=1}^sz_i\frac{\partial h_i(x)}{\partial x_j}\right)
  y^{\top}\nabla^2 \impfun_j(\uv)y
  %\frac{\partial \psi_{\mu}(x)}{\partial x_j}\right)y^{\top}\nabla^2 \impfun_j(\uv)y
  +\mu\|G(x)^{-\fr}\Delta G(x;d)G(x)^{-\fr}\|_{\rm F}^2\notag\\
=&
  d^{\top}\nabla_{xx}^2L(x,Y,z)d+
  \sum_{j=1}^n\left(\frac{\partial L(x,Y,z)}{\partial x_j}
  y^{\top}\nabla^2 \impfun_j(\uv)y\right)
+\Delta G(x;d)\bullet \calL_{G(x)}^{-1}\calL_{Y}\left(\Delta G(x;d)\right)\notag
\end{align}
where the first equality follows from \eqref{eq:231210-2-2} and \eqref{eqn:231211-2},
the second from \eqref{eq:231210-3}, the third from \eqref{eqn:231211-1} and 
%\begin{align}
$\mu\|G(x)^{-\fr}\Delta G(x;d)G(x)^{-\fr}\|_{\rm F}^2=\mu \Delta G(x;d)\bullet G(x)^{-1}\Delta G(x;d)G(x)^{-1}
=\Delta G(x;d)\bullet \calL_{G(x)}^{-1}\calL_{\mu G(x)^{-1}}\left(\Delta G(x;d)\right)
=\Delta G(x;d)\bullet \calL_{G(x)}^{-1}\calL_{Y}\left(\Delta G(x;d)\right).$
%\end{align}
%Notice the above $(x,Y,z)$ fulfills the assumptions in \cref{prop:231211-1} and \cref{prop:0531-1}.
$\hfill\Box$\vspace{0em}

\subsubsection*{Proof of \cref{prop:231219-1}}
First, let $K_2,M_3,M_4>0$ be the constants defined in \cref{prop:231211-1} and \eqref{eq:231217-1}, and $(\brhofour,\bmufour)\in (0,\brhothree]\times (0,\bmuthree]$
be such that
%\begin{align}
$\max(\brhofour,\bmufour)\le {\kappa M_3^2}/(8nK_2M_4^2M_5)$.
%\label{al:231229-1}
%\end{align}
Choose $(\rho,\mu)\in (0,\brhofour] \times (0,\bmufour]$ and $v\in V\cap \impfun^{-1}(\cl\calPc(\mu))$ arbitrarily. 
Let $x:=\impfun(v)\in\cl\calPc(\mu)$ and recall $\bmufour\le \bmuone$ and $\brhofour\le \brhoone$.
Then, from \eqref{eq:231202-1}, it holds that $G(x)\in \mathbb{S}^m_{++}$
and $(\nabla h(x)\nabla h(x)^{\top})^{-1}$ exists.  %for all $\mu\le \bmufour$.
Thus, we can define $Y:=\mu G(x)^{-1}$ and
$z:=-(\nabla h(x)\nabla h(x)^{\top})^{-1}\nabla h(x)^{\top}(\nabla f(x)-\calJ G(x)^{\ast}Y)$, and these $(x,Y,z)$ fulfills the conditions \eqref{eqn:0524-0}, \eqref{eqn:0524-1}, and \eqref{eqn:0524-2}.
Letting $d:=\nabla\impfun(\uv)^{\top}y$, we obtain
$\nabla h(x)^{\top}d=\sum_{i=1}^s\nabla h_i(x)^{\top}\nabla\impfun(\uv)^{\top}y=0$ from \eqref{eqn:231213-2}.
%Choose arbitrary $y\in \R^{n-s}$ such that $\|y\|=1$ and let $d:=\nabla\impfun(\uv)^{\top}y$.
%Note that $\nabla h(x)^{\top}d=\sum_{i=1}^s\nabla h_i(x)^{\top}\nabla\impfun(\uv)^{\top}y=0$ from \eqref{eqn:231213-2}.
%% Since $V$ is bounded and $\nabla^2\impfun_j$ is continuous on $\cl V$, there exists some $M_5>0$ such that
%% $\|\nabla^2\impfun_j(u)\|_{\rm F}\le M_5$ for any $j$ and $u\in V$
%% Then, it holds that for any $\ell$ sufficiently large,
It holds that 
\begin{align*}
  y^{\top}\nabla^2\Psi_{\mu}(\uv)y
  =&d^{\top}\nabla_{xx}^2L(x,Y,z)d
  +\Delta G(x;d)\bullet \calL_{G(x)}^{-1}\calL_{Y}\Delta G(x;d)
  +\sum_{j=1}^n\left(\frac{\partial L(x,Y,z)}{\partial x_j}d^{\top}\nabla^2 \impfun_j(\uv)d\right)\\
 \ge&d^{\top}\nabla_{xx}^2L(x,Y,z)d
  +\Delta G(x;d)\bullet \calL_{G(x)}^{-1}\calL_{Y}\Delta G(x;d)
  -\sum_{j=1}^n\|\nabla_xL(x,Y,z)\|\|\nabla^2\impfun_j(v)\|_{\rm F}\|d\|^2\notag\\
  \ge& \frac{\kappa}{2}M_3^2-nK_2M_4^2M_5(\mu+\rho)\ge \frac{\kappa M_3^2}{4},\notag 
\end{align*}
where
%$K_2,M_3,M_4>0$ are the constants defined in \cref{prop:231211-1} and \eqref{eq:231217-1},
the first equality follows from \cref{lem:231213-5}, the third inequality from \cref{prop:231211-1}, \cref{prop:0531-1}, and \eqref{eq:231217-1} with $d=\nabla\impfun(\uv)^{\top}y$, and the last inequality from the above-mentioned definitions of $\bmufour$ and $\brhofour$.
%Then, by re-taking $\brho(\ge \rho)$ so that $\brho \le \kappa M_3^2/(8nK_2M_4^2M_5)$,
%Then, , we obtain $y^{\top}\nabla^2\Psi_{\mu}(\uv)y\ge 3\kappa M_3^2/8-nK_2M_4^2M_5\mu$ and 
%Hence, by re-taking $(\mu\le )\bmu$ so that $\bmu\le \kappa M_3^2/(8nK_2M_4^2M_5)$ if necessary, we have
%$y^{\top}\nabla^2\Psi_{\mu}(\uv)y\ge \kappa M_3^2/4$.
The proof is complete.
$\hfill\Box$}

\end{document}